\documentclass[fontsize=10pt]{scrreprt}

\usepackage[square,numbers,nonamebreak]{natbib}
\usepackage{setspace}

\usepackage[english]{babel}

\usepackage{chngcntr}

\usepackage[hang]{footmisc}
 at 11truept
 at 11 truept
\newcommand\qedbox{$\rlap{$\sqcap$}\sqcup$}

\usepackage[hmarginratio=1:1, bottom=2cm, top=2.5cm]{geometry}

\usepackage[automark]{scrlayer-scrpage}
\cohead{}

\let\ceheadL\cehead
\renewcommand\cehead[1]{
\ceheadL{\textnormal{#1}}
}

\cehead{Brescia -- Di Siena -- Ingrosso -- Trombetti}
\lefoot{}
\rofoot{}

\usepackage{graphicx}
\usepackage[usenames,dvipsnames,svgnames,table]{xcolor}
\definecolor{Maroon}{cmyk}{0, 0.87, 0.68, 0.32}
\definecolor{RoyalBlue2}{cmyk}{80,100,0,0.1}

\renewcommand\abstract[1]{
\begin{center}
{\textbf{Abstract}}
\end{center}
{
\linespread{1.1}\fontsize{9pt}{-10pt}\selectfont #1}}

\usepackage[utf8]{inputenc}
\usepackage[sc,osf]{mathpazo}
\usepackage[euler-digits]{eulervm}
\usepackage[T1]{fontenc}
\usepackage{mathtools}

\DeclareSymbolFont{operators}{\encodingdefault}{ppl}{m}{n}
\DeclareMathAlphabet{\mathbf}{\encodingdefault}{ppl}{bx}{n}
\DeclareMathAlphabet{\mathit}{\encodingdefault}{ppl}{m}{it}

\usepackage{amsfonts}
\usepackage{amsmath}
\usepackage{ragged2e}
\usepackage{titlesec}
\usepackage{amssymb}
\usepackage{lipsum}

\renewcommand{\thesection}{\arabic{section}}
 
\titleformat{\section}{\medskip\bigskip\normalfont\Large\bf}{\thesection}{0.5em}{}
\titleformat{\subsection}{\smallskip\bigskip\normalfont\large\bf}{\thesubsection}{0.5em}{}

\usepackage[hidelinks]{hyperref}
\usepackage{amsthm}
\newtheoremstyle{dotless}{}{}{\itshape}{}{\bfseries}{}{1em}{}

\newtheorem*{theo*}{Theorem}
\newtheorem*{rem*}{Remark}
\theoremstyle{dotless}
\newtheorem{theorem}{Theorem}

\newtheorem{lemma}[theorem]{Lemma}
\newtheorem{definition}[theorem]{Definition}
\newtheorem{corollary}[theorem]{Corollary}
\newtheorem{claim}[theorem]{Claim}
\newtheorem{remark}[theorem]{Remark}
\newtheorem{fact}[theorem]{Fact}
\newtheorem{example}[theorem]{Example}
\newtheorem{quest}[theorem]{Question}
\newtheorem{conjecture}[theorem]{Conjecture}
\renewenvironment{proof}{\smallbreak\noindent {\sc Proof \;---\;}}{\hfill\qedbox}

\counterwithout{theorem}{chapter}
\numberwithin{theorem}{chapter}

\numberwithin{equation}{chapter}

\DeclareOldFontCommand{\rm}{\normalfont\rmfamily}{\mathrm}
\DeclareOldFontCommand{\sf}{\normalfont\sffamily}{\mathsf}
\DeclareOldFontCommand{\tt}{\normalfont\ttfamily}{\mathtt}
\DeclareOldFontCommand{\bf}{\normalfont\bfseries}{\mathbf}
\DeclareOldFontCommand{\it}{\normalfont\itshape}{\mathit}
\DeclareOldFontCommand{\sl}{\normalfont\slshape}{\@nomath\sl}
\DeclareOldFontCommand{\sc}{\normalfont\scshape}{\@nomath\sc}

\usepackage{cancel}

\usepackage{comment}

\usepackage{soul}
\DeclareSymbolFont{newfont}{OML}{cmm}{m}{it}%
\DeclareMathSymbol{\Varrho}{3}{newfont}{37}

\begin{document}

\begin{center}

{\LARGE \textbf{The Pseudocentre of a Group}\footnote{The authors are members of the non-profit association ‘‘AGTA --- Advances in Group Theory and Applications’’ (www.advgrouptheory.com), and are supported by GNSAGA (INdAM). Funded by the European Union - Next Generation EU, Missione 4 Componente 1 CUP B53D23009410006, PRIN 2022- 2022PSTWLB - Group Theory and Applications.}\\ {\large(with an appendix by Anthony Genevois)}}

\medskip

{M. Brescia -- B. Di Siena -- E. Ingrosso -- M. Trombetti}

\end{center}

\bigskip\bigskip

\begin{abstract}
In 1973, Jim Wiegold introduced the concept of pseudocentre $P(G)$ of a group~$G$ as the intersection of the normal closures of the centralizers of its elements. He proved that the pseudocentre of a non-trivial finite group is always non-trivial, giving a new variable on which one can use induction in finite group theory. In the same paper, Wiegold states that no obvious relations seem to hold between the pseudocentre and the ‘‘canonical’’ characteristic subgroups of a group. 

The aim of this work is to show that the pseudocentre is indeed much more involved in the structure of an arbitrary group then anyone could have expected. For example, we prove that a soluble group coincides with its pseudocentre if and only if it is abelian (see Corollary~\ref{pseudodercentr}), and that the structure of the commutator subgroup strongly influences the structure of the pseudocentre (see Theorem \ref{lowanupp} and its corollaries). And this is not the end of the story. In fact, the behaviour of the pseudocentre in arbitrary (possibly infinite) groups can be extremely wild:  sometimes it is very difficult even to understand whether the pseudocentre is trivial or not. This wilderness is exampled by some of our main results (see the introduction for a complete list):

\begin{itemize}
    \item There exists a polycyclic group of Hirsch length $3$ in which the pseudocentre is trivial (see Example \ref{expoly}).
    \smallskip
    \item The pseudocentre of the group of unitriangular matrices over any field is the largest term of the upper central series that is abelian (see~The\-o\-rem~\ref{theoincrediblemclain}).
    \smallskip
    \item Free products have a trivial pseudocentre (see Theorem \ref{theoremafinalefreeproduct}), but there exist amalgamated free products of non-trivial groups coinciding with their pseudocentre (acylindrically hyperbolic groups are discussed by Anthony Genovois in Section \ref{appendice}).
    \smallskip
    \item Weakly branch groups have a trivial pseudocentre,  so, in particular, Grigorchuk groups have trivial pseudocentre (see Theorem \ref{thpetto}). Note also that most free Burnside groups have a trivial pseudocentre as well (see Theorem \ref{burnsidegroups}).
    \smallskip
    \item The pseudocentre of the Thompson group is the derived subgroup (see~The\-o\-rem~\ref{thompson}). 
    \smallskip
    \item Wreath products can have a totally arbitrary pseudocentre (see~Chap\-ter~\ref{chapterwreathproducts}). In fact, using wreath products we can construct a non-trivial group $G$ which is isomorphic to both its pseudocentre~$P(G)$ and the factor~$G/P(G)$ (see Theorem \ref{exsuperincredible}).
\end{itemize}

Finally, we remark that the most beautiful thing about  this is that in order to prove many of the previous results, one needs techniques ranging from many different realms of mathematics. For example, the pseudocentre of $\operatorname{PSL}(2,\mathbb Z)$ is computed by using Pell's equations and Chebyshev polynomials of the first and second kind, while Example \ref{expoly} allows us to state a conjecture about the existence of certain Fibonacci numbers in terms of the pseudocentre of the group. Note that we have explicitly stated all relevant open problems in our work.

\end{abstract}

\medskip

\noindent{\it Mathematics Subject Classification \textnormal(2020\textnormal)}: Primary 20E05, 20E06, 20E08, 20E18, 20E22, 20E34, 20E45, 20F12, 20F14, 20F16, 20F18, 20F19, 20F24, 20F65, 20H20

\medskip

\noindent{\it Keywords}: pseudocentre, linear group, McLain group, wreath product, tree, branch group, free group, free product, amalgamated free product, Thompson group, integrability

\bigskip

\renewcommand{\contentsname}{\large\bf Contents\\[-0.3cm]}

\begingroup
\let\clearpage\relax
\tableofcontents
\endgroup

\chapter*{Notation}

\begin{tabular}{p{1.75cm}p{10cm}}

&{\bf set theory}\\
$\mathbb N,\mathbb N_0,\mathbb Z$ & set of positive integers, non-negative integers, integers\\
$\mathbb Q,\mathbb R,\mathbb C$ &  field of rational numbers, real numbers, complex numbers\\
$\mathbb R^+, A^+$ & set of all positive real numbers (of $A$)\\
$\mathbb P$ & set of all prime numbers\\
$\pi'$ & $\mathbb P\setminus\pi$\\
$p'$ & $\mathbb P\setminus\{p\}$\\
$\pi(n)$ & the set of all primes dividing the integer $n$\\
$\lfloor r\rfloor$ & the only integer $n$ such that $n\leq r<n+1$\\
$\lceil r\rceil$ & the only integer $n$ such that $n-1< r\leq n$\\
$X\subseteq Y$ & $X$ is a subset of $Y$\\
$X\subset Y$ & $X$ is a proper subset of $Y$\\
$X\setminus Y$ & set of all elements in $X$ which are not contained in $Y$\\
$X\cap Y$ & intersection of the sets $X$ and $Y$\\
$X\cup Y$ & union of the sets $X$ and $Y$\\
$(a,b)$ & greatest common divisor of $a$ and $b$, where $a,b\in\mathbb Z$\\
$a\,|\,b$ & $a$ divides $b$, where $a,b\in\mathbb Z$\\
$a\,\cancel{\mid}\, b$ & $a$ does not divide $b$, where $a,b\in\mathbb Z$\\
$|X|$ & cardinality of the set $X$\\
$\omega,\aleph_0$ & the first ordinal and cardinal number\\
$\operatorname{Id}_X$ & the identity function on the set $X$\\
$S^\circ$ & the interior of the subset $S$ of $\mathbb R$\\
$\operatorname{cl}(S)$ & the closure of the subset $S$ of $\mathbb R$\\[0.4cm]

\end{tabular}

\begin{tabular}{p{1.75cm}p{10cm}}

& {\bf group theory}\\

\{1\} & the trivial group\\
$\mathfrak X$ & a class of groups (or group class), that is a collection (not a set) of groups containing the trivial groups and closed with respect to forming isomorphisms --- note that every such a class can naturally be thought of as a property pertaining $G$ (e.g. the class of abelian groups is the property of being abelian).\\
$|G:H|$ & index of the subgroup $H$ in the group $G$\\
$C_G(S)$ & $\{g\in G\, :\, sg=gs\}$, centralizer of the set $S$ in the group $G$\\
$Z(G)$ & centre of the group $G$\\
$Z_\alpha(G)$ & the $\alpha$-th term of the upper central series of the group $G$, where $\alpha$ is any ordinal number\\

\end{tabular}

\begin{tabular}{p{1.75cm}p{10cm}}

$N_G(S)$ & $\{g\in G\,:\, Sg=gS\}$, normalizer of the set $S$ in the group $G$\\
$H^G$ & the normal closure of $H$ in the group $G$, that is, the smallest normal subgroup of $G$ containing $H$\\
$H_G$ & the normal core of $H$ in the group $G$, that is, the largest normal subgroup of $G$ contained in $H$\\
$P(G)$ & $\bigcap_{g\in G}C_G(g)^G$, the pseudocentre of the group $G$\\
$H\leq G$ & $H$ is a subgroup of the group $G$\\
$H<G$ & $H$ is a proper subgroup of the group  $G$\\
$H\trianglelefteq G$ & $H$ is a normal subgroup of $G$\\
$\langle S\rangle$ & the subgroup generated by the subset $S$ of the group $G$\\
$[x,y]$ & $x^{-1}y^{-1}xy$, the commutator of the elements $x$ and $y$ of a group\\
$[H,K]$ & the subgroup generated by the commutators $[h,k]$, where $h\in H$ and $k\in K$\\
$G'$ & $[G,G]$, the commutator subgroup (or derived subgroup) of the group $G$\\
$\gamma_\alpha(G)$ & the $\alpha$-th term of the lower central series of the group $G$ --- here, $\gamma_1(G)=G$ and $\gamma_2(G)=G'$\\ 
$[H,_nK]$ & $[H,K],\ldots,K]$, where $H$ and $K$ are subgroups of a group $G$, and~$K$ is repeated $n$ times\\
$x^y$ & $y^{-1}xy$, conjugate of $x$ by $y$, where $x,y$ are elements of a group\\
$G^n$ & $\langle x^n\,:\, x\in G\rangle$\\
$\operatorname{exp}(G)$ & the exponent of the group $G$, that is, the smallest positive integer $n$ such that $g^n=1$ for all $g\in G$.\\
$\operatorname{Aut}(G)$ & automorphism group of the group $G$\\
$\operatorname{Inn}(G)$ & inner automorphism group of the group $G$\\
$\operatorname{End}(G)$ & set of endomorphisms of the group $G$\\
$H\times K$ & direct product of the groups $H$ and $K$\\
$H\ltimes K$ & semidirect product of the groups $H$ and $K$\\
$H\simeq K$ & $H$ and $K$ are isomorphic groups\\
$o(x)$ & order of the element $x$ of the group $G$\\
$\pi(G)$ & set of all primes dividing the order of an element of the group $G$\\
$O^p(G)$ & intersection of all normal subgroups $N$ of the group $G$ such that $G/N$ is a $p$-group (here, $p$ is a prime number)\\
$O_p(G)$ & maximal normal $p$-subgroup of the group $G$ (here, $p$ is a prime number)\\
$O^\pi(G)$ & intersection of all normal subgroups $N$ of the group $G$ such that $G/N$ is a $\pi$-group (here, $\pi\subseteq\mathbb P$)\\
$O_\pi(G)$ & maximal normal $\pi$-subgroup of the group $G$ (here, $\pi\subseteq\mathbb P$)\\
$\operatorname{H}_n(G,M)$ & $n$-th homology group, where $G$ is a group and $M$ a right \hbox{$\mathbb ZG$-module}\\
$\operatorname{H}^n(G,M)$ & $n$-th cohomology group, where $G$ is a group and $M$ a right \hbox{$\mathbb ZG$-module}\\
$\operatorname{Sym}(S)$ & symmetric group on the set $S$\\
$\operatorname{Sym}(n)$ & symmetric group on the set $\{1,\ldots,n\}$\\
$\operatorname{Alt}(S)$ & alternating group on the set $S$\\
$\operatorname{Alt}(n)$ & alternating group on the set $\{1,\ldots,n\}$\\
$\mathbb Z$ & infinite cyclic group (endowed with addition)\\
$\mathbb Z[q]$ & the subring of $\mathbb Q$ generated by $\mathbb Z$ and $q\in\mathbb Q$\\
$\mathbb Z_n$ & cyclic group of order $n$\\
$\mathbb Z_{p^\infty}$ & Pr\"ufer $p$-group\\[0.4cm]

\end{tabular}

\begin{tabular}{p{1.75cm}p{10cm}}

& {\bf matrices}\\
$\operatorname{trace}(A)$ & trace of the matrix $A$\\
$\operatorname{Tr}(n,\mathbb F)$ & group of all triangular matrices of degree $n$ over the field $\mathbb F$\\
$\operatorname{Tr}_1(n,\mathbb F)$ & group of all unitriangular matrices of degree $n$ over the field $\mathbb F$\\
$\operatorname{GL}(n,\mathbb F)$ & general linear group of degree $n$ over the field $\mathbb F$\\
$\operatorname{GL}(n,q)$ & general linear group of degree $n$ over the field of order $q$\\
$\operatorname{SL}(n,\mathbb F)$ & special linear group of degree $n$ over the field $\mathbb F$\\
$\operatorname{SL}(n,q)$ & special linear group of degree $n$ over the field of order $q$\\
$\operatorname{PGL}(n,q)$ & projective general linear group of degree $n$ over the field of order~$q$\\
$\operatorname{PSL}(n,q)$ & projective special linear group of degree $n$ over the field of order~$q$\\
$\operatorname{PSL}(n,\mathbb Z)$ & projective special linear group of degree $n$ over the ring $\mathbb Z$\\
\end{tabular}

\bigskip\bigskip\bigskip

\noindent{\bf Relevant definitions}

\noindent Let $\mathfrak X$ and $\mathfrak Y$ be group classes, $G$ a group, $\pi$ a set of primes, and $p$ a prime. We define the {\it $\mathfrak X$-radical} of $G$ as the product of all normal $\mathfrak X$-subgroups of $G$. For example, if~$p$ is any prime, then $O_p(G)$ is the radical with respect to the property of being a~\hbox{$p$-group.} Now,~$G$ is:

\begin{itemize}
    \item {\it $\mathfrak X$-by-$\mathfrak Y$} if it contains a normal subgroup $N\in\mathfrak X$ such that $G/N\in\mathfrak Y$. For example, if $\mathfrak X$ is the class of all abelian groups, and $\mathfrak Y$ is the class of all finite groups, then $G$ is $\mathfrak X$-by-$\mathfrak Y$ if and only if it is {\it abelian-by-finite}, that is, if it contains a normal abelian subgroup of finite index. 
    \item {\it divisible} if it is abelian and the equation $\mathtt{x}^n=g$ is solvable for every $g\in G$ and $n\in\mathbb Z$.
    \item {\it \v Cernikov} if it is abelian-by-finite and satisfies the minimal condition on subgroups. This is equivalent to requiring that $G$ has an abelian subgroup of finite index which is the product of finitely many Pr\"ufer groups (this is the divisible radical of $G$).
    \item {\it polycyclic} if it is soluble and satisfies the maximal condition on subgroups.
    \item {\it Dedekind} if all its subgroups are normal, and {\it Hamiltonian} if it is Dedekind and not abelian.
    \item a {\it $\pi$-group} (resp. a {\it$p$-group}) if all its elements $x$ are \hbox{\it$\pi$-ele\-ments} (resp. \hbox{\it$p$-ele}\-{\it ments}), that is, $o(x)\in\pi$ (resp. $o(x)=p^n$ for some positive integer~$n$).
    \item an {\it $FC$-group} if all its elements have finitely many conjugates.
    \item is a {\it $BFC$-group} if its elements have boundedly many finite conjugates; note that a well-known theorem of B.H. Neumann characterizes the $BFC$-groups as those groups with a finite commutator subgroup. 
    \item {\it hypercyclic} if it has an ascending normal series with cyclic factors; of course, every supersoluble group is hypercyclic.
    \item {\it hypercentral} if it has an ascending central series; of course, every hypercentral group is hypercyclic.
    \item {\it $Z$-group} if it has a central series, that is, a chain $\mathcal S$ of subgroup connecting~$\{1\}$ to $G$ in which $[K,G]\leq H$ for every pair $(H,K)$ of consecutive subgroups of $\mathcal S$.
    \item {\it locally $\mathfrak X$} if every finitely generated subgroup of $G$ is contained in an~\hbox{$\mathfrak X$-sub}\-group. If the $\mathfrak X$ is closed with respect to forming subgroups, then this is equivalent to requiring that all finitely generated subgroups are $\mathfrak X$-group. Thus, for example, $G$ is {\it locally nilpotent} if and only if its finitely generated subgroups are nilpotent.
    \item a {\it Baer group} if its cyclic subgroups are subnormal; the {\it Baer radical} is defined as the largest normal subgroup that is a Baer group (note that a well-known theorem assures that the product of normal Baer subgroups is still a Baer subgroup).
    \item a {\it Fitting group} if it coincides with its {\it Fitting radical}, that is, the product of all the nilpotent normal subgroups of $G$.
    \item {\it residually $\mathfrak X$} if the intersection of the normal subgroups $N$ with $G/N\in\mathfrak X$ is trivial. The intersection of all normal subgroups with $\mathfrak X$-quotient is usually called the {\it $\mathfrak X$-residual} of $G$. Clearly, it is not always the case that the quotient by the $\mathfrak X$-residual is an $\mathfrak X$-group. Thus, for example, $G$ is {\it residually nilpotent} if its nilpotent residual is trivial, or, in other words, if the intersection of all the normal subgroups with nilpotent quotient is trivial. Moreover, if $G$ is \v Cernikov, then the finite residual of $G$ has finite index and coincides with the divisible radical.
    \item {\it minimal non-$\mathfrak X$} if $G\not\in\mathfrak X$ but all its proper subgroups are $\mathfrak X$-groups.
    \item {\it locally graded} if every non-trivial finitely generated subgroup has a proper subgroup of finite index.
    \item {\it just-infinite} if every non-trivial normal subgroup of $G$ has finite index.
    \item of {\it Hirsch length $n$} if there is a finite series $$\{1\}=G_0\trianglelefteq G_1\trianglelefteq\ldots\trianglelefteq G_m=G$$ in which the factors are either periodic or infinite cyclic, and there are precisely $n$ factors of the latter type.
\end{itemize}

Moreover, if $H$ is a subgroup of $G$, then we say that:
\begin{itemize}
    \item $H$ is a {\it Sylow $\pi$-subgroup} (resp. {\it Sylow $p$-subgroup}) of $G$ if it is a maximal element of the set of all $\pi$-subgroups (resp. $p$-subgroups) of $G$.
    \item $H$ is {\it subnormal} in $G$ if there is a finite series $$H=H_0\trianglelefteq H_1\trianglelefteq \ldots\trianglelefteq H_n=G$$ connecting $H$ and $G$.
    \item $H$ is {\it ascendant} in $G$ if there is an ascending series $$H=H_0\trianglelefteq H_1\trianglelefteq \ldots\,H_\alpha\trianglelefteq H_{\alpha+1}\trianglelefteq\ldots\, H_n=G$$ connecting $H$ and $G$. Clearly, if $H$ is subnormal, then it is also ascendant.
    \item $H$ is {\it $G$-invariant} if it is normal in $G$.
    \item $H$ is the {\it $\mathfrak X$-radical} of~$G$ if $H$ is the product of all normal $\mathfrak X$-subgroups of~$G$; in particular, if $\mathfrak X$ is the class of locally nilpotent groups (resp.~Baer groups, nilpotent groups), then we obtain the {\it Hirsch--Plotkin} radical (resp. Baer radical, Fitting radical).
\end{itemize}

\part{Introduction and structural results}

\chapter{Introduction and main results}

Let $G$ be a group. The {\it pseudocentre} $P(G)$ of $G$ is the intersection of the normal closures of the centralizers of the elements of $G$, that is, $$P(G) = \bigcap_{x \in G}C_G(x)^G.$$ This is a characteristic subgroup of $G$ introduced by Jim Wiegold in \cite{Weig}. In that paper, Wiegold proved that the pseudocentre of a non-trivial finite group is always non-trivial, a remarkable result that allows for a new variable on which to make induction arguments in proving results for finite groups. Clearly, $P(G)$ always contains the centre $Z(G)$ of $G$, but is also somewhat connected to the derived subgroup by its very definition. Despite this fact, in \cite{Weig}, Wiegold states that no obvious relations seem to hold between the pseudocentre and the ‘‘canonical’’ characteristic subgroups. The aim of this work is to show that the pseudocentre is actually much more involved in the structure of an arbitrary group (and in some other algebraic problems) than anyone could have guessed. 

This is achieved in three parts. In the first part, we prove some general results connecting the pseudocentre with other ‘‘canonical’’ characteristic subgroups. Although some of the main results in this part do not have a very long proof, they are quite surprising in many respects, and have many relevant and non-obvious consequences. In fact, the real difficulty about the pseudocentre is discerning what is true from what is false.

Here, we summarize the main structural results we have obtained:

\begin{itemize}

    \item The pseudocentre contains every {\bf minimal normal} subgroup (see Theorem \ref{normsempl}). In particular, groups satisfying the {\bf minimal condition} on normal subgroups have a non-trivial pseudocentre (see~Co\-rol\-la\-ry \ref{weigoldresult}). Moreover, the pseudocentre contains the square of any {\bf infinite cyclic normal} subgroup (see Lemma \ref{normcyclic}).

    \smallskip

    \item Let $A$ be a normal abelian subgroup of a group $G$. If $g$ has order $n$ modulo $C_G(A)$, then $A^n\leq C_G(g)^G$ (see Lemma \ref{divisiblepiece}). As a consequence, if~$G$ is periodic, then $P(G)$ contains every {\bf normal divisible abelian} subgroup (see~Co\-rol\-la\-ry~\ref{divisibleradical}), and a further consequence shows that every {\bf nilpotent-by-finite} group has a non-trivial pseudocentre (see~Co\-rol\-la\-ry~\ref{nilpbyfinite}). In particular, every {\bf \v Cernikov group} has a non-trivial pseudocentre (it actually is ‘‘pseudonilpotent’’), as well as $CC$-groups (see~The\-o\-rem~\ref{ccgroups}).

    \smallskip

    \item Any group inducing a Pr\"ufer group of automorphisms on an abelian normal subgroup has a non-trivial pseudocentre (see Theorem \ref{prufersopra}). As a consequence, every soluble {\bf minimal non-nilpotent} group and every group of {\bf Heieneken--Mohamed type} have a non-trivial pseudocentre (see~The\-o\-rem~\ref{hmnontrivial} and Corollary \ref{hmnontrivialcor}).

    \smallskip

    \item The pseudocentre of the {\bf direct product} is the direct product of the pseudocentres (see Lemma \ref{directproducts}). Also, it turns out that the class of groups coinciding with their pseudocentre is local (see Theorem \ref{localProperty}).

    \smallskip

    \item Let $G$ be a group. Then $P(G)'\leq G''$ and any $G'$-central section of $G$ is also centralized by $P(G)$ (see Theorem \ref{lowanupp}). This is probably one of the most non-obvious relations between the pseudocentre and the commutator subgroup. It implies for example that: (1) if $G$ is soluble and $G=P(G)$, then $G$ is abelian; (2) if $G'$ is hypercentral, then so is also $P(G)$ (see~Co\-rol\-la\-ry \ref{commutatorsubgrouphypercentral}); and that (3) finite groups with cyclic Sylow subgroups (for example finite groups of square-free order) have a cyclic pseudocentre (see~Co\-rol\-la\-ry \ref{cyclicpseudo}).

    \smallskip

    \item Let $H$ be a subgroup of a group $G$. If $H$ has a {\bf central supplement} in~$G$, then $P(H)$ is contained in $P(G)$ (this is something that rarely occurs) and actually has a central supplement in $P(G)$ (see Theorem \ref{Bern1}).
    
    \smallskip

    \item Let $G$ be a group. If $G=P(G)$ and $G'$ is finite, then $G/Z(G)$ is finite (see Theorem \ref{analogschur}). This shows that the converse of the celebrated Schur's theorem holds in the context of groups coinciding with their pseudocentre.

\end{itemize}

In the second part of the work, we focus on computing the pseudocentre of certain relevant examples of groups and constructions. In this part, we show that sometimes the pseudocentre coincides with certain very well-known subgroups (thus offering an alternative description of these subgroups) and at other times it allows to describe new features of the groups in question. We should also note that this is the (most technically difficult) part of the work in which many techniques from different realms of mathematics are employed, and which shows how the pseudocentre is related to other types of problems, for example in number theory.

\begin{itemize}
    \item The pseudocentre of a finite {\bf symmetric group} is the {\bf alternating group}, while infinite symmetric groups coincide with their pseudocentre (see Theorem \ref{thsymalg}).

    \smallskip

    \item There exists a {\bf polycyclic} group of Hirsch length $3$ with trivial pseudocentre (see Example \ref{expoly}). This example allows us to re-state a conjecture about the existence of certain Fibonacci numbers in terms of the pseudocentres of its finite quotients (see Remark \ref{remarkfibonacci}).

    \smallskip

    \item The {\bf Baumslag--Solitar} group $BS(1,n)$, with $|n|\geq2$, has a trivial pseudocentre (see Theorem \ref{baumslagsolitar}).
    
    \smallskip

    \item The pseudocentre of the group of all {\bf unitriangular matrices} over any field is the largest term of the upper central series that is abelian (see~The\-o\-rem~\ref{theoincrediblemclain}). More generally, the pseudocentre of the {\bf McLain's groups} over any field and any linearly ordered poset is computed in Theorems~\ref{theoincrediblemclain} and~\ref{mclaingeneral}. In particular, such groups coincide with their pseudocentre if and only if their associated poset is dense (see Corollary \ref{densemclain}). Applications of this construction are used to produce strange examples of groups with various types of properties in connection with the pseudocentre (see~Examples \ref{norelation}, \ref{exmclainno}, \ref{exmctrivialpseudo} and \ref{exampleminp}).

    \smallskip

    \item The pseudocentre of {\bf wreath products} is investigated in Chapter \ref{chapterwreathproducts}. It is proved for example that the wreath product of infinitely many groups often coincides with its pseudocentre (see Theorem \ref{theowr}), obtaining (among other results) that every group can be subnormally embedded in a group that coincides with its pseudocentre (see Corollary \ref{pseudocentralembedding}). On the other hand, it turns out that the pseudocentre of a wreath product of two groups can be rather arbitrary. This is exemplified by the pseudocentre of wreath products of the form $H\operatorname{\it wr}\operatorname{Alt}(n)$ (see Theorems~\ref{A_5Pseudo} and~\ref{a5easycases}), $H\operatorname{\it wr}\operatorname{Sym}(n)$ (see~The\-o\-rem~\ref{Sym}), and $H\wr K$, where $K$ is infinite cyclic (see~The\-o\-rems~\ref{similarlamp} and~\ref{exsuperincredible}). In particular, Theorem \ref{exsuperincredible} provides an example of a group in which the pseudocentre is isomorphic to the whole group, while Theorem \ref{similarlamp} shows that the {\bf lamplighter group} has a trivial pseudocentre.

    \smallskip

    \item The pseudocentre of the {\bf Rubik's cube group} and of the symmetry group of the Rubik's cube obtained by disassembling and reassembling is computed in Remark \ref{cuborubik}, partially as a consequence of our results on the pseudocentre of wreath products.
    
    \smallskip

    \item The pseudocentre of a non-dihedral {\bf free product} is  trivial (see~The\-o\-rem~\ref{theoremafinalefreeproduct}). In particular, {\bf free groups} have a trivial pseudocentre. Note that some of the proofs here use Chebyshev polynomials of first and second kind, and~Pell's equations. {\bf Amalgamated} free products are also discussed a bit, but the situation here is much more complicated (see the end of Chapter~\ref{freeproductsect}). {\bf Free solvable groups} have a trivial pseudocentre too (see~The\-o\-rem~\ref{wiegoldcrepa}). In Section \ref{appendice}, Anthony Genevois shows that the pseudocentre of an \textbf{acylindrically hyperbolic} group is contained in the finite radical (see Theorem \ref{thm:AcylHyp}). 

    \smallskip

    \item In Chapter \ref{sectmatrice}, we compute the pseudocentre of the ({\bf affine}) {\bf general linear group}, the ({\bf affine}) {\bf special linear group}, and of some of their relevant quotients.

    \smallskip
    
    \item The {\bf free Burnside group} $B(m,n)$, where $m$ is any cardinal number~$\geq2$ and $n$ is an odd integer $\gg 1$, has a trivial pseudocentre (see The\-o\-rem~\ref{burnsidegroups}).

    \smallskip

    \item The full automorphism group of an {\bf infinite rooted tree} has a trivial pseudocentre and the same is true for every {\bf weakly branch} group (see Theorem \ref{thpetto}). In particular, the~{\bf Grigorchuk}, {\bf Gupta--Sidki}, and {\bf Basilica} groups have a trivial pseudocentre (see the end of~Chapter \ref{treesect}).

    \smallskip

    \item The pseudocentre of the {\bf Thompson group} and of the groups of all orien\-ta\-tion-preserving piecewise linear {\bf homeomorphisms} over any open set is the commutator subgroup (see The\-o\-rem~\ref{finaltheoremthompson}).
\end{itemize}

\medskip

The third and final part of the work is devoted to the study of a couple of well-known problems that always come up in the literature when a new relevant characteristic subgroup or property is defined: 
\begin{itemize}
    \item the {\it inverse problem}, that is, understanding which groups (up-to-iso\-mor\-phism) are allowed to be the pseudocentre of a group (see Chapter~\ref{secpseudoint}).
    
    \smallskip
    
    \item the {\it minimality problem}, that is, understanding the structure of groups in which all proper subgroups coincide with their pseudocentre (see Chapter~\ref{propersubgroupsarepseudocentral}).
\end{itemize}

\medskip

We hope that this work will serve as a basic reference in the subject area, as a text for postgraduate studies, and also as a source of new research ideas. This is in fact the reason why we have highlighted the most relevant problems in this new subject, the solution of which will indubitably offer new, precious insights on the role played by the pseudocentre in the structure of a group, and on many other related problems in group theory and in algebra as a whole.

\chapter{The pseudocentre}


The aim of this section is to study the pseudocentre of an arbitrary group. Our first basic result is a huge generalization of the main result of \cite{Weig}, which actually provides a very short and easy proof of the same result.

\begin{theorem}\label{normsempl}
Let $G$ be a group. If $N$ is a minimal normal subgroup of $G$, then $N\leq P(G)$.
\end{theorem}
\begin{proof}
Let $g \in G$. If $C_G(g)^G\cap N=\{1\}$, then $[C_G(g)^G,N]=\{1\}$, so \hbox{$N\leq C_G(g)\leq C_G(g)^G$,} a contradiction. Therefore, $N$ is always contained in $C_G(g)^G$. The arbitrariness of $g$ yields that $N\leq P(G)$ and completes the proof.
\end{proof}

  \begin{corollary}
Let $G$ be a group. If $N$ is a non-trivial normal subgroup of $G$ satisfying the minimal condition on $G$-invariant subgroups, then $N\cap P(G)\neq\{1\}$.
 \end{corollary}

 \begin{corollary}\label{weigoldresult}
Let $G$ be a non-trivial group satisfying the minimal condition on normal subgroups. Then $P(G)$ is not trivial. In particular, every non-trivial finite group has a non-trivial pseudocentre.
 \end{corollary}

\begin{remark}\label{remWiegold}
{\rm Actually, Wiegold's proof of Corollary \ref{weigoldresult} relies on the fact that the intersection of finitely many normal closures of centralizers of elements of an arbitrary non-trivial group~$G$ is not trivial. This is easily proved by induction using the fact that the normal closure of the centralizer of an element $g$ of $G$ cannot have trivial intersection with any non-trivial normal subgroup $N$ of $G$ (otherwise $N$ would centralize $g$ and $N\leq C_G(g)^G$).}
\end{remark}

\begin{remark}
{\rm It is not possible to replace normality with subnormality in Theorem \ref{normsempl}: in fact, the pseudocentre of the dihedral group of order $8$ coincides with the centre.}
\end{remark}

Corollary \ref{weigoldresult} has a couple of relevant special cases: groups with finitely many normal subgroups, and groups with finitely many conjugacy classes of elements. The relevance of these classes of groups stems from the results and constructions in~\cite{hnn},\cite{osin} --- it is, in fact, possible to prove that any torsion-free group can be embedded in a group with only two conjugacy classes of elements.

\medskip

Now, we describe the basic facts and closure properties of the pseudocentre.

\begin{lemma}\label{directproducts}
The pseudocentre of a direct \textnormal(resp. cartesian\textnormal) product is the product of the pseudocentres of the factors.
\end{lemma}
\begin{proof}
We only concern ourselves with proving the case in which a group \hbox{$G=H\times K$} is the direct product of the groups $H$ and $K$, the other cases being proved in a similar fashion. Let~\hbox{$g\in G$.} Then there are $h\in H$ and $k\in K$ such that~\hbox{$g=hk$.} Obvious\-ly,~\hbox{$C_G(g)=C_H(h)\times C_K(k)$}, so $C_G(g)^G=C_H(h)^H\times C_K(k)^K$. Therefore $$
\begin{array}{c}
\displaystyle P(G)=\bigcap_{g\in G}C_G(g)^G=\!\!\!\!\!\bigcap_{h\in H,k\in K}\!\!\!\!\big(C_H(h)^H\times C_K(k)^K\big)\\[0.6cm]
\displaystyle=\!\!\bigcap_{h\in H}C_H(h)^H\!\times\!\bigcap_{k\in K}C_K(k)^K=P(H)\times P(K)
\end{array}
$$ and we are done.
\end{proof}

\medskip

Although the semidirect products of pseudocentral groups need not be pseudocentral (as shown by the consideration of $\operatorname{Sym}(3)$), on some occasions, the semidirect product of non-abelian simple groups may be pseudocentral. This is for example the case of a finite non-abelian simple group acting as a group of inner automorphisms on itself. More generally, we have the following result.

\begin{lemma}
    Let $N$ be a group and $K$ a group of automorphisms of $N$. If $K$ and $N$ are simple, the following conditions are equivalent:
    \begin{itemize}
        \item[\textnormal{(i)}] For every $h \in N$ there is a non-trivial $k \in K$ such that $h^k = h$. 
        \item[\textnormal{(ii)}] $G = K \ltimes N$ is a pseudocentral group. 
    \end{itemize}
\end{lemma} 
\begin{proof}
Clearly, (ii) implies (i). Assume (i). Let $g$ be an element of $G$, write $g=kx$ with~\hbox{$k\in K$} and $x\in N$, and put $C=C_G(g)^G$. By Theorem \ref{normsempl}, $N\leq C$ so if $k$ is not trivial, we obtain $G=C$. On the other hand, if $k=1$, by (i) there is a non-trivial element $h\in K$ centralizing $g$. It follows that $K\leq\langle h\rangle^G\leq C$ and again~\hbox{$G=C$.} 
\end{proof}

\begin{lemma}\label{homomorphicimages}
Let $G$ be a group. If $N$ is a normal subgroup of a group $G$, then $P(G)N/N\leq P(G/N)$.
\end{lemma}
\begin{proof}
This follows at once from the fact that $C_G(g)N/N\leq C_{G/N}(gN)$ for every $g\in G$.
\end{proof}

\medskip

It turns out that the pseudocentre of a group is not necessarily pseudocentral. In fact, the pseudocentre of the symmetric group of degree $4$ is the alternating group of degree $4$, whose pseudocentre has order $4$.

\begin{theorem}\label{thsymalg}
\begin{itemize}
    \item[\textnormal{(1)}] For every $n\geq3$, we have that $P\big(\operatorname{Sym}(n)\big)=\operatorname{Alt}(n)$.
    \item[\textnormal{(2)}] For any infinite set $I$, $\operatorname{Sym}(I)$ is pseudocentral.
    \item[\textnormal{(3)}] For any infinite set $I$, $\operatorname{FSym}(I)$ is pseudocentral.
\end{itemize}
\end{theorem}
\begin{proof}
(1)\quad Let $P$ be the pseudocentre of $G=\operatorname{Sym}(n)$. If $n=4$, then $P=\operatorname{Alt}(4)$ because every element of the Klein subgroup of $\operatorname{Alt}(4)$ is centralized by some odd permutation. Assume $n>4$. If $n$ is odd, then $X=\langle (1\ldots n)\rangle$ is self-centralizing and so $P\leq C_G(X)^G=X^G\leq\operatorname{Alt}(n)$. If $n$ is even, then \hbox{$X=\langle (1\ldots n-1)\rangle$} is self-centralizing, and again $P\leq\operatorname{Alt}(n)$. Since $\operatorname{Alt}(n)$ is simple, the result follows from~Co\-rol\-lary~\ref{weigoldresult}.

\medskip

\noindent(2)\quad If $x$ is any element of $\operatorname{Sym}(I)$, then either the support of $x$ has the same cardinality of $I$, or its complement has that cardinality. In any case, the centralizer of $x$ contains elements whose support has cardinality $I$, and so whose normal closure is $\operatorname{Sym}(I)$. The statement follows.

\medskip

\noindent(3)\quad The proof is similar to that of (2).
\end{proof}

\begin{remark}
{\rm The consideration of $\operatorname{Alt}(4)$ and $\operatorname{Sym}(4)$ shows that in a finite soluble group $G$ we do not have $N\cap P(G)\leq P(N)$ for any normal subgroup $N$ of~$G$.}
\end{remark}

\medskip

Even though there are exceptions (see Theorem \ref{theoincrediblemclain} for the most extraordinary one), the pseudocentre of a group may have no relation with the terms of the upper and lower central series in general.

\begin{example}
There exists a finite \textnormal(metabelian\textnormal) nilpotent group $G$ of class $3$ in which the pseudocentre does not coincide with any term of the upper and lower series of $G$.
\end{example}
\begin{proof}
Recall that a nilpotent group of class $3$ is always metabelian. Now, let $X=\langle x\rangle\ltimes\langle y\rangle$ be the subgroup of $\operatorname{GL}(3,2)$, where     $$x = \begin{pmatrix}
            1 & 1 & 0 \\
            0 & 1 & 0 \\
            0 & 0 & 1
        \end{pmatrix}\quad\textnormal{and}\quad y = \begin{pmatrix}
            1 & 1 & 1 \\
            0 & 1 & 1 \\
            0 & 0 & 1
        \end{pmatrix}.$$ Let $E=\langle a\rangle\times\langle b\rangle\times\langle c\rangle$ be an elementary abelian $2$-group of order $8$, and define \hbox{$G=X\ltimes E$} to be the natural semidirect product. Clearly, $Z(G)=\gamma_3(G)=\langle c\rangle$, $Z_2(G)=G'=\langle y^2,b,c\rangle$, and $Z_3(G)=G$. Since $C_G(abc)=E$, so $$P(G)\leq E\cap Z_2(G)=\langle b,c\rangle.$$ On the other hand, it is easy to see that $b\in P(G)$ and so $P(G)=\langle b,c\rangle$. 
\end{proof}

\begin{example}\label{norelation}
For every pair of non-negative integers $n<m$, there exists a group $G$ such that:
\begin{itemize}
    \item $Z_n(G)\leq P(G)\leq Z_m(G)$.
    \item If $Z_{n'}(G)\leq P(G)\leq Z_{m'}(G)$, then $n'\leq n$ and $m\leq m'$.
\end{itemize}
Moreover, a similar example exists if we consider the lower central terms instead of the upper ones.
\end{example}
\begin{proof}
This example is based on an independent result that we obtain in~Section \ref{secmclain}. Clearly, the consideration of any pseudocentral group allows us to assume $1\leq n$. Let $H_i$ be the group of all upper unitriangular matrices of degree $i$ over a given field; by Theorem \ref{theoincrediblemclain}, $P(H_i)=Z_{\lfloor i/2\rfloor}(H_i)$. Put $G=H_{2n}\times H_{2m}$. Then $P(G)=Z_n(H_{2n})\times Z_m(H_{2m})$ by Lemma \ref{directproducts}. Since $Z_i(G)=Z_i(H_{2n})\times Z_i(H_{2m})$, so we are done. 

The moreover part follows with similar examples by noting that the upper central terms of~$H_i$ coincide with the lower central ones.
\end{proof}

\begin{example}
Let $G=\langle a\rangle\ltimes A$ be a locally dihedral group, where $a$ has order~$2$ and $A$ is abelian with $|A|\geq2$. Then~\hbox{$P(G)=A^2$.} In particular, $P(G)=G'=Z_{n-1}(G)$ whenever $A$ is cyclic of order $2^n\geq4$.
\end{example}

\medskip\smallskip

However, despite the previous examples, there are some obvious and non-obvious relations between the central series and the pseudocentre. For example, the centre is always contained in the pseudocentre. Moreover, if $G$ is a nilpotent group of class $2$, then the centralizers of elements of $G$ are normal, so $P(G)=Z(G)$. This remark easily implies that $P(G)\leq Z_{n-1}(G)$ whenever $G$ is a nilpotent group of class $n$ (see also~Lem\-ma~\ref{homomorphicimages}). Our next result proves one of the non-obvious relations of the pseudocentre with the lower central series and has many relevant consequences, among which the fact that a soluble group coincides with its pseudocentre if and only if it is abelian (see Corollary \ref{pseudodercentr}).

\begin{theorem}\label{lowanupp}
    Let $G$ be a group. 
    \begin{itemize}
        \item[\textnormal{(1)}] $P(G)' \le G''$.
        \item[\textnormal{(2)}] If $H/K$ is a section of $G$ such that $[H,G']\leq K$, then $[H,P(G)]\leq K$. In particular, we have that $C_G(G') \le C_G(P(G))$ and $P(G)\leq C_G\big(Z_2(G)\big)$.
    \end{itemize}
\end{theorem}
\begin{proof}
(1)\quad Let $a,b \in P(G)$. Since $C_G(a)^G = C_G(a)[C_G(a),G]$, we can write $b = z z_1$, where $z \in C_G(a)$ and $z_1 \in [C_G(a),G]$.
    Now, $$[a,zz_1] = [a,z_1][a,z]^{z_1} = [a,z_1].$$ Now, if $a = yy_1$, where $y \in C_G(z_1)$ and $y_1 \in [C_G(z_1),G]$, then 
    $$[a,z_1] = [y,z_1]^{y_1}[y_1,z_1] = [y_1,z_1]$$ and hence
     $[a,b] = [y_1,z_1] \in G''$.

\medskip

(2)\quad  Let $g\in H$ and $x\in P(G)$. Since $C_G(g)^G = [C_G(g),G]C_G(g)$, so there are \hbox{$z_1 \in [C_G(g),G]$} and $z\in C_G(g)$ with $x = zz_1$. It follows that $$[g,x] =[g,zz_1]=[g,z]^{z_1}[g,z_1]\in K.$$ Therefore \hbox{$[H,P(G)]\leq K$}. In particular,  $C_G(G') \le C_G(P(G))$ and, since $Z_2(G)\leq C_G(G')$, we also have that \hbox{$P(G)\leq C_G\big(Z_2(G)\big)$.}~\end{proof}

\medskip

Theorem \ref{lowanupp} allows us to show that nilpotency-like properties on the commutator subgroup of a group strongly influence the pseudocentre (see the preliminary section ‘‘Notation’’ for some of the definitions).

\begin{corollary}\label{zgroupcor}
Let $G$ be a group. If $G'$ is a \hbox{$Z$-group}, then $P(G)$ is a \hbox{$Z$-group}.
\end{corollary}

\begin{corollary}\label{corlocallynilp}
Let $G$ be a locally finite group. If $G'$ is locally nilpotent, then $P(G)$ is locally nilpotent.
\end{corollary}
\begin{proof}
It follows from Theorem 5.27 of \cite{Rob72} that $G'$ is a $Z$-group, so $P(G)$ is a $Z$-group by~Co\-rol\-la\-ry~\ref{zgroupcor}, and hence $P(G)$ is locally nilpotent because it is locally finite.
\end{proof}

\begin{quest}
Let $G$ be a group with a locally nilpotent commutator subgroup. Is $P(G)$ locally nilpotent?
\end{quest}

\begin{corollary}\label{commutatorsubgrouphypercentral}
Let $G$ be a group. If $G'$ is hypercentral \textnormal(resp., nilpotent\textnormal), then $P(G)$ is hypercentral \textnormal(resp., nilpotent\textnormal).
\end{corollary}
\begin{proof}
Let $$\{1\}=Z_0\leq Z_1\leq\ldots\, Z_\alpha\leq Z_{\alpha+1}\leq\ldots\, Z_\mu=G'$$ be the upper central series of $G'$. It follows from Theorem \ref{lowanupp} that $$
\begin{array}{c}
\{1\}\!=\!Z_0\cap P(G)\leq Z_1\cap P(G)\leq\ldots\, Z_\alpha\cap P(G)\qquad\qquad\qquad\qquad\qquad\\[0.2cm]
\qquad\qquad\qquad\qquad\qquad\leq Z_{\alpha+1}\cap P(G)\leq\ldots\,\, Z_\mu\cap P(G)=G'\cap P(G)\leq P(G)
\end{array}
$$ is a central series of $P(G)$, so $P(G)$ is hypercentral of length at most $\mu+1$.
\end{proof}

\begin{corollary}\label{corsupersoluble}
Let $G$ be a hypercyclic \textnormal(resp., supersoluble\textnormal) group. Then~$P(G)$ is hypercentral \textnormal(resp., nilpotent\textnormal).
\end{corollary}
\begin{proof}
Since the derived subgroup of a hypercyclic (resp., supersoluble) group is hypercentral (resp., nilpotent), the result follows at once from~Co\-rol\-la\-ry~\ref{commutatorsubgrouphypercentral}.~\end{proof}

\begin{corollary}\label{cyclicpseudo}
Let $G$ be a finite group whose Sylow subgroups are cyclic. Then $P(G)$ is cyclic.
\end{corollary}
\begin{proof}
This follows at once from Corollary \ref{corsupersoluble} by recalling that $G$ is supersoluble.
\end{proof}

\begin{corollary}
Any finite group of square free order has a cyclic pseudocentre.
\end{corollary}

\begin{example}
There exist finite groups whose Sylow subgroups are abelian but whose pseudocentre is not.\footnote{We thank Anthony Pisani for this observation.}
\end{example}
\begin{proof}
This is easily seen in GAP \cite{GAP4} by considering the groups \texttt{SmallGroup(168,43)}
and \texttt{SmallGroup(300,25)}.
\end{proof}

\medskip

Finally, we point out a useful reduction of the pseudocentre in case there exists a central supplement.

\begin{theorem}\label{Bern1}
Let $G$ be a group and $Z$ a central subgroup of $G$. If $G = HZ$ for some subgroup $H$ of~$G$, then $P(G) =P(H)\cdot Z$, so in particular $P(G)\cap H = P(H)$.
\end{theorem}
\begin{proof}
Let $g$ be any element of $G$, and write $g=zh$, where $z\in Z$ and $h\in H$. Then $C_G(g) = C_G(h) = C_H(h)\cdot Z$, so $C_G(h)^G=C_H(h)^H\cdot Z$ and hence $$C_G(h)^G\cap H=C_H(h)^H\cdot (Z\cap H)=C_H(h)^H.$$ Finally, $$
P(G)=Z\cdot(P(G)\cap H)=Z\cdot \big(\bigcap_{g\in G}C_G(g)^G\cap H\big)
=Z\cdot \bigcap_{x\in H} C_H(x)^H=Z\cdot P(H),
$$ and the statement is proved.
\end{proof}

\medskip

\section{Pseudocentral groups}

\begin{definition}
{\rm A group $G$ is said to be {\it pseudocentral} (or {\it pseudoabelian}) if $G=P(G)$, that is, if it coincides with its pseudocentre.}
\end{definition}

It follows from Lemmas \ref{directproducts} and \ref{homomorphicimages} that the class of pseudocentral groups is closed under taking homomorphic images and direct products. On the other hand, it is easy to see that this class is not closed with respect to forming subgroups: in fact, one just needs to observe that every simple group is pseudocentral and that every group embeds in a simple group.

\begin{theorem}\label{localProperty}
    The class of pseudocentral groups is local.
\end{theorem}
\begin{proof}
Let $a,b\in G$. Then there exists a pseudocentral subgroup $X$ of $G$ that contains $a$ and~$b$. Now, $a\in X=C_X(b)^X\leq C_G(b)^G$. The arbitrariness of $a$ yields that $G=C_G(b)^G$. Hence $G$ is pseudocentral. 
\end{proof}

\medskip

We should also mention that Wiegold \cite{Weig} defined pseudonilpotent groups starting from the pseudocentre in the same way as nilpotent groups are defined starting from the centre. Let $G$ be a group. Define $P_0(G)=G$ and $P_1(G)$. Then if $\alpha$ is any ordinal number, then $P_{\alpha+1}(G)$ is defined by $$P_{\alpha+1}(G)/P_\alpha(G)=P\big(G/P_{\alpha}(G)\big),$$ while we put $P_\lambda(G)=\bigcup_{\gamma<\lambda}P_\gamma(G)$ for every limit ordinal $\lambda$. This is the {\it upper pseudocentral series} of $G$. Then $G$ is said to be {\it pseudohypercentral} if it concindes with some term of the upper pseudocentral series, while $G$ is {\it pseudonilpotent} if~\hbox{$G=P_{n}(G)$} for some finite $n$. The smallest ordinal number $\mu$ for $G=P_\mu(G)$ is the {\it pseudohypercentral length}, and in case $G$ is pseudonilpotent is the {\it pseudonilpotent class}. Clearly, every pseudonilpotent group is pseudohypercentral, but the converse does not hold (see for example Theorem \ref{exsuperincredible}). Also, it immediately stems from~Co\-rol\-la\-ry \ref{weigoldresult} that every finite group is pseudonilpotent, while every group with the minimal condition on normal subgroups is pseudohypercentral. However, the classes of pseudonilpotent and pseudohypercentral groups are so vast that it is very difficult to say something about them. For example, ~Co\-rol\-la\-ry \ref{pseudodercentr} cannot be extended to arbitrary pseudonilpotent groups. In other words, there is almost no relation between the pseudonilpotent class and the length of the upper/lower central series. In fact, for every positive integer $n\geq2$, the dihedral group of order~$2^{n+1}$ has nilpotency class $n$ and pseudonilpotency class $2$. It seems very difficult even to prove an analogue of Theorem \ref{localProperty}. A further closure property that would be nice to prove (but which seems to be even harder) is the inheritance by normal products.  Of course, the consideration of the quaternion group of order $8$ shows that the product of two abelian normal subgroups is not pseudocentral in general, but the following appears to be a natural problem in the infinite context.

\begin{quest}
Let $G$ be a group which is the product of two normal pseudocentral subgroups. Is $G$ pseudonilpotent?
\end{quest}

\medskip

Theorem \ref{lowanupp} gives some interesting insights into the structure of arbitrary pseudocentral groups, as shown by the following result.

\begin{corollary} \label{pseudodercentr}
    Let $G$ be a group. If $G$ is pseudocentral, then $G' = G''$ and $Z_2(G) = Z(G)$. In particular, if $G$ is hypocentral, hypoabelian or hypercentral, then~$G$ is abelian. 
\end{corollary}
\begin{proof}
Since $[Z_2(G),G']=\{1\}$, the statement follows at once from~The\-o\-rem~\ref{lowanupp}.
\end{proof}

\begin{corollary}
Let $G$ be a pseudocentral group. If $G$ is hypercyclic, then $G$ is abelian.
\end{corollary}
\begin{proof}
By Corollary \ref{pseudodercentr}, we only need to show that $G$ is hypercentral. Now, let $\langle x\rangle$ be any normal cyclic subgroup of $G$. Then $C_G\big(\langle x\rangle\big)=C_G(x)$ is normal in~$G$, so it must coincide with $G$ because $G$ is pseudocentral. Thus, $x\in Z(G)$. Since pseudocentrality is inherited by homomorphic images, $G$ is hypercentral and the statement is proved.
\end{proof}

\medskip

In connection with the above results, we provide an example of a perfect group with a proper pseudocentre.

\begin{example}\label{perfectnonpseudo}
There exists a finite perfect group that is not pseudocentral.
\end{example}
\begin{proof}
Let $H=\langle a,b\rangle$ be the subgroup of $\operatorname{SL}(2,11)$ generated by the matrices $$a= \begin{pmatrix}
      4 & 1 \\
    0 & 3
\end{pmatrix}\quad\textnormal{and}\quad b=\begin{pmatrix}
    0 & 3 \\
    7 & 10
\end{pmatrix}.$$
By Proposition 13.7 of \cite{passman}, we have that $\operatorname{SL}(2,5)$ has the following presentation $$\langle x,y,z \,|\,  x^5=y^3=z^2=1,\, x^z=x,\, y^z=y,\, (xy)^2=z\rangle.$$ Let $R$ be the set of these relations. It is easily seen that these relations are satisfied by $a$, $b^2$ and~$(ab^2)^2$. Note that $$(ab^2)^2 = \begin{pmatrix} -1 & 0 \\0 & -1 \end{pmatrix},$$ so $H$ has non-trivial central elements. Thus, $H\simeq\operatorname{SL}(2,5)$. 

Now, let $N=\langle c\rangle\times\langle d\rangle$ be the direct product of two copies of $\mathbb Z_{11}$. Let $$G=H\ltimes \langle c\rangle\times\langle d\rangle$$ be the natural semidirect product of $H$ and $\langle c\rangle\times\langle d\rangle\simeq\mathbb Z_{11}\times\mathbb Z_{11}$. Since $[c,a] = c^{3}$ and $[d,a]=cd^{2}$, so $N=\langle c,d\rangle\le G'$. But $G/N\simeq SL(2,5)$ and hence $G=G'$. Finally, the centralizer of $c$ in $\operatorname{SL}(2,11)$ is $\operatorname{Tr}_1(2,11)$. However, since $11$ does not divide the order of $\operatorname{SL}(2,5)$, so~\hbox{$C_G(c)=N$} and consequently,~\hbox{$P(G)\leq N$.}
\end{proof}

\begin{remark}
{\rm By using GAP, one can easily check that the group constructed in Example \ref{perfectnonpseudo} is the smallest perfect finite groups that is not pseudocentral.}
\end{remark}

\medskip

Finally, we highlight some of the relations between pseudocentral groups and {\it $FC$-groups}, that is groups with finite conjugacy classes of elements. Recall also that a group is {\it $BFC$-group} if it has boundedly finite conjugacy classes of elements, and that a well known theorem of B.H. Neumann characterizes $BFC$-groups as those with a finite commutator subgroup. For results concerning the class of~\hbox{$FC$-groups} and $BFC$-groups the reader is referred to \cite{tomkinson}. In general, a central-by-finite group is~$BFC$ which in turns is $FC$, and none of these implications reverts. Clearly, the direct product of infinitely many finite non-abelian simple groups is a pseudocentral~\hbox{$FC$-group} that is not a central-by-finite, and not even a $BFC$-group. However, as shown by our next result, in the universe of pseudocentral groups, the classes of central-by-finite and $BFC$-groups coincide.

\begin{theorem}\label{analogschur}
Let $G$ be a pseudocentral group. Then $G$ is a $BFC$-group if and only if $G$ is central-by-finite.
\end{theorem}
\begin{proof}
If $G$ is central-by-finite, then $G'$ is finite by a celebrated theorem of Schur, and hence $G$ is $BFC$. Conversely, suppose $G'$ is finite. In this case,~$G/Z_2(G)$ is finite by a theorem of Philip Hall (see for example \cite{Rob72}, Theorem 4.25), and consequently ~The\-o\-rem~\ref{lowanupp} gives that $G/Z(G)$ is finite.
\end{proof}

\begin{remark}\label{remarkfcgroups}
{\rm The pseudocentre of an $FC$-group is always non-trivial, and we can actually say the same for a larger class of groups (see Theorem \ref{ccgroups}). Actually, if we use the standard arguments that are employed in \cite{tomkinson} to show that for an~\hbox{$FC$-group} the property of being locally nilpotent (locally soluble) and hypercentral (hyperabelian) are equivalent, then we have that $G=P_\omega(G)$ for any~\hbox{$FC$-group} $G$. 

The consideration of a direct product of infinitely many finite groups $G_n$ with $P_n(G)\neq G$ shows that the previous result cannot be improved (see Theorem \ref{highpseudonilp}).}
\end{remark}

\chapter{Groups with a non-trivial pseudocentre}

Knowing that the pseudocentre of certain classes of (infinite) groups  is always non-trivial produces (as in the finite case) a new variable on which one may possibly apply (transfinite) induction. Thus, the problem of determining if the pseudocentre is trivial or not is a very relevant one. Unfortunately, as the reader may have already guessed (and as will certainly see in the second part of the work), it is not at all clear under which conditions the pseudocentre of an infinite group is trivial or not.  In this section, we concern ourselves with giving sufficiently broad conditions for the pseudocentre of an infinite group to be non-trivial, and exhibiting relevant examples in which the pseudocentre is trivial. We start by giving two examples of the latter type. The first is a polycyclic one, but first we need the following auxiliary result on certain sequences of Fibonacci numbers.

\begin{lemma}\label{damettere}
Let $T,U$ be positive integers with $T\ge 2$.  
Let $\{f(n)\}_{n=0}^\infty$ denote the Fibonacci sequence defined by
$f(0)=0$, $f(1)=1$, $f(n+1)=f(n)+f(n-1)$ for $n\ge1$.  
For each integer $U\ge1$ define
\[
c_1(U):=\sum_{j=1}^{T-1} f(jU),\qquad
d_1(U):=1+\sum_{j=1}^{T-1} f(jU+1),\qquad
d_2(U):=1+\sum_{j=1}^{T-1} f(jU-1),
\]
and set
\[
D(U):=d_1(U)d_2(U)-c_1(U)^2.
\]

Then for every integer $M$ there exists an integer $N>M$ such that
\[
D(N)>T^2.
\]
\end{lemma}

\begin{proof}
Let $\varphi=\dfrac{1+\sqrt5}{2}$ be the golden ratio and set
$\psi=-\varphi^{-1}$. By Binet's formula, for every $n\ge0$,
\[
f(n)=\frac{\varphi^n-\psi^n}{\sqrt5}.
\]
Define the auxiliary geometric sums
\[
A(U):=\sum_{j=1}^{T-1}\varphi^{jU},\qquad
B(U):=\sum_{j=1}^{T-1}\psi^{jU}.
\]
Using Binet's formula termwise we obtain the exact equalities
\[
c_1(U)=\frac{A(U)-B(U)}{\sqrt5},\qquad
\sum_{j=1}^{T-1}f(jU+1)=\frac{\varphi A(U)-\psi B(U)}{\sqrt5},
\]
\[
\sum_{j=1}^{T-1}f(jU-1)=\frac{\varphi^{-1}A(U)-\psi^{-1}B(U)}{\sqrt5}.
\]
Set
\[
a_1(U):=\frac{\varphi A(U)-\psi B(U)}{\sqrt5},\qquad
a_2(U):=\frac{\varphi^{-1}A(U)-\psi^{-1}B(U)}{\sqrt5},
\]
so that $d_1(U)=1+a_1(U)$ and $d_2(U)=1+a_2(U)$. A straightforward algebraic
computation (expand and simplify using $\psi=-\varphi^{-1}$) yields the identity
\[
D(U)=\big(1+A(U)\big)\big(1+B(U)\big).
\]
(One can verify this identity by direct expansion; the cancellation of all
terms of order $A(U)^2$ is exact.)

We now produce explicit, uniform bounds for the factor $1+B(U)$. Since
$|\psi|=\varphi^{-1}<1$, we have for every $U\ge1$
\[
|B(U)|=\Big|\sum_{j=1}^{T-1}\psi^{jU}\Big|\le\sum_{j=1}^{T-1}|\psi|^{jU}
\le\sum_{j=1}^{T-1}|\psi|^{j}=\frac{|\psi|\big(1-|\psi|^{T-1}\big)}{1-|\psi|}.
\]
Hence $B(U)$ is uniformly bounded (independently of $U$). More usefully,
observe the exact geometric-sum formula
\[
1+B(U)=1+\sum_{j=1}^{T-1}\psi^{jU}=\frac{1-\psi^{UT}}{1-\psi^U}.
\]
Because $|\psi|<1$ we have $1-\psi^{UT}\neq0$ and $1-\psi^U\neq0$ for every $U$,
and, furthermore,
\[
1+B(U)=\frac{1-\psi^{UT}}{1-\psi^U}
\ge \frac{1-|\psi|^{T}}{1+|\psi|}
\stackrel{\text{def}}{=} c_T>0,
\]
where the constant $c_T$ depends only on $T$ (and not on $U$). Thus
$1+B(U)\ge c_T>0$ for every $U\ge1$.

On the other hand, the sum $A(U)$ satisfies the trivial lower bound
\[
A(U)=\sum_{j=1}^{T-1}\varphi^{jU}\ge \varphi^{U(T-1)}.
\]
Combining these estimates with the identity for $D(U)$ we obtain, for all $U\ge1$,
\[
D(U)=(1+A(U))(1+B(U)) \ge A(U)\,c_T \ge c_T\,\varphi^{U(T-1)}.
\]

Since $\varphi>1$ and $T\ge2$, the right-hand side $c_T\,\varphi^{U(T-1)}$
tends to $+\infty$ as $U\to\infty$. Therefore, given any integer $M$,
choose $N>M$ large enough so that
\[
c_T\,\varphi^{N(T-1)}>T^2.
\]
For such an $N$ we have $D(N)>T^2$, which proves the statement.
\end{proof}

\begin{example}\label{expoly}
There exists a polycyclic metabelian group $G$ with trivial pseudocentre. Furthermore, this group is residually pseudonilpotent of class at most~$2$. 
\end{example}
\begin{proof}
Let $A_1 = \langle a_1 \rangle$ and $A_2 = \langle a_2 \rangle$ be  infinite cyclic groups, and put $A = A_1 \times A_2$. Let $x$ be the automorphism of $A$ defined by the matrix $$X=\begin{pmatrix}
1 & 1\\
1 & 0
\end{pmatrix},$$ so $a_1^x = a_1 a_2$ and $a_2^x = a_1$. Let $\{f(n)\}_{n=0}^\infty$ be the Fibonacci sequence. It is easy to see that for every integer $n\ge 1$ $$X^{n} = 
    \begin{pmatrix}
        f(n+1) & f(n) \\
        f(n) & f(n-1)
    \end{pmatrix}$$ and that $$X^{-n}=
    \begin{cases}
        \begin{pmatrix}
        f(n-1) & -f(n) \\
        -f(n) & f(n+1)
    \end{pmatrix} & \textnormal{$n$ is even}\\[0.6cm]
    \begin{pmatrix}
        -f(n-1) & f(n) \\
        f(n) & -f(n+1)
    \end{pmatrix} & \textnormal{$n$ is odd}
    \end{cases}
    $$

    Let $G = \langle x \rangle \ltimes A$ be the natural semidirect product. First of all, we observe that $C_G(x) = \langle x \rangle$, so $P(G)\leq C_G(a)^G=A$ for all $a\in A$.

    We claim that $\big[A,\langle x^n \rangle\big]\leq A^{f(n/2)}$, whenever $12$ divides $n\neq0$. We start by noticing that $$\big[A,\langle x^n \rangle\big] = \big\langle [a_1,x^n], [a_2,x^n] \big\rangle= \big\langle a_1^{f(n+1)-1} a_{2}^{f(n)}, a_{1}^{f(n)} a_{2}^{f(n-1)-1} \big\rangle.$$ Now, it follows from Theorem 4.7 of \cite{DanGuy} that \[f(n/2) = \operatorname{gcd}\big(f(n+1)-1,f(n+2)-1\big).\] On the other hand, the definition of Fibonacci numbers easily yields that \[\operatorname{gcd}\big(f(n+1)-1,f(n+2)-1\big)= \operatorname{gcd}\big(f(n),f(n+1)-1\big) = \operatorname{gcd}\big(f(n),f(n-1)-1\big),\] so $f(n/2)$ divides $f(n+1)-1$, $f(n)$ and $f(n-1)-1$. Therefore \[[A,\langle x^{n}\rangle] \le A^{f(n/2)} \] and the claim is proved.

    \medskip

Now, let $p$ be an odd prime. We want to find an element  $w=x^{\Xi(p)}a\in G$ such that:
\begin{itemize}
    \item $a\in A^p$;
    \item $12$ divides $\Xi(p)$;
    \item $C_G(w)=\langle w\rangle$;
    \item $p$ divides $f(\Xi(p)/2)$.
\end{itemize}

\noindent If we manage to find such an element, then we have $P(G)\leq [\langle w\rangle,G]\leq A^{f(\Xi(p)/2)}A^p=A^p$ and the arbitrariness of $p$ yields $P(G)=\{1\}$.

\medskip

Start by recalling that for any integer $\ell$, the sequence of Fibonacci numbers taken modulo~$\ell$ is periodic (the period being called the Pisano period). Recall also that $r$ divides $s$ if and only if $f(r)$ divides $f(s)$. Thus, if we want to find a large enough positive integer $s$ such that~$f(s)$ is divisible by a number $x$, we only need to find an $r$ such that $f(r) \equiv_x 0$ and take $s$ to be any multiple of $r$. Then $r$ divides $s$, so $f(r)$ divides $f(s)$, and thus $x$ divides $f(s)$. These observations easily allow us to choose positive integers $n$ and $r_p$ such that $f(r_p)\equiv_p0$, $n\equiv_{12}0$, and $f(n/2)\equiv_{f(12r_p)p^2}0$. Put $w_1=x^na_1^{p}$. Clearly, $C_G(w_1)\cap A=\{1\}$, so~$C_G(w_1)$ is cyclic.

%
%


Let $w_1=x^na_1^p$, and $y_1=x^ua_1^{v_1}a_2^{v_2}$, $u>0$, an element centralizing $w_1$ and such that $y_1^t=w_1$ for some $t>1$ (note that $C_G(w_1)$ is cyclic); in particular, $n=ut$. Then $$
\begin{array}{c}
x^na_1^p=(x^na_1^p)^{x^ua_1^{v_1}a_2^{v_2}}=x^naa_1^{pf(u+1)}a_2^{pf(u)},
\end{array}
$$ where $a\in A^{f(n/2)}$. Now, modulo $A^{f(n/2)}$, the previous equation gives that $f(n/2)$ divides~\hbox{$p\cdot f(u)$.} Then $f(12r_p)p^2$ divides $pf(u)$, so $pf(12r_p)$ divides $f(u)$, so in particular $12r_p$ divides~$u$ and $f(u)\equiv_p0$. Similarly,  $pf(u+1)-p\equiv_{p^2}0$, so $f(u+1)\equiv_p1$ and hence $f(u-1)\equiv_p1$.

We want to show that $p$ divides $f(u/2)$. Since $u$ divides $n$, then $\operatorname{gcd}(n/2,u)=u$ or $u/2$. In the latter case $p$ divides $f(u/2)$ because $f(\operatorname{gcd}(n/2,u))=\operatorname{gcd}(f(n/2),f(u))$. Assume the former. Since $2r_p$ divides $u$, so $r_p$ divides $u/2$. Consequently, $f(r_p)$ divides $f(u/2)$, and hence~$p$, which divides $f(r_p)$, must also divide $f(u/2)$. The claim is proved.

\medskip



Next, we try to obtain some information on the pair $(t,u)$. Using that $y_1^t = w_1$, we have that $$x^n a_1^p = x^{tu} (a_1^{v_1} a_2^{v_2})^{x^{(t-1)u}} \ldots a_1^{v_1} a_2^{v_2}.$$ Now, $$
\begin{array}{c}
a_1^p = \big(a_1^{f((t-1)u +1)v_1} a_1^{f((t-2)u+1)v_1} \ldots a_1^{f(u +1)v_1}a_1^{v_1}\big)\cdot \big(a_2^{f((t-1)u)v_1} \ldots a_2^{f(u)v_1}\big)\cdot \\[0.2cm]
\cdot\big(a_1^{f((t-1)u)v_2} \ldots a_1^{f(u)v_2}\big)\cdot\big(a_2^{f((t-1)u-1)v_2} \ldots a_2^{f(u-1)v_2}a_2^{v_2}\big),
\end{array}
$$ so $$
\begin{cases}
    p = \big(f((t-1)u +1) + f((t-2)u +1) + \ldots + f(u+1)+ 1\big)v_1 + \big(f((t-1)u) + \ldots + f(u)\big)v_2\\[0.2cm]
    0 = \big(f((t-1)u) + \ldots + f(u)\big)v_1+ \big(f((t-1)u-1) + \ldots + f(u-1) +1\big)v_2
\end{cases}$$ and hence $$
\begin{cases}
    p = d_1v_1 + c_1v_2\\
    0 = c_1v_1 + d_2v_2
\end{cases},$$ where $d_1=f((t-1)u +1) + f((t-2)u +1) + \ldots + f(u+1) + 1$, $c_1=f((t-1)u) + \ldots + f(u)$, and $d_2=f((t-1)u-1) + \ldots + f(u-1) +1$. Clearly, $d_1 > c_1$, $d_2>0$, $d_1 = c_1 +d_2$, $$v_1 = - \frac{d_2p}{c_1^2 - d_1d_2}\quad\textnormal{ and }\quad v_2 = \frac{c_1p }{c_1^2 - d_1d_2}.$$

Moreover, $c_1\geq d_2$ when $(t,u)\neq(2,2)$, while $c_1>d_2$ when $(t,u)\neq(2,2),(3,1),(3,2)$. Since~$12$ divides $u$, none of these cases arises and hence $c_1>d_2$.

Now, the equation $p=d_1v_1+c_1v_2$ yields that the greatest positive common divisor of $d_1$ and $c_1$ divides $p$, so it is either $1$ or $p$. On the other hand, $d_1 = d_2 + c_1$, so the greatest common divisor of $d_2$ and $c_1$ divides $p$. Thus, $c_1^2 -d_1d_2$, which divides the greatest common divisor of $d_2p$ and $c_1p$, must also divide $p^2$. It follows that $c_1^2 - d_1d_2  \in \{ \pm 1, \pm p, \pm p^2\}$, in particular $|c_1^2 - d_1d_2| \le p^2$.


Since $f(u)\equiv_p0$ and $f(u+1)\equiv_p1$, so $f(ku+1)\equiv_p1$ and $f(ku)\equiv_p0$ for every integer~$k$, and consequently $d_1\equiv_pt$ and $c_1\equiv_p0$. Similarly, $d_2\equiv_pt$. Suppose that $v_1$ and $v_2$ are not divided by~$p$. Then $p$ must divide $c_1^2-d_1d_2$, so it must also divide either $d_1$ or $d_2$ (and so both), and hence it must divide $t$. But then we may assume $t=p$, so $u=n/p$, and hence we obtain a contradiction by Lemma \ref{damettere}. This contradiction shows that $v_1$ and $v_2$ (by their expressions above) are $\equiv_p0$. Our choice for $w$ is the generator of~$C_G(w_1)$. The claim (and the first half of the statement) is thus proved.

%
%
%
%
%

\medskip



\medskip

Finally, if $p$ is any prime for which the polynomial $x^2-x-1$ is irreducible over $F_p$, then $A/A^p$ is easily seen to be a minimal normal subgroup of~$G/A^p$, so $P(G/A^p)= A/A^p$ and $G/A^p$ is pseudonilpotent of class~$2$. Since one may find infinitely many primes with this property, we have that the group $G$ is residually pseudonilpotent of class at most~$2$.~\end{proof}

\begin{remark}\label{remarkfibonacci}
{\rm Actually, Example \ref{expoly} may even have all its finite homomorphic images pseudo\-nil\-potent of class at most $2$. In fact, this can be proved to be equivalent to the non-existence of a prime $p$ such that the $p^2$-th Fibonacci number is equivalent to $1$ modulo $p^2$, while the $(p^2-1)$-th Fibonacci number is equivalent to~$0$ modulo $p^2$. Unfortunately, the existence of such Fibonacci numbers seems to be an open question.}
\end{remark}

\medskip

Our second example concerns the Baumslag--Solitar groups. Recall that the {\it Baumslag--Solitar groups} $BS(m,n)$, $n,m\in\mathbb Z$, are examples of two-generator one-relator groups that play an important role in combinatorial group theory and geometric group theory as (counter)examples and test-cases (see for example \cite{sampling}). They are given by the group presentation $$\langle a,b\ |\ ba^{m}b^{-1}=a^{n}\rangle.$$

\begin{theorem}\label{baumslagsolitar}
For any integer $|n|\geq2$, the Baumslag--Solitar group $BS(1,n)$ has a trivial pseudocentre.
\end{theorem}
\begin{proof}
Clearly, $G=\langle x\rangle\ltimes A$, where $A=\mathbb Z[1/n]$ and $x$ acts on $A$ as the product by $n$.  Let $p$ be any prime, and define $g=x^{p-1}p$. Note that $C_A(g)=\{1\}$, so $C_G(g)$ is infinite cyclic. Suppose there is $h\in G$ such that $h^q=g$ for some prime~$q$. Write $h=x^ru$ for some positive $r$ and $u\in A$. Then $$h^q=x^{rq}(n^{(q-1)r}u+\ldots+n^ru+u),$$ and hence $rq=p-1$ and $(n^{(q-1)r}+\ldots+n^r+1)u=p$. Write $u=a/b$, where $\pi(b)\subseteq\pi(n)$. Then the previous equation shows that $b=1$, otherwise there would exist a prime dividing $b,n$ and $n^{(q-1)r}+\ldots+n^r+1$. Consequently, $u$ is an integer. Since $n^{(q-1)r}+\ldots+n^r+1\neq1$, so $u=1$ and we obtain a contradiction. This shows that $C_G(g)=\langle g\rangle$.

Now, if $h=x^k\ell$ is any element of $G$ with $\ell\in A$, then $$[g,h]=
[x^{p-1},\ell][p,x^k]=(1-n^{p-1})\ell+(n^k-1)p
$$ belongs to $pA$ because $n^{p-1}\equiv_p1$. Therefore $$C_G(g)^G=\langle g\rangle[g,G]\leq \langle g\rangle pA.$$ Since the normal closure of the centralizer of any non-trivial element of $A$ is $A$, so~\hbox{$P(G)\leq A$.} It follows that $$P(G)\leq \bigcap_{p\in\mathbb P} pA=\{0\},$$ and so $P(G)=\{1\}$.
\end{proof}

\begin{remark}
{\rm Obviously, $BS(1,0)$, $BS(1,\pm1)$ have non-trivial pseudocentres. It is in fact very easy to check to compute them.}
\end{remark}

\begin{quest}
Describe the pseudocentre of all Baumslag--Solitar groups.
\end{quest}

\medskip

Example \ref{expoly} is a polycyclic group of Hirsch length $3$ and trivial pseudocentre. Our next lemma shows that this is the smallest possible Hirsch length for such an example.

\begin{lemma}\label{normcyclic}
Let $G$ be a group. If $N$ is a non-trivial paracentral subgroup of~$G$, then $N\cap P(G)\neq\{1\}$. More precisely, if $N$ is non-periodic, then $N^2\leq P(G)$, while if $N$ is periodic, then $P(G)$ contains the socle of~$N$.
\end{lemma}
\begin{proof}
If $N$ is periodic, then the result follows from Theorem \ref{normsempl}. If $N$ is non-periodic, then the only non-trivial power automorphism of $N$ is the inversion, and hence $N^2\leq C_G(g)^G$ for all $g\in G$. The statement is proved.
\end{proof}

\begin{corollary}
Non-trivial hypercyclic groups have a non-trivial pseudocentre.
\end{corollary}

Our next results deal with group classes that are somewhat duals of that of polycyclic groups.

\begin{lemma}\label{divisiblepiece}
Let $G$ be a group, and $A$ a normal abelian subgroup of $G$. If $g\in G$ has order $n$ modulo~$C_G(A)$, then $A^n\leq C_G(g)^G$. In particular, if $A$ is divisible, then $A\leq C_G(g)^G$.
\end{lemma}
\begin{proof}
If $a\in A$, then $$x=aa^{g} \ldots a^{g^{n-1}}\in C_G(g)\quad\textnormal{ and }\quad y=[a,g][a,g^2]\ldots[a,g^{n-1}]\in\langle g\rangle^G,$$ so $a^n=xy^{-1}\in C_G(g)^G$ and the statement is proved.
\end{proof}

\begin{corollary}\label{Centralfinito} 
Let $G$ be a group and $A$ a non-trivial normal abelian subgroup of $G$. If $G/C_G(A)$ is finite, then~$P(G)\neq\{1\}$.
\end{corollary}
\begin{proof}
Assume by contradiction that $P(G)=\{1\}$.  By Theorem \ref{normsempl}, $G$ has no finite minimal normal subgroups, so the periodic subgroup of $A$ is trivial, so $A$ is torsion-free. Now, the statement follows from Lemma \ref{divisiblepiece}.
\end{proof} 

\begin{corollary}\label{nilpbyfinite}
Any non-trivial nilpotent-by-finite groups has a non-trivial pseudocentre.
\end{corollary}

\begin{corollary}\label{divisibleradical}
Let $G$ be a periodic group. Then $P(G)$ always contains the divisible abelian radical of~$G$.  In particular, any divisible abelian ascendant subgroup of $G$ is contained in $P(G)$.
\end{corollary}

\begin{corollary}
Let $G$ be a \v Cernikov group. Then $G/P(G)$ is finite. 
\end{corollary}

Theorem \ref{similarlamp} shows that Corollary \ref{divisibleradical} does not hold in general.

\begin{corollary}
Let $G$ be a group and $A$ a non-trivial abelian normal subgroup of finite rank of $G$. If $G/C_G(A)$ is periodic, then $P(G)\neq\{1\}$. 
\end{corollary}
\begin{proof}
Suppose by contradiction that $P(G)=\{1\}$. By Theorem \ref{normsempl}, $T=\{1\}$. Now,~$G/C_G(A)$ is a periodic subgroup of $\operatorname{GL}(n,\mathbb Q)$, and hence it is finite. Therefore the statement follows at once from Corollary \ref{Centralfinito}.
\end{proof}

\medskip

Corollary \ref{Centralfinito} shows that when some large centralizers have finite index, then the group has a non-trivial pseudocentre. In fact, it is not difficult to prove that this is the case for every group in which the centralizer of the elements have finite index, the so-called {\it $FC$-groups}. The class of~\hbox{$FC$-groups} has been generalized by Polovickij~\cite{polovickij} to the class of {\it$CC$-groups}, that is, to the class of groups~$G$ with~\v Cernikov conjugacy classes (equivalently,~$G/C_G(x^G)$ is \v Cernikov for every~\hbox{$x\in G$)}. Clearly, every $FC$-group is a $CC$-group, but the consideration of the locally dihedral $2$-group shows that there are $CC$-groups which are not \hbox{$FC$-groups.} Now, using a result of Polovickij we show that non-trivial $CC$-groups have non-trivial pseudocentre.

\begin{theorem}\label{ccgroups}
Let $G$ be a $CC$-group, then $P(G)\neq\{1\}$.
\end{theorem}
\begin{proof}
It follows from Theorem 4.36 of \cite{Rob72} that $[g,G]$ is \v Cernikov for every~\hbox{$g\in G$.} Thus, if $G$ is non-abelian, then there is an element $g\in G$ such that~$[g,G]$ is a non-trivial normal subgroup of $G$, and so $G$ contains a minimal normal subgroup. Thus, $P(G)\neq\{1\}$ by Theorem \ref{normsempl}.
\end{proof}

\begin{remark}
{\rm Let $G$ be a $CC$-group, and let $D$ be the divisible radical of $G$. By Theorem 4.36 of \cite{Rob72} and Lemma \ref{divisiblepiece}, $D\leq P(G)$. Moreover, $G/D$ is covered by finite-by-cyclic normal subgroups, so Lemma \ref{normcyclic} yields that $G/P_2(G)$ is covered by finite normal subgroups and hence is an $FC$-group. By Remark \ref{remarkfcgroups}, $G=P_\omega(G)$. Note that, as for $FC$-groups, this cannot be improved.}
\end{remark}

\begin{remark}
{\rm We cannot expect the pseudocentre of an infinite $FC$-group to be infinite. In fact, the consideration of any infinite extra-special $p$-group for any prime $p$ shows that the pseudocentre of an infinite $FC$-group can be finite.}
\end{remark}

\medskip

The following question concerns with extensions of Corollary \ref{nilpbyfinite}.

\begin{quest}
Do non-trivial nilpotent-by-\v Cernikov groups have a non-trivial pseudocentre?
\end{quest}

In the context of the above question, we have the following partial result.

\begin{theorem}\label{prufersopra}
Let $G$ be a group with an abelian normal subgroup $A$ such that $G/C_G(A)$ is a~Pr\"u\-fer~\hbox{$p$-group} for some prime $p$. Then $P(G)$ is not trivial.
\end{theorem}
\begin{proof}
Set $C=C_G(A)$ and $G=\langle s_i,C\,:\, i \in \mathbb{N}\rangle$, where $o(s_iC)=p^i$ for every~\hbox{$i\in\mathbb N$.}
Let $g \in G\setminus C$. Then $g = s_k^nc$ for some $c \in C$ and $k,n\in\mathbb N$ with $(n,p)=1$.

Clearly, $[\langle g\rangle,A]=[g,A]=[s_k^n,A]$ because $A$ is abelian and normal in $G$. Also,~$n$ is prime with $p$, so there is $m$ such that $mn\equiv_{p^k}1$. Thus, $s_k=s_k^{nm}c_1$ for some~\hbox{$c_1\in C$,} so $$[s_k,A]=[s_k^{nm}c_1,A]=[s_k^{nm},A]\leq [s_k^n,A],$$ and hence $[s_k,A]=[s_k^n,A]$. It follows that $C_A(g)^G\geq [s_k,A]\geq [s_1,A]$. Consequently, $$P(G)\geq \bigcap_{g\in G}C_G(g)^G\geq [s_1,A]\neq\{1\}$$ and the statement is proved.
\end{proof}

\medskip

Theorem \ref{prufersopra} has some rather interesting consequences concerning the so-called groups of Heineken--Mohamed type. Recall that a group is said to be of {\it Heineken--Mohamed  type} if it is not nilpotent and all of its proper subgroups are nilpotent and subnormal. The first (metabelian) example of such a group has been discovered by Heineken and Mohamed \cite{casolo48}, but similar groups and construction were subsequently studied by many other authors (see \cite{casolo9},\cite{casoloB},\cite{casolo41},\cite{casolo42},\cite{casolo73},\cite{casolo76}). It is well-known that a group $G$ of Heineken--Mohamed type is a (countable) soluble~\hbox{$p$-group} for some prime $p$ such that $G/G'$ is a Pr\"ufer $p$-group. It should be noted that the example constructed in \cite{casolo48} is actually one in which the centre is trivial. Our next result shows that no-matter which group of Heineken--Mohamed type we choose, its pseudocentre is always non-trivial.

\begin{theorem}\label{hmnontrivial}
Let $G$ be a group of Heineken--Mohamed type. Then \hbox{$P(G)\!\neq\!\{1\}$.}
\end{theorem}
\begin{proof}
Let $Z$ be the centre of $G'$. If $Z\leq Z(G)$, then $Z\leq P(G)\neq\{1\}$. Otherwise, $G/C_G(Z)$ is a Pr\"ufer $p$-group for some prime $p$ and the statement follows from Theorem \ref{prufersopra}.
\end{proof}

\begin{corollary}\label{hmnontrivialcor}
Let $G$ be a locally graded minimal non-nilpotent group. Then $P(G)\neq\{1\}$.
\end{corollary}
\begin{proof}
This follows at once from Theorem 3.4.2 of \cite{casoloB} and Theorem \ref{hmnontrivial}.
\end{proof}

\begin{quest}
Is there a metabelian $p$-group $G$ with trivial pseudocentre for some prime $p$?
\end{quest}

\part{The pseudocentre of relevant groups and constructions}

\chapter{Unitriangular matrices and McLain's constructions}\label{secmclain}

Let $\Lambda$ be a partially ordered set, $\mathbb F$ a field, and $V$ a vector space on $\mathbb F$ with basis $\{v_\lambda : \lambda \in \Lambda\}$. For all $\lambda, \mu \in \Lambda$ with $\lambda < \mu$, define the linear application $e_{\lambda,\mu}$ as follows:
\begin{center}
    $e_{\lambda,\mu} (v_{\lambda}) = v_{\mu}$ and $e_{\lambda,\mu} (v_{\zeta}) = 0$ for every $\zeta \neq \lambda$. 
\end{center}
 Clearly, $1 + a e_{\lambda,\mu}$, where $a \in \mathbb F$, is a non-singular linear application and $$(1 + a e_{\lambda,\mu})^{-1} = 1 - a e_{\lambda,\mu}.$$ Furthermore, the following properties hold:
\begin{itemize}
    \item[(i)] $[1+ae_{\lambda,\mu}, 1+ be_{\mu,\zeta}] = 1+abe_{\lambda,\zeta}$ for all $\lambda,\mu,\zeta \in \Lambda$. 
    \item[(ii)] $[1+ ae_{\lambda,\mu}, 1 + be_{\nu,\zeta}] = 1$ for all $\lambda,\mu,\nu,\zeta\in\Lambda$ with $\zeta \neq \lambda$ and $\mu \neq \nu$.
\end{itemize}

The group $M=M(\Lambda,\mathbb F) = \langle 1 + ae_{\lambda,\mu} : \lambda,\mu \in \Lambda, \lambda < \mu, a \in \mathbb F \rangle$ is the {\it McLain's group determined by~$\Lambda$ and $\mathbb F$} (see \cite{Rob72}, Section 6.2). McLain's groups have been introduced by McLain in \cite{McLain},\cite{McLain2} and provide a great source of examples of locally nilpotent groups that are in some sense locally isomorphic to unitriangular groups of matrices. In fact, if $\Lambda$ is finite of order $n$ and linearly ordered, then $M$ is isomorphic to the group $\operatorname{Tr_1}(n,\mathbb F)$ of all upper unitriangular matrices on $\mathbb F$. Moreover, $M$ is always a Fitting group  because it is generated by the abelian normal subgroups $\langle 1 + a e_{\lambda,\mu}\rangle^M$. Note also that~$M(\Lambda,\mathbb F)$ is perfect if and only if $\Lambda$ is dense, and that~$M(\Lambda,\mathbb F)$ has trivial centre if and only if there are no two distinct elements~$\lambda<\mu$ of $\Lambda$ such that $\lambda$ is minimal and $\mu$ is maximal.

Our first main result provides  a charming description of the pseudocentre of~$\operatorname{Tr}_1(n,\mathbb F)$ as the last abelian term (resp., the first abelian term) of the upper central series (resp., of the lower central series), but first we need the following remark.

\begin{remark}\label{remconjmclain}
{\rm Let $\mathbb F$ be a field, $\Lambda$ a poset, and let $$g=1+a_{\alpha_1,\beta_1}e_{\alpha_1,\beta_1}+\ldots+a_{\alpha_t,\beta_t}e_{\alpha_t,\beta_t},$$ where $a_{\alpha_i,\beta_i}\neq0$ and $(\alpha_i,\beta_i)\neq(\alpha_j,\beta_j)$ for all $i\neq j$. If $x=1+ae_{u,v}$, then $g^x=g+h$, where~$h$ is a sum of terms $be_{r,s}$ for which there exists $i$ with $r<\alpha_i$ or~\hbox{$\beta_i<s$.}}
\end{remark}

\begin{theorem}\label{theoincrediblemclain}
Let $\mathbb F$ be a field, $\Lambda$ a finite linearly ordered set of cardinality $n$, and put $M=M(\mathbb F,\Lambda)$. Then $P(M)=Z_\rho(M)$, where $\rho=\lfloor n/2\rfloor$.
\end{theorem}
\begin{proof}
Without loss of generality, we may assume that $\Lambda=\{1,\ldots,n\}$ with its natural ordering. First, we show that $P(M)\geq Z_\rho(M)$. Let $g$ be an arbitrary element of $M$. We claim that~$C_M(g)^M$ contains $Z_\rho(M)$, that is, $C_M(g)^M$ contains all the elements of the form $1+ae_{u,v}$, where $a\neq0$, $1\leq u\leq \lfloor n/2\rfloor$ and $v\geq u+\lceil n/2\rceil$. Clearly, \hbox{$Z_1(M)\leq C_M(g)^M$.} Assume by induction that $C_M(g)^M\geq Z_{\ell}(M)$, where $1\leq \ell<\rho$: we need to prove that \hbox{$C_M(g)^M\geq Z_{\ell+1}(M)$.} In the following, whenever the reader runs into a commutator of the form $[1+e_{\xi,\xi},g]$, with $g\in M$ and $\xi\in\Lambda$, they must consider it equal to $g$.

Write $$g=1+a_{\alpha_1,\beta_1}e_{\alpha_1,\beta_1}+\ldots+a_{\alpha_t,\beta_t}e_{\alpha_t,\beta_t},$$ where $(\alpha_i,\beta_i)\neq(\alpha_j,\beta_j)$ if $i\neq j$, and $a_{\alpha_i,\beta_i}\neq0$ for every $i$. Of course, if $g$ centralizes $1+ae_{x,y}$ with $u\leq x<y\leq v$, then $$1+ae_{u,v}=\big[1+e_{u,x},[1+ae_{x,y},1+e_{y,v}]\big]\in C_M(g)^M.$$ Assume~$g$ does not centralize any element of the form $1+ae_{x,y}$ with $u\leq x<y\leq v$. This means that for every such $x$ and $y$, we can find $i$ and $j$ such that $\alpha_i=y$ or~\hbox{$\beta_j=x$.} This condition, which we refer to as condition ($\star$), will be applied several times in the following arguments without always being mentioned. The proof will be carried out in a number of steps.

\medskip

\noindent Case 1a: $\beta_i\neq u$ for every $i$.\\
First of all, we re-order the terms in the expression of $g$ as follows: $$
\begin{array}{c}
g=1+a_{\alpha_1,\beta_{1,1}}e_{\alpha_1,\beta_{1,1}}+\ldots+a_{\alpha_1,\beta_{1,n_1}}e_{\alpha_1,\beta_{1,n_1}}+\ldots\qquad\qquad\qquad\qquad\\[0.2cm]
\qquad\qquad\qquad\qquad\ldots+a_{\alpha_t,\beta_{t,1}}e_{\alpha_t,\beta_{t,1}}+\ldots+a_{\alpha_t,\beta_{t,n_t}}e_{\alpha_t,\beta_{t,n_t}},
\end{array}
$$ where $\alpha_i<\alpha_j$ and $\beta_{k,i}<\beta_{k,j}$ if $i<j$. Note that one can write $g$ as follows $$
\begin{array}{c}
(1+a_{\alpha_t,\beta_{t,n_t}}e_{\alpha_t,\beta_{t,n_t}})\ldots(1+a_{\alpha_t,\beta_{t,1}}e_{\alpha_t,\beta_{t,1}})\ldots\qquad\qquad\qquad\qquad\qquad\qquad\\[0.2cm]
\qquad\qquad\qquad\qquad\qquad\qquad\qquad\ldots(1+a_{\alpha_1,\beta_{1,n_1}}e_{\alpha_1,\beta_{1,n_1}})\ldots(1+a_{\alpha_1,\beta_1}e_{\alpha_1,\beta_1}).
\end{array}$$

By condition ($\star$), there is $\ell$ with $\alpha_{\ell+j}=u+1+j$ for $j=0,\ldots,v-u-1$. Suppose first that there are $s\in\{u+2,\ldots,v-1\}$ and $r\in\{\ell,\ldots,\ell+v-u-2\}$ such that $\beta_{r,j_r}=s$ for some $1\leq j_r\leq n_r$. Then the above product expression for $g$ yields that $$
\begin{array}{c}
[g,1-a_{\alpha_r,\beta_{r,j_r}}^{-1}ae_{u,\alpha_r}]=\\[0.2cm]
=1+a_{\alpha_r,\beta_{r,j_r}}^{-1}a_{\alpha_r,\beta_{r,1}}ae_{u,\beta_{r,1}}+\ldots+ae_{u,\beta_{r,j_r}}+\ldots+a_{\alpha_r,\beta_{r,j_r}}^{-1}a_{\alpha_r,\beta_{r,n_r}}ae_{u,\beta_{r,n_r}}=\\[0.2cm]
=1+a_{\alpha_r,\beta_{r,j_r}}^{-1}a_{\alpha_r,\beta_{r,1}}ae_{u,\beta_{r,1}}+\ldots+ae_{u,s}+\ldots+a_{\alpha_r,\beta_{r,j_r}}^{-1}a_{\alpha_r,\beta_{r,n_r}}ae_{u,\beta_{r,n_r}}
\end{array}
$$ because of the condition $\beta_i\neq u$ for all $i$. It follows that $$\big[[g,1-a_{\alpha_r,\beta_{r,j_r}}^{-1}ae_{u,\alpha_r}],1+e_{s,v}\big]=1+ae_{u,v}\in C_M(g)^M$$ and we are done.

Suppose that there are no such $s$ and $r$. Conjugating by elements of the form $1+e_{\alpha_{\ell+j_1},\alpha_{\ell+j_2}}$, where $j_1<j_2$, we may replace $g$ by a conjugate in such a way that the elements $\beta_{i,1}$ are distinct. Now, since \hbox{$v-u\geq\lceil n/2\rceil\geq n/2$,} so there is~\hbox{$r\in\{\ell,\ldots,\ell+v-u-2\}$} with $\beta_{r,1}=v$. Then $$
\begin{array}{c}
w=[g,1-a_{\alpha_r,\beta_{r,1}}^{-1}ae_{u,\alpha_r}]=\\[0.2cm]
=1+ae_{u,v}+a_{\alpha_r,\beta_{r,1}}^{-1}a_{\alpha_r,\beta_{r,2}}ae_{u,\beta_{r,2}}+\ldots+a_{\alpha_r,\beta_{r,1}}^{-1}a_{\alpha_r,\beta_{r,n_r}}ae_{u,\beta_{r,n_r}}.
\end{array}$$
Finally, replacing $w$ by a suitable conjugate, we can get rid of the extra term and obtain $1+ae_{u,v}\in C_M(g)^M$, as required.

\bigskip

\noindent Case 1b: $\alpha_i\neq v$ for every $i$.\\
This case can be proved symmetrically with respect to the previous case by re-ordering the terms in the expression of $g$ as follows: $$
\begin{array}{c}
g=1+a_{\alpha_{1,1},\beta_1}e_{\alpha_{1,1},\beta_1}+\ldots+a_{\alpha_{1,n_1},\beta_1}e_{\alpha_{1,n_1},\beta_1}+\ldots\qquad\qquad\qquad\qquad\qquad\\[0.2cm]
\qquad\qquad\qquad\qquad\qquad\ldots+a_{\alpha_{t,1},\beta_t}e_{\alpha_{t,1},\beta_t}+\ldots+a_{\alpha_{t,n_t},\beta_t}e_{\alpha_{t,n_t},\beta_t},
\end{array}
$$ where $\alpha_{k,j}<\alpha_{k,i}$ and $\beta_j<\beta_i$ if $i<j$.

\medskip

The proof now proceeds by induction by assuming that there exist $p_1,\ldots,p_t$ and $q_1,\ldots,q_t$ such that $\alpha_{p_1}=v,\ldots,\alpha_{p_t}=v-t+1$, $\beta_{q_1}=u,\ldots,\beta_{q_t}=u+t-1$, and $\beta_{q_t}<\alpha_{p_t}+1$. We prove that we can find $p_{t+1}$ and $q_{t+1}$ with $\alpha_{p_{t+1}}=v-t$ and~\hbox{$\beta_{q_{t+1}}=u+t$,} respectively. 

\medskip

\noindent Case 2a: $\beta_i\neq u+t$ for every $i$.\\
By condition ($\star$), there is $\ell$ such that $\alpha_{\ell+j}=u+t+1+j$ for $j=0,\ldots,v-u-2t-1$. The number of these objects $\alpha_{\ell+j}$, plus the elements $\alpha_{p_1},\ldots,\alpha_{p_t}$, $\beta_{q_1},\ldots,\beta_{q_t}$, gives a total of~$v-u\geq \lceil n/2\rceil$ objects. Rewrite the sum of the terms of $g$ of the form $a_{\alpha_*,\beta_{q_i}}e_{\alpha_*,\beta_{q_i}}$ with fixed~$i$ as follows: $$a_{\alpha_{i,1},\beta_{q_i}}e_{\alpha_{i,1},\beta_{q_i}}+\ldots+a_{\alpha_{i,m_i},\beta_{q_i}}e_{\alpha_{i,m_i},\beta_{q_i}},$$ where $a_{\alpha_{i,j},\beta_{q_i}}\neq0$ and $\alpha_{i,k}>\alpha_{i,l}$ for $k<l$. Similarly, rewrite the sum of the terms of $g$ of the form $a_{\alpha_{p_i},\beta_*}$ with fixed $i$ as follows: $$a_{\alpha_{p_i},\beta_{i,1}}e_{\alpha_{p_i},\beta_{i,1}}+\ldots+a_{\alpha_{p_i},\beta_{i,l_i}}e_{\alpha_{p_i},\beta_{i,l_i}}$$ with $\beta_{i,k}<\beta_{i,l}$ for $k<l$. As in Case 1a, by suitable conjugations  we can assume that $\alpha_{i,1}\neq\alpha_{j,1}$ and $\beta_{i,1}\neq\beta_{j,1}$ for $i\neq j$; note that if some of the previous terms disappears from the expression of $g$, then we are in a previous case, and we are done.

If there is $\alpha_{p_i}\in\{v-t+1,\ldots,v\}$ such that $\beta_{i,1}\in\{u,\ldots,v\}$, then we can argue as in Case~1a by considering a commutator with an element of the form $1+be_{v-t,\alpha_{p_i}}$ and then removing the terms of the kind $ba_{\alpha_{p_i},\beta_{i,j}}e_{v-t,\beta_{i,j}}$ for $j\neq 1$ by means of another suitable conjugation; the same argument can be used also in case we find $\alpha_r\in\{u+t+1,\ldots,v-t\}$ such that $\beta_{r,1}\in\{u,\ldots,v\}$. Hence, suppose there is no such element. Since $v-u\geq \lceil n/2\rceil$, there must be an $i$ such that $\alpha_{i,1}=u$. Then we take a commutator of $g$ with an element of the form $1+be_{\beta_{q_i},v-1}$, and then the commutator of the result with an element of the form $1+be_{v-1,v}$. In this way we obtain an element containing the term $ae_{u,v}$ and only terms of the form $be_{\ast,v}$ in which the first subscript is always $\leq u$ (to see this factorise $g$ as in Case 1a and apply the usual commutator identities).  Finally, suitable conjugations allow us to remove the extra terms and get $1+ae_{u,v}\in C_M(g)^M$.

\medskip

\noindent Case 2b: $\alpha_i\neq v-1$ for every $i$.\\
This is symmetrical to Case 2a in the same way Case 1b is symmetrical to Case 1a.

\medskip

By induction, we can find $t$ in such a way that $\alpha_{p_t}=\beta_{q_t}+1$. Applying the last part of the argument of Case 2a to the elements whose subscripts contains either an~$\alpha_{p_i}$ or a $\beta_{q_j}$, we obtain the result.

\bigskip

Now, it remains to prove that $P(M)$ coincides with $Z_\rho(M)$. Suppose not. Then there exists $$x=1+a_{u_1,v_1}e_{u_1,v_1}+\ldots+a_{u_w,v_w}e_{u_w,v_w}\in \big(P(M)\cap Z_{\rho+1}(M)\big)\setminus Z_\rho(M)$$ where $a_{u_i,v_i}\neq0$ and $v_i-u_i=\lceil n/2\rceil-1$, for every $1\leq i\leq w$. Put $(u_1,v_1)=(u,v)$ and $r=u+(u-1)$. Then there exist injective mappings $$\sigma:\,\{u,\ldots,r-1\}\rightarrow\{1,\ldots,u-1\}\quad\textnormal{and}\quad\tau:\,\{r+1,\ldots,v\}\rightarrow\{v+1,\ldots,n\}$$ such that $\sigma(a)>\sigma(b)$ and $\tau(c)>\tau(d)$, when $a<b$ and $c<d$. Note that if $u=1$ or $r+1>v$, then we do not take into account the function $\sigma$ or $\tau$, respectively. Now, set
$$g=1+\sum_{\substack{u\leq i<r}} e_{\sigma(i),i}+\sum_{\substack{r<i\leq v}} e_{i,\tau(i)}$$ and suppose that the centralizer of $g$ in $M$ contains an element of the form $$y=1+a_{u',v'}e_{u',v'}+a_{\alpha_1,\beta_1}e_{\alpha_1,\beta_1}+\ldots+a_{\alpha_t,\beta_t}e_{\alpha_t,\beta_t}$$ with $a_{u',v'},a_{\alpha_i,\beta_i}\neq0$, $(u',v')\neq (\alpha_i,\beta_i)\neq(\alpha_j,\beta_j)$ for $i\neq j$, and \hbox{$u\leq u'<v'\leq v$.} If~\hbox{$u'<r$}, then~$gy$ contains an extra term $e_{\sigma(u'),v'}$ with respect to $yg$. Thus, $C_M(g)$ does not contain any such element. Since $C_M(g)^M$ consists of products of conjugates of elements of~$C_M(g)$, so $P(M)\leq C_M(g)^M$ does not contain any element of the previous type (see Remark~\ref{remconjmclain}) and so does not contain $x$. We similarly deal with the other cases, so the proof is complete.
\end{proof}

\begin{remark}
{\rm The relevance of the above result can also be seen by the fact that if $p$ is a prime, then every finite $p$-group $G$ can be embedded in $M(\mathbb F_p,\Lambda)$, where~$\mathbb F_p$ is the field of order $p$, and $\Lambda$ is a linearly ordered set of order at least~$|G|$. In fact,~$G$ can be embedded in $\operatorname{GL}(|G|,p)$ using permutation matrices, and thus~$G$ embeds in~$\operatorname{SL}(|G|,p)$ because $\operatorname{GL}(|G|,p)/\operatorname{SL}(|G|,p)$ is isomorphic to the multiplicative group of $\mathbb F_p$; by Sylow's theorem,~$G$ is then isomorphic to a subgroup of $\operatorname{Tr}_1(|G|,p)$. Note, however, that not all groups generated by involutions are isomorphic to McLain's groups over a certain poset, an obvious example being the dihedral group of order~$16$.}
\end{remark}

Now, we move to the case of infinite linearly ordered sets. In the following, if~\hbox{$a,b\in \Lambda$,} then we put $[a,b]=\{c\in\Lambda\,:\, a\leq c\leq b\}$, $(-\infty,a)=\{c\in\Lambda\,:\, c< a\}$ and $(b,+\infty)=\{c\in\Lambda\,:\, b< c\}$.

\begin{theorem}\label{mclaingeneral}
Let $\mathbb F$ be a field and $\Lambda$ an infinite linearly ordered set. Consider an arbitrary element  $$u=1+a_{\alpha_1,\beta_1}e_{\alpha_1,\beta_1}+\ldots+a_{\alpha_m,\beta_m}e_{\alpha_m,\beta_m}$$ of $M=M(\mathbb F,\Lambda)$ with $a_{\alpha_i,\beta_i}\neq0$ for  $i=1,\ldots,m$ and $(\alpha_i,\beta_i)\neq(\alpha_j,\beta_j)$ if $i\neq j$. Then $u\in P(M)$ if and only if the interval $[\alpha_i,\beta_i]$ is infinite for every $i=1,\ldots,m$.

In particular, $P(M)=M(\mathbb F,\Sigma)$, where $\Sigma$ is a poset having the same underlying set as $\Lambda$ but in which $\alpha\leq\beta$ if and only if $[\alpha,\beta]$ is infinite.
\end{theorem}
\begin{proof}
Let $x=1+ae_{\mu,\lambda}$, $a\neq0$ and $g= 1 + a_{\mu_1,\lambda_1}e_{\mu_1,\lambda_1} + \ldots + a_{\mu_n,\lambda_n}e_{\mu_n,\lambda_n}$ for some elements $a_{\mu_i,\lambda_i}\in\mathbb F\setminus\{0\}$ and $(\mu_i,\lambda_i)\neq(\mu_j,\lambda_j)$ if $i\neq j$. We start by showing that if the interval~$[\mu,\lambda]$ is infinite, then $x\in C_M(g)^M$, so in particular~$x$ belongs to the pseudocentre of~$M$ by the arbitrariness of $g$.  To see this, take $r_1,r_2\in\Lambda\setminus\{\mu_i, \lambda_i : i = 1,\ldots,n \}$ with \hbox{$\mu<r_1<r_2<\lambda$.} Thus, $1+ae_{r_1,r_2}\in C_M(g)$, so $$1+ae_{\mu,\lambda}=[1 + ae_{\mu,r_2},1+e_{r_2,\lambda}] = \big[[1 + e_{\mu,r_1},1+ ae_{r_1,r_2}],1+e_{r_2,\lambda}\big] \in C_M(g)^M.$$ It follows from the arbitrariness of $x$ that also $u$ belongs to $P(M)$, provided that the interval $[\alpha_i,\beta_i]$ is infinite for every $i=1,\ldots,m$, and we are done.

Conversely, consider the element $$x=1+a_{u_1,v_1}e_{u_1,v_1}+\ldots+a_{u_w,v_w}e_{u_w,v_w},$$ where $a_{u_i,v_i}\neq0$, $(u_i,v_i)\neq(u_j,v_j)$ for $i\neq j$, and the interval $[u_1,v_1]$ is finite. Since~$\Lambda$ is infinite, so either $(-\infty,u_1)$ or $(v_1,+\infty)$ is infinite. Suppose the former. Then there exists an injective mapping $$\sigma:\,[u_1,v_1]\rightarrow (-\infty,u_1)$$ such that $\sigma(b)<\sigma(a)$ for $a<b$. Now, set
$$g=1+\sum_{\substack{u_1\leq i\leq v_1}} e_{\sigma(i),i}$$ and suppose that the centralizer of $g$ in $M$ contains an element of the form $$y=1+a_{u',v'}e_{u',v'}+a_{\alpha_1,\beta_1}e_{\alpha_1,\beta_1}+\ldots+a_{\alpha_t,\beta_t}e_{\alpha_t,\beta_t}$$ with $a_{u',v'},a_{\alpha_i,\beta_i}\neq0$, $(u',v')\neq (\alpha_i,\beta_i)\neq(\alpha_j,\beta_j)$ for $i\neq j$, and \hbox{$u_1\leq u'<v'\leq v_1$.} Then~$gy$ contains an extra term $e_{\sigma(u'),v'}$ with respect to $yg$, and we have a contradiction. Thus, $C_M(g)$ does not contain any such element. Since $C_M(g)^M$ consists of products of conjugates of elements of~$C_M(g)$, so $P(M)\leq C_M(g)^M$ does not contain any element of the previous type (see Remark~\ref{remconjmclain}) and so does not contain~$x$. 
\end{proof}

\begin{corollary}\label{densemclain}
Let $\mathbb F$ be a field, and $\Lambda$ an infinite linearly ordered set. Then $M=M(\Lambda,\mathbb F)$ is pseudocentral if and only if $\Lambda$ is dense.
\end{corollary}

\begin{example}\label{exmclainno}
There exists an infinite group $G$ satisfying the following properties:
\begin{itemize}
    \item $G$ is a perfect Fitting group.
    \item $G$ has trivial centre.
    \item $G=P(G)$.
\end{itemize}
\end{example}
\begin{proof}
Consider an infinite dense linearly ordered set $\Lambda$ which does not have a minimum or a maximum, and apply Theorem \ref{mclaingeneral}.
\end{proof}

\begin{example}\label{exmctrivialpseudo}
There exists an infinite Fitting group $G$ with trivial pseudocentre. Moreover, $G$ can be chosen to be either a $p$-group for any choice of the prime~$p$, or torsion-free.
\end{example}
\begin{proof}
Take $\Lambda=\mathbb Z$ with the usual ordering and apply Theorem \ref{mclaingeneral}.
\end{proof}

\begin{example}\label{exampleminp}
For each prime $p$, there exists a locally finite group $G$ satisfying $\operatorname{min}$-$p$ and having a trivial pseudocentre.
\end{example}
\begin{proof}
Let $q$ be a prime such that $p\leq q-1$, consider the $q$-group $H$ constructed in~Example~\ref{exmctrivialpseudo}, and put $G=H\wr\langle x\rangle$, where $\langle x\rangle$ is a cyclic group of order $p$. Let $B=H_1\times\ldots\times H_p$ be the base group of~$G$, where $H_i\simeq H$ for every $i=1,\ldots,p$. If $g$ is any non-trivial element of $H_1$, then the automorphism group of $\langle g\rangle$ has at least $p$ elements $\alpha_1,\ldots,\alpha_p$ of $q'$-order. Put $$g_1=g,\, g_2=\big(g^{\alpha_1}\big)^x,\,\ldots,\, g_p=\big(g^{\alpha_p}\big)^{x^{p-1}}.$$ Then, we easily see that the centralizer $C$ of $g_1g_2\ldots g_p$ in $G$ is contained in $B$, so $$C=C_{H_1}(g_1)\times \ldots \times C_{H_p}(g_p)$$ and hence $$C^G=C_{H_1}(g_1)^{H_1}\times \ldots \times C_{H_p}(g_p)^{H_p}.$$ Since each $H_i$ has trivial pseudocentre, it follows that $P(G)$ is trivial as well.
\end{proof}

\medskip

In case of arbitrary posets $\Lambda$, the situation seems much more difficult to deal with, but we conjecture that the following result holds. 

\begin{conjecture}
Let $\mathbb F$ be a field, and $\Lambda$ a poset. Let $$u=1+a_{\alpha_1,\beta_1}e_{\alpha_1,\beta_1}+\ldots+a_{\alpha_m,\beta_m}e_{\alpha_m,\beta_m}$$ be an element of $M=M(\mathbb F,\Lambda)$ with $a_{\alpha_i,\beta_i}\neq0$ for  $i=1,\ldots,m$ and such that $(\alpha_i,\beta_i)\neq(\alpha_j,\beta_j)$ if $i\neq j$. Then $u\in P(M)$ if and only if for every $i=1,\ldots,m$ one of the following conditions holds: \begin{itemize}
    \item[\textnormal{(1)}] There are infinitely many pairs $(x_j,y_j)\in\Lambda^2$, $j\in\mathbb N$, with \hbox{$\alpha_i\leq x_j<y_j\leq\beta_i$} and $|\{x_j,y_j,x_k,y_k\}|=4$ for every $j\neq k$.
    \item[\textnormal{(2)}] There exist finitely many distinct elements $x_1,\ldots,x_t$ of $[\alpha_i,\beta_i]$ such that if $\alpha_i\leq y<z\leq\beta_i$, then~$y$ or~$z$ belongs to $\{x_1,\ldots,x_t\}$. If $t$ is the smallest possible such number of elements, then $\big|(\beta_i,+\infty)\big|+|(-\infty,\alpha_i)|< t$.
\end{itemize}

Moreover, $P(M)=M(\mathbb F,\Sigma)$, where $\Sigma$ is a poset having the same underlying set as $\Lambda$ but in which $\alpha\leq\beta$ if and only if $[\alpha,\beta]$ satisfies one of the above two conditions.
\end{conjecture}

\begin{remark}
    {\rm If $\Lambda$ is a linearly ordered set of finite cardinality $n$, then condition (2) in the above statement reads $$(n-\beta_i)+(\alpha_i-1)<\beta_i-\alpha_i.$$ This is equivalent to $2(\beta_i-\alpha_i)>n-1$ and hence to $\beta_i-\alpha_i\geq\lceil n/2\rceil$, which is precisely the condition in the statement of Theorem \ref{theoincrediblemclain}.}
\end{remark}

\medskip

Finally, we wish to observe that Mclain's examples can be easily employed to construct infinite locally nilpotent $p$-groups, $p$ a prime, such that $$G\simeq H\leq\underbrace{P(P(\ldots P}_{\textnormal{$n-1$ times}}(G))\ldots)\neq \underbrace{P(P(\ldots P}_{\textnormal{$n$ times}}(G))\ldots)$$ for any $n$ (choose for example as a poset an infinite grid in which every point in a horizontal line is connected to all the points above and below). Thus, the following question is very natural.

\begin{quest}
Let $n$ be a positive integer. Is there any finite nilpotent group~$G$ such that $$\underbrace{P(P(\ldots P}_{\textnormal{$n$ times}}(G))\ldots)\neq \underbrace{P(P(\ldots P}_{\textnormal{$n-1$ times}}(G))\ldots)?$$ 
\end{quest}

In this context, we would also like to note that it would be very interesting to give a description of the pseudocentre of an arbitrary finite $p$-groups.

\begin{quest}
Describe the pseudocentre of relevant families of finite $p$-groups, for any prime $p$.
\end{quest}

\chapter{General wreath products}\label{chapterwreathproducts}

Before studying wreath products in relation to the pseudocentre, we briefly wish to recall one of the most general definitions of wreath product for the reader's convenience. This is mainly taken from Section 6.2 of \cite{Rob72}. Let $\Lambda$ be a linearly ordered set. For each $\lambda\in\Lambda$, let~$H_\lambda$ be a transitive permutation group acting on the non-empty set $X_\lambda$, and fix an element $1_\lambda\in X_\lambda$. Let~$X$ be the {\it restricted} set product of all the~$X_\lambda$'s, that is, $X$ is the set of all  $(x_\lambda)_{\lambda\in\Lambda}$, where~\hbox{$x_\lambda\in X_\lambda$} and $x_\lambda\neq1_\lambda$ for finitely many $\lambda$'s. Put $(1_\lambda)=1$. Two elements $x$ and $y$ of $X$ are said to be~\hbox{\it$\lambda$-equi\-va\-lent} if $x_\mu=y_\mu$ for all $\mu>\lambda$; in this case we write $x\equiv_\lambda y$. Then $H_\lambda$ acts faithfully as a group of permutations of $X$ as follows. If $\xi\in H_\lambda$ and $x\in X$, then $(x)\xi=x$ if~\hbox{$x\not\equiv_\lambda1$,} otherwise $\xi$ only affects the $\lambda$-component of $x$ as $\xi$ would do. The {\it wreath product} $W=\operatorname{Wr}_{\lambda\in\Lambda}\,H_\lambda$ of the groups $H_\lambda$, $\lambda\in\Lambda$, is by definition the group of permutations of $X$ generated by the $H_\lambda$'s. Of course, $W$ depends on the groups $H_\lambda$ and on the way in which $H_\lambda$ acts on $X_\lambda$ but does not depend on the choice of the elements $1_\lambda$. Moreover, $W$ acts transitively on $X$. 

If $\Gamma$ and $\Delta$ are subsets of $\Lambda$ such that $\gamma<\delta$ for all $\gamma\in\Gamma$ and $\delta\in\Delta$, we  write $\Gamma\prec\Delta$. A non-empty subset~$\Gamma$ of $\Lambda$ is a {\it segment} of $\Lambda$ if for each $\lambda\in\Lambda$ exactly one of the following conditions holds: $\lambda\in\Gamma$, $\{\lambda\}\prec\Gamma$, and $\Gamma\prec\{\lambda\}$. If $\Lambda\prec\Lambda\setminus\Gamma$ ($\Lambda\setminus\Gamma\prec\Lambda$, resp.), then $\Gamma$ is said to be a {\it lower} segment (an {\it upper} segment, resp.). A collection $\{\Gamma_i\}_{i \in I}$ of mutually disjoint segments of $\Lambda$ is said to be a {\it segmentation of $\Lambda$} if $\Lambda = \bigcup_{i\in I}\Gamma_i$. If $W_i = \operatorname{Wr}_{\lambda\in\Gamma_i}H_\lambda$ for every~\hbox{$i\in I$,} then one can show that $W = \operatorname{Wr}_{i\in I}W_i$. If  $\{\Gamma_1, \Gamma_2\}$ is a segmentation of $\Lambda$, then we put $W = \operatorname{Wr}_{i=1,2}W_i=W_1 \operatorname{\it wr} W_2$. Moreover, if $\Lambda=\{1,2\}$ and $X_2=H_2$ with the action of $H_2$ on $X_2$ given by the product on the right, then $W$ is just the {\it standard wreath product} $H_1\wr H_2$ of $H_1$ and $H_2$.

\bigskip

\begin{theorem}\label{theowr}
Let $\Lambda$ be a linearly ordered with no maximum, $\{ H_\lambda : \lambda \in \Lambda\}$ a collection of non-trivial groups and $W = \underset{\lambda \in \Lambda}{\operatorname{Wr}} H_{\lambda}$. Then $W$ is pseudocentral.
\end{theorem}
\begin{proof}
Let $x,g\in W$. Then there is $\mu\in\Lambda$ such that $g,x\in U=\operatorname{Wr}_{\lambda<\mu} H_\lambda$. Let $B$ be the base group of $U\operatorname{\it wr} H_\mu$. Now, $B$ is a direct product of conjugates of~$U$ by elements of $H_\mu$. Since $C_W(x)$ contains all the factors of $B$ distinct from $U$, so $g\in C_W(x)^W$ by transitivity of the action. The arbitrariness of $g$ and $x$ yields that~$W$ is pseudocentral.
\end{proof}

\begin{corollary}\label{pseudocentralembedding}
Every group can be subnormally embedded in a pseudocentral group.
\end{corollary}

\begin{remark}
{\rm A dual result to Corollary \ref{pseudocentralembedding} is given in Corollary \ref{trivialembedding}.}
\end{remark}

\begin{corollary}
There exists a pseudocentral locally nilpotent group with trivial Baer radical.
\end{corollary}
\begin{proof}
Let $p$ be a prime, let $\Lambda$ be the set of all integers in their natural order, and for each $\lambda\in\Lambda$, let $H_\lambda$ be the cyclic group of order $p$. Then $W=\operatorname{Wr}_{\lambda\in\Lambda}H_\lambda$ is easily seen to have a trivial Baer radical and to be locally nilpotent. By~The\-o\-rem~\ref{theowr}, $W$ is pseudocentral.
\end{proof}

\medskip

In connection with the previous result, we state the following problem.

\begin{quest}
Is there a pseudocentral Baer group which is not a Fitting group? 
\end{quest}

\medskip

In the following, we deal with wreath products of finitely many groups, and we show how difficult and variegated the situation can be with respect to the computation of the pseudocentre. We start with the following simple remark: the pseudocentre of the standard wreath product~\hbox{$\mathbb Z_2\wr\mathbb Z_2$} of two cyclic groups of order~$2$ coincides with the centre, so it does not coincide with the direct product of the pseudocentres of the factors of the base group. Similar constructions can be used to produce examples of finite pseudonilpotent $p$-groups of arbitrarily high pseudonilpotent class.

\begin{theorem}\label{highpseudonilp}
Let $p$ be a prime. There exists finite $p$-groups of arbitrarily high pseudonilpotent class.
\end{theorem}
\begin{proof}
Let $G_1=\mathbb Z_p$, and assume $G_n$ has already been defined in such a way that $P_{n}(G_n)\neq G_n$. Put $G_{n+1}=\mathbb Z_p\wr G_n$. Clearly, $P=P(G_{n+1})$ is contained in the base group $B$, and $P_{n}(G_{n+1}/P)\neq G_{n+1}/P$ because $G_{n+1}/B\simeq G_n$.
\end{proof}

\begin{remark}
{\rm It is not clear which is the exact pseudonilpotent class of the previous examples. If we want finite groups of a given pseudonilpotent class, then we can start with $G_1$ being a finite simple non-abelian group in the proof of~The\-o\-rem~\ref{highpseudonilp}.}
\end{remark}

%


Now, we study the pseudocentre of wreath products of the forms $H\operatorname{\it wr}\operatorname{Alt}(n)$ (see Theorems~\ref{A_5Pseudo} and~\ref{a5easycases}), $H\operatorname{\it wr}\operatorname{Sym}(n)$ (see~The\-o\-rem~\ref{Sym}), and $H\wr K$, where~$K$ is infinite cyclic (see~The\-o\-rems~\ref{similarlamp} and~\ref{exsuperincredible}). In particular, Theorem \ref{exsuperincredible} provides an example of a group in which the most reasonable pseudocentral series have factors that are isomorphic to the whole group.

\begin{lemma}\label{wrLem}
Let $G = H \wr K$, where $H,K$ are non-trivial groups. If $B$ is the base group of $G$, then $B' \le [\langle kx\rangle, B]$ for every $1\neq k\in K$ and $x\in B$. In particular, if $K$ is infinite, then $B'\leq P(G)$.
\end{lemma}
\begin{proof}
Write $B=\operatorname{Dr}_{k\in K}H_k$, where $H_k\simeq H$ and the action of $K$ on the~$H_i$'s is the obvious one. If $a,b \in B$, then $b = b_{k_1} \ldots b_{k_t}$ and $a = a_{k_1} \ldots a_{k_t}$, where $a_{k_i},b_{k_i}\in H_{k_i}$ for all $i$, and $k_i\neq k_j$ if $i\neq j$. In order to prove the statement, it is sufficient to show that $[b_{k_i},a_{k_i}]\in [\langle kx \rangle,B]$ for every $1\neq k \in K$ and for all $k_i$'s. Thus, without loss of generality, we may assume $a,b \in H_1$.

Set $y=kx$ and let $c\in H_{k^{-1}}$ with $c^y = a$. Then
$$
\begin{array}{c}
[a,b]=[c^y,b]=[y^{-1},b]^{cy}[cy,b]=[y^{-1}[c,y],b^{cy}][y,b]\\[0.2cm]
=[y^{-1},b^{cy}]^{[c,y]}[[c,y],b^{cy}][y,b]\in [\langle y\rangle,B]
\end{array}
$$
and the statement is proved.
\end{proof}

\begin{remark}
{\rm Example \ref{exampleminp} shows that one may have standard wreath products $G=H\wr K$ in which $K$ is finite, $H$ is non-abelian, and $P(G)=\{1\}$.}
\end{remark}

\begin{lemma}\label{quotientwr}   
Let $G = H\wr K$, where $H,K$ are non-trivial groups. Let $B$ be the base group of $G$, and $X$ a non-trivial normal subgroup of $K$. If $\mathcal T$ is a left transversal of $X$ in $K$, then $$G/[X,B] \simeq K\ltimes\underset{t\in \mathcal T}{\operatorname{Dr}} (H_t/H_t'),$$ where $H_t\simeq H$ for each $t\in \mathcal T$, and $K$ acts by naturally permuting the factors of the direct product by subscript multiplication on the right.
\end{lemma}
\begin{proof}
Write $B=\operatorname{Dr}_{k\in K}H_k$, where $H_k\simeq H$ and the action of $K$ on the $H_i$'s is the obvious one. For each $b\in B$ and $k\in K$, let $\varphi_k(b)$ be the projection of $b$ onto~$H_k$. Now, for each $t\in\mathcal T$, let $$\psi_t(b)=\prod_{x\in X}\varphi_{tx}(b)^{x^{-1}}H_t'.$$ Moreover, if $g=kb\in G$, where $k\in K$ and $b\in B$, then put $\psi_{t}(g)=\psi_t(b)$ and $\xi(g)=k$ for every $t\in\mathcal T$. Finally, define $$\varphi: g\in G \mapsto \big(\xi(g),\big(\psi_t(g)\big)_{t\in\mathcal T}\big)\in K\ltimes\underset{t\in\mathcal T}{\operatorname{Dr}}(H_t/H_t’),$$ where the action is the one described in the statement. It is easy to see that $\varphi$ is an epimorphism. Suppose that $\varphi(g)=1$. Write $g=kb$, where $k\in K$ and $b\in B$. Then $k=1$ and $\psi_t(b)=H_t'$ for every $t \in \mathcal{T}$. Since $$\psi_t(b) = \prod_{x\in X}\varphi_{tx}(b)^{x^{-1}}H_t' = \prod_{x \in X}\varphi_{tx}(b)\prod_{x \in X}[\varphi_{tx}(b),x^{-1}]H_t',$$ we have that $\prod_{x \in X}\varphi_{tx}(b) \in H_t'[X,B]\leq [X,B]$ by Lemma \ref{wrLem}. The arbitrariness of $t \in \mathcal{T}$ yields that $b\in [X,B]$ and so that $g\in [X,B]$. Conversely, $$\psi_t(b^y) = \prod_{x \in X}\varphi_{tx}(b^y)^{x^{-1}} H_t' = \prod_{x \in X}\varphi_{txy^{-1}}(b)^{yx^{-1}}H_t'=\psi_t(b),$$ where $y\in X$ and $b\in B$, so  $\varphi([b,y]) =\varphi(b)^{-1}\varphi(b^y)=(1,H_t')$. Therefore, $[X,B]=\operatorname{Ker}(\varphi)$ and the statement is proved.
\end{proof}

\medskip

Let $\Omega$ be a non-empty set and $K$ a transitive group of permutations on $\Omega$. Let $x\in K$. If~\hbox{$\alpha\in \Omega$,} then we denote by $[\alpha]_{x}$ the \textit{orbit} of $\alpha$ under the action of $\langle x\rangle$, and with $\mathcal R_x$ a set of representatives for the orbits under the action of $\langle x\rangle$. Now, let $H$ be a non-trivial group, put~\hbox{$G = H\,\operatorname{\textit{wr}}\, K$}, and let $B$ be the base group of~$G$. If we put~\hbox{$H_{[\alpha]_x} = \operatorname{Dr}_{\beta\in [\alpha]_{x}} H_{\beta}$}, where $H_\beta \simeq H$, then it is easy to see that $H_{[\beta]_x} = H_{[\alpha]_x}$ if and only if $[\alpha]_x=[\beta]_x$. Clearly, $H_{[\alpha]_x}$ is the base group of $H_{\alpha}\wr\langle x\rangle$ and $B = \operatorname{Dr}_{\alpha\in\mathcal{R}_x} H_{[\alpha]_x}$. 
\medskip

\begin{corollary}
Let $K$ be a transitive permutation group on the non-empty set $\Omega$,  $H$ a non-trivial group, and $B$ the base group of $G = H\,\operatorname{\textit{wr}}\, K$. If $x\in K$, then $$B/[\langle x \rangle, B] \simeq \underset{\substack{\alpha\in\mathcal{R}_{x}\\ |[\alpha]_x|\neq1}}{\operatorname{Dr}} H_\alpha/H_\alpha'\times \underset{\substack{\alpha\in\mathcal{R}_{x}\\ |[\alpha]_x|=1}}{\operatorname{Dr}} H_\alpha,$$  where $H \simeq H_\alpha$ for every $\alpha\in \Omega$.
\end{corollary}

\begin{lemma}\label{A_5wr}
Let $p$ be a prime, $K$ a finite group of order $n$, and $H$ a cyclic group of order $p$. Set~\hbox{$G=H\wr K$,} and let $B$ be the base group of $G$.
    \begin{itemize}
        \item[\textnormal{(1)}] If $p$ does not divide $n$, then $B = [K,B] \times Z(G)$.
        \item[\textnormal{(2)}] If $p$ divides $n$, then $Z(G) \le [K,B]$.
    \end{itemize}
\end{lemma}
\begin{proof}
Write $B=\operatorname{Dr}_{k\in K}H_k$, where $H_k\simeq H$ and the action of $K$ on the~$H_i$'s is the obvious one. Consider the mappings $\varphi_1$, $\varphi$ and $\psi_1$ as in Lemma \ref{quotientwr}, where $X=K$ and $\mathcal T=\{1\}$. Now, if $x\in Z(G)$, then $x$ is contained in $B$ and \hbox{$\psi_1(x)=\varphi_1(x)^n$.} Thus, if $p$ divides $n$, then $x$ belongs to $\operatorname{Ker}(\psi_1)=\operatorname{Ker}(\varphi)=[K,B]$. If $p$ does not divide $n$, then $\psi_1(x)\neq1$, so $Z(G)\cap[K,B]=\{1\}$; also, order considerations yield that $B=Z(G)\times [K,B]$. The statement is proved.
\end{proof}

\begin{corollary}\label{A_5wrlevelpro}
Let $K$ be a transitive permutation group on the non-empty set $\Omega$,  $H$ a cyclic group of order $p$, and $B$ the base group of $G = H\,\operatorname{\textit{wr}}\, K$. If $x\in K$, then $[\langle x\rangle,B]C_B(x)$ is a proper subgroup of $B$ if and only if $p$ divides the length of some of the cycles in the disjoint cycle decomposition of $x$. Moreover, if this is not the case, then $B=[\langle x\rangle,B]\times C_B(x)$.
\end{corollary}
\begin{proof}
Let $\alpha\in\Omega$, and let $n$ be the order of $[\alpha]_x$ if $x$ has finite order and~$\infty$ otherwise. If~\hbox{$n=\infty$,} then $p$ divides $\infty$, while $H_\alpha\not\leq\langle[\langle x\rangle,H_\alpha], C_B(x)\rangle$. On the other hand, if $x$ is periodic, then the result follows at once from the fact that $\langle H_{\alpha},x\rangle/\langle x^n\rangle\simeq H_\alpha\wr\big(\langle x\rangle/\langle x^n\rangle\big)$ by applying~Lem\-ma~\ref{A_5wr}.
\end{proof}

\medskip

The next results will deal with the cases of  non-standard wreath products $H_1\operatorname{\it wr}H_2$, where either $H_2\simeq\operatorname{Alt}(n)$ or $H_2\simeq\operatorname{Sym}(n)$ acting on $X_2=\{1,\ldots,n\}$  in the natural way. Here, when speaking of {\it a cycle of an element} $x$ of $\operatorname{Sym}(n)$, we refer to a cycle in the disjoint cycle decomposition of $x$.

\begin{lemma}\label{centrA_5}
 Let $n\geq3$ be an integer, $p$ a prime, and $G=H\operatorname{\it wr}K$, where $H$ is cyclic of order $p$, and~\hbox{$K\simeq\operatorname{Alt}(n)$.} Then the following conditions are equivalent:
 \begin{itemize}
     \item[\textnormal{(1)}] $C_K(b)\neq\{1\}$ for every $b$ in the base group $B$ of $G$.
     \item[\textnormal{(2)}] $p+1 < n$.
 \end{itemize}
\end{lemma}
\begin{proof}
Write $B=H_1\times\ldots\times H_n$, where $H_i=\langle b_i\rangle\simeq H$ for every $i$. If $n=3$, then (2) is never satisfied and $C_G(b_1b_2)=B$. Assume $n\geq4$.

\medskip

\noindent(1)$\implies$(2)\quad Suppose $p+1\geq n$. Then $b_1b_2^2\ldots b_{n-2}^{n-2}$ has trivial centralizer in $K$.

\medskip

\noindent(2)$\implies$(1)\quad  Let $b \in B$. Then $b = b_{1}^{m_1}\ldots b_n^{m_n}$ where $0 \le m_j<p$. Define $\equiv_b$ over $\{1,\ldots,n\}$ with $$i\equiv_b j\iff m_i = m_j.$$ The number of equivalence classes of $\equiv_b$ is at most $p$. Let $k$ be the maximum of the orders of the equivalence classes of $\equiv_b$. If~$k$ is odd, then there is a $k$-cycle in $K$ centralizing $b$. If~$k$ is even and $k \ge 4$, then there is a $(2,2)$-cycle in $K$ centralizing~$b$. Suppose $k = 2$. Since~\hbox{$p+1<n$,} there exist at least two equivalence classes of order~$2$. Again, $K$ has a~\hbox{$(2,2)$-cycle} centralizing~$b$.~\end{proof}

\begin{theorem}\label{A_5Pseudo}
Let $n\geq5$ be an integer, $p$ a prime, and $G=H\operatorname{\it wr}K$, where~$H$ is cyclic of order $p$, and~\hbox{$K\simeq\operatorname{Alt}(n)$.}  Let $B$ be the base group of $G$.
    \begin{itemize}
       \item[\textnormal{(1)}] If $p +1 < n$ and $p$ does not divide $n$, then $P(G) = G$.
        \item[\textnormal{(2)}] If $p > n$ or $p +1 = n$, then $P(G)=B$.
        \item[\textnormal{(3)}] If $p+1 < n$ and $p$ divides $n$, then $P(G) = K\ltimes [K,B]$.
        \item[\textnormal{(4)}] If $p=n$, then $P(G)=[K,B]$.
    \end{itemize}
\end{theorem}
\begin{proof}
Before proving the four points of the statement, we need some preliminary remarks on the wreath product $G=H\operatorname{\it wr}K$. Write $B=H_1\times\ldots\times H_n$, where $H_i=\langle b_i\rangle\simeq H$, and~$K$ acts as $\operatorname{Alt}(n)$ on the subscripts of the $H_i$'s. It is easy to see that $$[K,B]=\big\{b_1^{m_1}\ldots b_n^{m_n}\in B\,:\, b_1^{m_1+\ldots +m_n}=1\big\}.$$ Thus, $B/[K,B]\simeq H$. In particular, the only normal subgroups of $G$ containing $[K,B]$ are $[K,B]$, $K[K,B]$, $B$ and $G$. Now, if $x\in K$ and $b\in B$, then $C_G(xb)^G$ contains $[\langle xb\rangle^G,B]=[K,B]$ and so even $K[K,B]$.

\medskip

\noindent(1)\quad Let $x \in K$ and $b \in B$. By Lemma \ref{centrA_5}, $C_G(b)^G=G$. Now, since $p$ does not divide~$n$, we may assume that there is $i\in\{1,\ldots,n\}$ whose orbit $[i]_x$ has order not divided by $p$. It follows from~Lem\-ma~\ref{A_5wr} applied to $U=\langle H_i,xb\rangle$ that $C_U(xb)^U\geq H_i$, so $C_G(xb)^G\geq B$ and consequently~\hbox{$C_G(xb)^G=G$.} Thus, $P(G)=G$.

\medskip

\noindent(2)\quad  By Lemma \ref{centrA_5}, $P(G) \le B$. Let $b \in B$ and $x \in K$. Since we have also in this case that $p$ does not divide $n$, as in case (1), an application of~Lem\-ma~\ref{A_5wr} yields that $C_G(xb)^G \ge B$. This shows that $P(G)=B$.

\medskip

\noindent(3)\quad By Lemma \ref{centrA_5}, $C_G(b)^G = G$ for every $b \in B$. By hypothesis, $n=pm$ for some positive integer~$m$. Then $K$ contains the element $x$ which is a product of~$m$ disjoint $p$-cycles. It easily follows from Corollary \ref{A_5wrlevelpro} that $C_G(x)^G\neq G$, so~\hbox{$C_G(x)=K[K,B]$.} Consequently, \hbox{$P(G)=K[K,B]$.}

\medskip

\noindent(4)\quad By Lemma \ref{centrA_5}, $P(G) \le B$. Since $n \ge 5$, so $n$ is odd and hence $K$ contains all the $p$-cycles. Let $x$ be one of these. By Corollary \ref{A_5wrlevelpro}, $C_G(x)^G = K[K,B]$. Thus,~\hbox{$P(G) = [K,B]$.}
\end{proof}

\begin{corollary}
The class of pseudocentral groups is not closed with respect to forming subdirect products. In particular, this class is not a formation.
\end{corollary}
\begin{proof}
Let $G=H\ltimes (K_1\times K_2)$, where $H=\operatorname{Alt}(5)$, $K_1\simeq K_2\simeq \mathbb Z_3^5$ and {$\langle H, K_1\rangle\simeq\langle H,K_2\rangle\simeq \mathbb Z_3\operatorname{\it wr}\operatorname{Alt}(5)$.} By Theorem \ref{A_5Pseudo}, $G/K_i$ is pseudocentral, while clearly $K_1\cap K_2=\{1\}$. Write $K_i=\langle a_{i,1}\rangle\times\ldots\times\langle a_{i,5}\rangle$, where $i\in\{1,2\}$, so that $H$ naturally acts permutationally on the $a_{i,j}$'s. Set $a_1=a_{1,1}a_{1,2}a_{1,3}^2a_{1,4}^2$ and $a_2=a_{2,2}a_{2,3}a_{2,4}^2a_{2,5}^2$. Then $C_H(a_1)=\langle (12)(34)\rangle$ and $C_H(a_2)=\langle (23)(45)\rangle$. Therefore $C_H(a_1a_2)=\{1\}$, so $P(G)\leq K_1\times K_2$, and hence $G$ is not pseudocentral.
\end{proof}

\begin{theorem}\label{a5easycases}
Let $n\in\{3,4\}$, $p$ a prime, and $G=H\operatorname{\it wr}K$, where $H$ is cyclic of order $p$, and~\hbox{$K\simeq\operatorname{Alt}(n)$.}  Let $B$ be the base group of $G$. 
   \begin{itemize}
       \item[\textnormal{(1)}] If $n = 3$ and $p\neq 3$, then $P(G)=B$.
       \item[\textnormal{(2)}] If $n=3$ and $p=3$, then $P(G)=[K,B]$.
       \item[\textnormal{(3)}] If $n=4$ and $p=2$, then $P(G)=V_4 \ltimes[K,B]$. 
       \item[\textnormal{(4)}] If $n=4$ and $p\geq3$, then $P(G)=B$.
   \end{itemize}
\end{theorem}
\begin{proof}
If $n = 3$, or $n=4$ and $p\geq3$, then Lemma \ref{centrA_5} yields  $P(G) \le B$. The argument in cases (1), (2) and (4) is similar to the one employed in the proof of~The\-o\-rem~\ref{A_5Pseudo}.

Assume $n =4$ and $p=2$. Let $x\in K$ and $b \in B$. If $x$ has order $3$, then $C_G(xb)^G \ge K[K,B]$. Also note that $[V_4,B]=[K,B]$ where $V_4$ is the~Sy\-low~\hbox{$2$-sub}\-group of $K$, so if $x$ has order $2$, then $C_G(xb)^G \ge V_4[K,B]$ and $C_G(x)^G=V_4[K,B]$. By Lemma \ref{centrA_5}, $C_G(b)^G\geq V_4$. Thus, $P(G) = V_4[K,B]$. 
\end{proof}

\begin{remark}
{\rm One can replace the cyclic group $H$ of order $p$ in Theorems \ref{A_5Pseudo} and \ref{a5easycases} by some other group. For example, if we choose a divisible abelian group, then the base group is always contained in the pseudocentre (see Corollary \ref{divisibleradical}) and one can easily extend the previous results.}
\end{remark}

Now, we move on to the study of wreath products of the form $H\,\operatorname{\it wr}\,\operatorname{Sym}(n)$.

\begin{lemma}\label{centrsym}
    Let $n \ge 3$ be an integer, $p$ a prime, and $G = H \operatorname{\textit{wr}} K$, where $K\simeq\operatorname{Sym}(n)$. Then the following conditions are equivalent: 
    \begin{itemize}
        \item[\textnormal{(1)}] $C_{K}(b) \neq \{1\}$ for every $b$ in the base group $B$ of $G$.
        \item[\textnormal{(2)}] $p < n$.  
    \end{itemize}
    More precisely, if $p <n$, then for each $b \in B$ there exists a transposition commuting with $b$.
\end{lemma}
\begin{proof}
The proof is essentially the same as that of Lemma \ref{centrA_5}.  
\end{proof}

\begin{lemma}[see \cite{Suz}, Chapter II, pp.295--296]\label{suzukicentralizer}
    Let $n \ge 3$ be an integer and $x\in\operatorname{Sym}(n)$. 
    \begin{itemize}
        \item[\textnormal{(1)}] If $x$ is an $n$-cycle, then $\langle x\rangle$ is self-centralizing.
        \item[\textnormal{(2)}] If $x$ has empty stabilizer and $x = x_1\ldots x_k$, where the $x_i$'s are product of cycles of the same order~$n_i$, and $n_i\neq n_j$ if $i\neq j$, then $$C_{\operatorname{Sym}(n)}(x) \simeq C_{\operatorname{Sym}(\Omega_1)}(x_1) \times \ldots \times C_{\operatorname{Sym}(\Omega_k)}(x_k),$$ where $\Omega_i$ is the support of $x_i$.
        \item[\textnormal{(3)}] If $x$ has empty stabilizer and $x = x_1\ldots x_k$, where the $x_i$'s are pairwise disjoint $m$-cycles, then $C_{\operatorname{Sym}(n)}(x)\simeq C_m \,\operatorname{\textit{wr}}\,\operatorname{Sym}(k)$. 
        \item[\textnormal{(4)}] If $\Omega$ is the support of $x$, then $C_{\operatorname{Sym}(n)}(x) \simeq C_{\operatorname{Sym}(\Omega)}(x) \times\operatorname{Sym}(\overline{\Omega})$ where $\overline{\Omega} = \{1,\ldots,n\} \setminus \Omega$. 
    \end{itemize}
\end{lemma}
\begin{proof}
(1)\quad Write $x = (a_1,\ldots,a_n)$, and let $g\in C_{\operatorname{Sym}(n)}(x)$. Since $$x^g = (g(a_1),\ldots,g(a_n)) = (a_1,\ldots,a_n),$$ we have that~\hbox{$g(a_1) = a_k$} for some $k = 1,\ldots,n$. Consequently, $g(a_i)=a_{\varepsilon(k+i-1)}$, where $\varepsilon(k+i-1)$ is the integer $k+i-1$ reduced modulo $n$, and hence $g=x^k$.

\medskip

\noindent(2)\quad Let $g \in C_{\operatorname{Sym}(n)}(x)$, so in particular $x_1^g \ldots x_k^g =x^g= x$. Since the decomposition of $x$ is unique, so $x_i^g = x_i$ for $i = 1,\ldots,k$. If $g_i$ is the restriction of $g$ to the support $\Omega_i$ of~$x_i$, then $g_i \in\operatorname{Sym}(\Omega_i)$ and the result follows,

\medskip

\noindent(3)\quad Let $B$ be the base group of $G=C_m\,\operatorname{\it wr}\,\operatorname{Sym}(k)$. If $g=ab\in G$, where \hbox{$a\in\operatorname{Sym}(k)$} and~\hbox{$b=(b_1,\ldots,b_k)\in B$,} then we can map $g$ to the element of $C_{\operatorname{Sym}(n)}(x)$ that first permutes the subscripts of the elements $x_i$ as prescribed by the permutation $a$, and then acts as the element $x_1^{b_1}\ldots x_k^{b_k}$. It is easy to see that such an assignation defines an isomorphism.

\medskip

\noindent(4)\quad This is obvious.
\end{proof}

\begin{lemma}\label{lemsymalt}
Let $n$ be a positive integer, $p$ a prime, and $G = \mathbb Z_p \,\operatorname{\it wr}\,\operatorname{Sym}(n)$. If $n \ge 3$ and $B$ is the base group of $G$, then $[\operatorname{Sym}(n),B] = [\operatorname{Alt}(n),B]$. Moreover, if $n=4$, then $[\operatorname{Alt}(n),B]=[V_4,B]$. 
\end{lemma}
\begin{proof}
Write $B=B_1\times\ldots\times B_n$, where $B_i=\langle b_i\rangle\simeq\mathbb Z_p$. If $b=b_1^{m_1}\ldots b_n^{m_n}\in B$ and $x\in\operatorname{Sym}(n)$, then $$[x,b]=[x,b_1]^{m_1}\ldots [x,b_n]^{m_n},$$ so it is enough to show that $[x,b_i]\in [\operatorname{Alt}(n),B]$ for every $i=1,\ldots,n$. Now, since $n\geq3$, there is $\ell\in\{1,\ldots,n\}$ such that $y=\big(\ell,i,x(i)\big)\in\operatorname{Alt}(n)$, and hence $[y,b_i]=[x,b_i]$. If $n=4$, then instead of $\big(\ell,i,x(i)\big)$, we can use the $(2,2)$-cycle $(\ell,\mu)\big(i,x(i)\big)$.
\end{proof}

\begin{theorem}\label{Sym}
Let $n \ge 3$ be an integer, $p$ a prime and $G = \mathbb Z_p\,\operatorname{\textit{wr}}\,\operatorname{Sym}(n)$. Moreover, let $B$ be the base group of $G$.
    \begin{itemize}
        \item[\textnormal{(1)}] If $p < n$ and $p$ does not divide $n$, then $P(G) =\operatorname{Alt}(n) \ltimes B$.
        
        \item[\textnormal{(2)}] If $p < n$ and $p$ divides $n$, then $P(G) = \operatorname{Alt}(n) \ltimes [\operatorname{Alt}(n),B]$.

        \item[\textnormal{(3)}] If $p = n$, then $P(G) = [\operatorname{Alt}(n),B]$.

        \item[\textnormal{(4)}] If $p > n$, then $P(G) = B$.  
    \end{itemize}
\end{theorem}
\begin{proof}
Before proving the four points of the statement, we observe that $$C_G(xb)^G\geq[\langle xb\rangle^G,B]\geq[\operatorname{Alt}(n),B]$$ for every $1\neq x\in\operatorname{Sym}(n)$ and $b\in B$ (see Lemma \ref{lemsymalt}). Note also that $$[\operatorname{Sym}(n),B]=[\operatorname{Alt}(n),B]$$ by~Lem\-ma~\ref{lemsymalt}, so $G/[\operatorname{Alt}(n),B]\simeq \operatorname{Sym}(n)\times\mathbb Z_p$. Write $B\!=\!H_1\times\ldots\times H_n$, where \hbox{$H_i\!=\!\langle b_i\rangle\simeq\mathbb Z_p$.} 

\medskip

\noindent(1)\quad Suppose first $n\neq 4$. Let $x\in \operatorname{Sym}(n)\setminus\{1\}$ and $b \in B$. By Lemma \ref{centrsym}, \hbox{$C_G(b)^G=G$.} Moreover, since $p$ does not divide $n$, we may assume that there is~\hbox{$i\in\{1,\ldots,n\}$} whose orbit~$[i]_x$ has order not divided by~$p$. It follows from~Lem\-ma~\ref{A_5wr} applied to $U=\langle H_i,xb\rangle$ that \hbox{$C_U(xb)^U\geq H_i$,} so~\hbox{$C_G(xb)^G\geq B$} and consequently $$C_G(xb)^G\geq\operatorname{Alt}(n)B.$$ Finally, put \hbox{$y=(1,\ldots,n)$} for $n$ odd, and $y=(1,\ldots,n-1)$ for $n$ even. Then $C_G(y)^G\leq\operatorname{Alt}(n) B$ and the statement is proved.

Now, assume that $n=4$, so $p=3$. Clearly, $$C_G(g)^G\geq V_4\ltimes[\operatorname{Alt}(4),B]$$ for every $g\in G\setminus B$. Moreover, Lemma \ref{centrsym} yields that $C_G(b)^G=G$ for every~\hbox{$b\in B$.}  Now, if $b\in B$ and $x$ is any element of $\operatorname{Sym}(4)$ that is not a~\hbox{$3$-cycle,} then \hbox{$C_G(xb)^G\geq B$} by Lemma \ref{A_5wr}. But of course if $x$ is a~\hbox{$3$-cycle,} then $$B\leq C_G(xb)^G\leq\operatorname{Alt}(4)\ltimes B,$$ so in any case $C_G(g)^G\geq B$ for $g\in G\setminus B$. Finally, we need to show that $$C_G(xb)^G\geq \operatorname{Alt}(4)\ltimes B$$ for every $1\neq x\in V_4$ and $b\in B$. This is certainly true if $b=1$. Assume $b\neq1$, and without loss of generality $x=(1,2)(3,4)$. If the projections of~$b$ on~$H_1$ and~$H_2$ are equal, then $(1,2)\in C_G(xb)$ and we are done.  If the projection of~$b$ on~$H_1$ is~$1$, while the projection on~$H_2$ is $b_2$ (resp., $b_2^{-1}$), then $(1,2)b_1^{-1}\in C_G(xb)$ (resp., $(1,2)b_1\in C_G(xb)$. If the projection of $b$ on~$H_1$ is $b_1$ (resp., $b_1^{-1}$) and that of $b$ on~$H_2$ is $b_2^{-1}$ (resp., $b_2$), then $(1,2)b_1b_2^{-1}\in C_G(xb)$ (resp., $(1,2)b_1^{-1}b_2\in C_G(xb)$). This essentially accounts for all cases and shows that $P(G)=\operatorname{Alt}(4)\ltimes B$.

\medskip

\noindent(2)\quad Suppose first that $p=2$ and $n=4$. Lemma \ref{centrsym} yields that $C_G(b)^G=G$ for every~\hbox{$b\in B$.} Let $x\in V_4$ and $b\in B$. We need to show that $C_G(xb)^G\geq\operatorname{Alt}(4)$. This is certainly true if there is a transposition commuting with $x$ and $b$. Assume not. Then without loss of generality we may assume  that the projection of~$b$ on $H_1$ is $1$, while the projection of~$b$ on $H_2$ is $b_2$. But in this case $(1,2)b_1$ centralizes~$xb$, and the claim is proved. Therefore $C_G(g)^G\geq\operatorname{Alt}(4)$ for every~\hbox{$g\in G\setminus B$.} Now, if~\hbox{$y=(1,2,3,4)$,} then $C_G(y)=C_{\operatorname{Sym}(4)}(y) C_B(y)$, so~$$C_G(y)^G\leq \operatorname{Sym}(4)\ltimes[\operatorname{Alt}(4),B]$$ by~Lem\-ma~\ref{A_5wr}. If $y=(1,2,3)$, then $C_G(y)^G\leq\operatorname{Alt}(4)\ltimes B$. It follows that $$P(G)=\operatorname{Alt}(4)\ltimes [\operatorname{Alt}(4),B].$$

Assume now that $n\neq 4$, so $n\geq5$ and consequently $$C_G(g)^G\geq\operatorname{Alt}(n)\ltimes[\operatorname{Alt}(n),B]$$ for any~\hbox{$g\in G\setminus B$.} Let~\hbox{$b\in B$.} By Lemma \ref{centrsym}, $C_G(b)^G=G$. Let $x=(1,\ldots,n)$. If $ya\in C_G(g)$, then the condition~\hbox{$x^{ya}=x$} gives $x^y=x$ and so $x^a=x$. Thus, $a\in Z(G)$. By Lemma \ref{A_5wr}, $Z(G)\leq [\operatorname{Alt}(n),B]$, and hence $$C_G(x)^G\leq\operatorname{Sym}(n)\ltimes[\operatorname{Alt}(n),B].$$ Finally, if $n$ is odd, then $C_G(x)^G\leq\operatorname{Alt}(n)\ltimes B$, while if $n$ is even, then $$C_G(y)^G\leq\operatorname{Alt}(n)\ltimes B,$$ where $y=(1,\ldots,n-1)$. In any case, it follows that~\hbox{$P(G)=\operatorname{Alt}(n)\ltimes[\operatorname{Alt}(n),B]$.}

\medskip

\noindent(3)\quad By Lemma \ref{centrsym}, $P(G) \le B$. Let $x=(1,\ldots,n)$. As in the previous point,~Lem\-ma~\ref{A_5wr} shows that $C_G(x)\leq\operatorname{Sym}(n)\ltimes [\operatorname{Alt}(n),B]$. Consequently, $P(G)=[\operatorname{Alt}(n),B]$.

\medskip

\noindent(4)\quad By Lemma \ref{centrsym}, $P(G) \le B$. As in case (1), since $p$ does not divide $n$, we have that $C_G(g)^G\geq B$ for every $g\in G\setminus B$. Consequently, $P(G)=B$.
\end{proof}

\begin{remark}
{\rm Let $p$ be any prime. The pseudocentre of $\mathbb Z_p\,\operatorname{wr}\,\operatorname{Sym}(2)$ is easily computed. If~\hbox{$p=2$,} then it coincides with the centre, while if $p>2$, then it is the base group.}
\end{remark}

\begin{remark}\label{cuborubik}
{\rm Using Theorem \ref{Sym} we can easily compute the pseudocentre of the symmetry group of the Rubik's cube obtained by disassembling and reassembling. In fact, this group is $$\mathbb Z_4^6\times\big(\mathbb Z_3\,{\it wr}\,\operatorname{Sym}(8)\big)\times \big(\mathbb Z_2\,{\it wr}\,\operatorname{Sym}(12)\big),$$ so its pseudocentre is the product of the pseudocentres of the factors (see Lem\-ma~\ref{directproducts}), which are all known by Theorem \ref{Sym}.

The Rubik's cube group, that is $$\big(\mathbb Z_3^7\times\mathbb Z_2^{11}\big)\rtimes\big(\big(\operatorname{Alt}(8)\times\operatorname{Alt}(12)\big)\rtimes\mathbb Z_2\big),$$ is a bit harder to tackle directly. Nevertheless, it is possible to compute it using~GAP as follows:

\smallskip

\noindent\texttt{gap> cube:=Group(\\
( 1, 3, 8, 6)( 2, 5, 7, 4)( 9,33,25,17)(10,34,26,18)(11,35,27,19),\\ ( 9,11,16,14)(10,13,15,12)( 1,17,41,40)( 4,20,44,37)( 6,22,46,35),\\ (17,19,24,22)(18,21,23,20)( 6,25,43,16)( 7,28,42,13)( 8,30,41,11),\\ (25,27,32,30)(26,29,31,28)( 3,38,43,19)( 5,36,45,21)( 8,33,48,24),\\ (33,35,40,38)(34,37,39,36)( 3, 9,46,32)( 2,12,47,29)( 1,14,48,27),\\ (41,43,48,46)(42,45,47,44)(14,22,30,38)(15,23,31,39)(16,24,32,40));\\
gap> l:=ConjugacyClasses(cube);\\
gap> R:=cube;\\
gap> for i in [1..81120] do\\
> c:=Representative(l[i]);\\
> d:=NormalClosure(cube,Centralizer(cube,c));\\
> R:=NormalIntersection(R,d);\\
> od;\\
gap> Size(R);
}

\smallskip

This analysis shows that the order of the pseudocentre of the Rubik's cube group is half of the order of the group itself, so the pseudocentre must coincide with the derived subgroup --- note that the normal subgroups of the Rubik's cube group have been computed in \cite{rubikcite}.}
\end{remark}

\medskip

In the remainder of this section, we compute the pseudocentre of certain infinite standard wreath products, among which the lamplighter group, and we construct some rather interesting examples of groups (see Theorem \ref{exsuperincredible}).

\begin{lemma}\label{lemexsuperincredible}
Let $K$ be an infinite group, $H$ a non-trivial perfect group, and set~\hbox{$G = H\wr K$.} Then \hbox{$P(G) = B$,} where $B$ is the base group of $G$.
\end{lemma}
\begin{proof}
Write $B = \operatorname{Dr}_{k \in K} H_k$, where $H_k\simeq H$ for every $k\in K$. If $b_1$ is any non-trivial element of $H_1$, then $C_G(b_1)^G = B$, and so $P(G) \le B$. On the other hand, by~Lem\-ma~\ref{wrLem}, $P(G) \ge B'=B$. Thus, $P(G) = B$.
\end{proof}

\begin{theorem}\label{exsuperincredible}
There exists a group $G$ satisfying the following properties: 
    \begin{itemize}
         \item[\textnormal{(i)}] $P(G) < G$ and $P(G) \simeq G$. 
        \item[\textnormal{(ii)}] $G/P(G) \simeq G$. 
    \end{itemize}

\noindent Moreover, $G$ can be chosen to be either torsion-free or locally finite.
\end{theorem}
\begin{proof}
Let $H$ be a countably infinite pseudocentral perfect group (for example a torsion-free or a locally finite simple group). Set $G_0 = H$, $G_1 = H \wr H$ and~\hbox{$G_{n+1} = G_n \wr G_n$} for~\hbox{$n>1$.} For each positive integer $n$, $G_n$ is an infinite perfect group. Now, define $L_n$ as the direct product of countably many copies of~$G_n$. It follows from~Lemma~\ref{lemexsuperincredible} that $P(L_{n+1})\simeq L_n$. Finally, set $G = \operatorname{Dr}_{n \in \mathbb{N}}L_n$. Then~\hbox{$P(G)\simeq G$} and~\hbox{$G/P(G) \simeq G$.} 
\end{proof}

\begin{remark}
{\rm With respect to the statement of Theorem \ref{exsuperincredible}, we remark that, among others, it is possible to construct the following other types of pathologies:
\begin{itemize}
    \item[(a)] groups $H$ such that the descending series $\{P(P(\ldots P(H))\ldots)\}$ stabilizes after finitely many steps but the upper pseudocentral series does not.
    \item[(b)] groups $H$ such that the upper pseudocentral series stabilizes at $H$ after finitely many steps but the descending series $\{P(P(\ldots P(H))\ldots)\}$ does not.
\end{itemize}

In fact, if $G$ and $G_0$ are the groups given in the proof of Theorem \ref{exsuperincredible}, then we can consider standard wreath products like $(\ldots ((G_0\wr G_0)\wr\ldots )\wr G_0)\wr G$ and~\hbox{$G\wr (G_0\wr (G_0\wr(\ldots (G_0\wr G_0)\ldots)))$.}}
\end{remark}

\medskip

Now, we deal with the pseudocentre of the lamplighter group and of similar other examples.

\begin{lemma}\label{commutaors}
Let $A$ be a non-trivial abelian group and $G = A \wr \langle x \rangle$ where $\langle x \rangle$ is infinite cyclic. Then, for each positive integer $n$, $[B,_n\langle x \rangle]/[B,_{n+1} \langle x \rangle] \simeq A$, where~$B$ is the base group of $G$. Moreover, $\bigcap_{n\in\mathbb N}[B,_n\langle x\rangle]=\{1\}$.
\end{lemma}
\begin{proof}
Write $B=\operatorname{Dr}_{n\in\mathbb Z}A_n$, where $A_n\simeq A$. Since $A$ is abelian, so $$[B,\langle x \rangle] = \underset{n \in \mathbb{Z}}{\operatorname{Dr}}D_n,$$ where $$D_n = \langle a_n^{-1}a_{n+1} : a_n \in A_n,a_{n+1} \in A_{n+1}, a_n^x=a_{n+1}\rangle,$$ and hence $\langle x,[B,\langle x\rangle]\rangle\simeq A\wr\langle x\rangle$. Thus, an application of Lemma \ref{quotientwr} shows that $$[B,\langle x\rangle]/[B,\langle x\rangle,\langle x\rangle]\simeq A.$$ The first half of the statement now follows by induction having observed that the base group of $\langle x,[B,\langle x\rangle]\rangle$ is~$[B,\langle x\rangle]$.

Finally, note that $[B,_m\langle x\rangle]$ is generated by  elements of the form $$(a)_{i_1}^{\pm}(a)_{i_{2}}^{n_2}\ldots (a)_{i_{m}}^{n_m}(a)_{i_{m+1}}^{\pm},$$ where $(a)_{i_j}$ is the image of $a$ in $A_{i_j}$, $i_1<i_2<\ldots<i_{m+1}$, and $n_2,\ldots,n_m\in\mathbb Z$. Thus,~$[B,_m\langle x\rangle]$ does not contain elements of the form $b_{\ell_1}b_{\ell_2}\ldots b_{\ell_m}$, where $b_{\ell_j}\in A_{\ell_j}$. It follows that $$\bigcap_{n\in\mathbb N}[B,_n\langle x\rangle]=\{1\}$$ and the statement is proved.
\end{proof}

\medskip

\begin{remark}
{\rm In the statement of Lemma \ref{commutaors}, we cannot require that $\langle x\rangle$ be finite. In fact, let $G=\mathbb Z_2\wr \mathbb Z_3$ be the wreath product of a cyclic group of order $2$ and a cyclic group $\langle x\rangle$ of order~$3$. If $B$ denotes the base group of $G$, then~\hbox{$\langle x,[\langle x\rangle, B]\rangle\simeq\operatorname{Alt}(4)$} and $[B,\langle x\rangle]=[B,\langle x\rangle,\langle x\rangle]$. 

More generally, as a consequence of Lemma \ref{quotientwr} (and its proof), we have that if $\langle x\rangle$ is finite, then $G=A\wr\langle x\rangle$ is a homomorphic image of $H=A\wr \langle y\rangle$, where $\langle y\rangle$ is infinite, so \hbox{$[C,_n\langle y \rangle]/[C,_{n+1} \langle y \rangle] \simeq A$,} where $C$ is the base group of $H$, and hence $[B,_n\langle y \rangle]/[B,_{n+1} \langle y \rangle]$ is a homomorphic image of $A$, where $B$ is the base group of~$G$.}
\end{remark}

\begin{theorem}\label{similarlamp}
Let $H$ be a non-trivial abelian group, $K = \langle x \rangle$ and infinite cyclic group, and put $G=H\wr K$. Then $P(G)=\{1\}$.
\end{theorem}
\begin{proof}
Let $B = \operatorname{Dr}_{n \in \mathbb{Z}} H_n$, where $H_n \simeq H$, be the base group of $G$ --- here, $H_n^x=H_{n+1}$. Let $p$ be an odd prime, and set $g=g_p= x^pa_0a_{p}^{-1}$, where $1\neq a_0 \in H_1$ and~\hbox{$a_0^{x^i}=a_i$} for all $i\in\mathbb Z$.

We claim that $C_G(g) = \langle g \rangle$. Clearly, $C_B(g)=\{1\}$, so~$C_G(g)$ is cyclic. Suppose~$C_G(g)$ properly contains $\langle g \rangle$. Then there is a prime $q$ and $h\in G$ with $h^q=g$. Write $h=x^kb_{i_1} \ldots b_{i_t}$, where $k$ is a positive integer, $b_{i_j}\in H_{i_j}\setminus\{1\}$ for $j\in\{1,\ldots,t\}$, and $i_1<i_2<\ldots<i_t$. Since $h^q=g$, so $kq=p$, which means $q=p$ and~\hbox{$k=1$.} Moreover, the same equality also yields that $i_t=1$ and $b_{1}=a_1^{-1}$ because $b_{i_t}^{x^{p-1}}=a_p^{-1}$. On the other hand, $i_1$ must be equal to $0$ because any positive power of $h$ preserves the $i_1$-term; in particular, $b_0=a_0$ and $t=2$. This means that $x^{p-1}\in C_G(g)$, a contradiction. 

Finally, note that $\langle g_p\rangle^G\cap B\leq [x^p,B]$ and that $P(G)\leq B=C_G(a_0)^G$. Therefore $$P(G)\leq B\cap \bigcap_{p\in\mathbb P\setminus\{2\}}\langle g_p\rangle^G\leq \bigcap_{p\in\mathbb P\setminus\{2\}}[x^p,B]=\{1\}$$ by Lemma \ref{commutaors}, and the result is proved.  
\end{proof}

\begin{corollary}\label{perfectgroupincredible}
Let $K$ be an infinite cyclic group, $H$ a non-trivial group, and $G = H\wr K$. Then \hbox{$P(G) = B'$,} where $B$ is the base group of $G$.
\end{corollary}
\begin{proof}
Note first that the consideration of any non-trivial element of one of the factors of~$B$ yields that $P(G)\leq B$. Thus, it follows from Theorem \ref{similarlamp} that $P(G)\leq B'$. On the other hand,~Lem\-ma~\ref{wrLem} shows that $B'\leq P(G)$, and hence $P(G)=B'$.
\end{proof}

\begin{corollary}
The lamplighter group has trivial pseudocentre.
\end{corollary}

 \begin{corollary}
     There exists a residually nilpotent group with trivial pseudocentre. 
 \end{corollary}
 \begin{proof}
Let $p$ be a prime and $G =\mathbb Z_p \wr \langle x \rangle$, where $\langle x \rangle $ is infinite cyclic. By~Lem\-ma~\ref{commutaors}, $G$ is residually nilpotent, while Theorem \ref{similarlamp} yields that \hbox{$P(G)=\{1\}$.}
 \end{proof}

\chapter{Free products}\label{freeproductsect}

It is well-known that the free product of two groups of order $2$ is the infinite dihedral group, and it is easy to see that the pseudocentre of the infinite dihedral group is the derived subgroup. However, it turns out that this is the only exception. In fact, we are now going to prove that the pseudocentre of any free product of two groups $H$ and $K$, which do not have both order $2$, is trivial (see Theorem \ref{theoremafinalefreeproduct}). This is accomplished in two steps, according to whether or not the free product is the {\it modular group} (that is, $\operatorname{PSL}(2,\mathbb Z)=\mathbb Z_2\ast\mathbb Z_3$). The reason for this being as follows: in the ‘‘non-modular’’ case, $H$ and $K$ have enough elements to allow us to find commutators generating self-centralizing cyclic subgroups, while this seems to be not possible in the ‘‘modular case’’. As we shall see, in the case of $\operatorname{PSL}(2,\mathbb Z)$, we need to find self-centralizing cyclic subgroups of the congruence subgroups. This is accomplished in two different ways, each interesting in its own respect: one uses~Pell's equation and the other Chebyshev polynomials of first and second kind.

\begin{remark}
{\rm In Section \ref{appendice}, Anthony Genevois, using deep methods from geometric group theory, proves a more general result about acynlindrically hyperbolic groups. Differently from his proofs, ours are constructive, and give further insights on other types of groups (like the free soluble ones).}
\end{remark}

\begin{theorem}\label{freeproducttrivial}
Let $G=H\ast K$ be the free product of the non-trivial groups $H$ and $K$. If $\{|H|,|K|\}\neq\{2\},\{2,3\}$, then $P(G) = \{1\}$. 
\end{theorem}
\begin{proof}
Suppose $y$ is an element of $G$ that satisfies the following properties:
\begin{itemize}
    \item[(1)] there exist $u\in G$, $h\in H\setminus\{1\}$ and $k\in K\setminus\{1\}$ such that $y=[h,k]ux$, where $x=h_1k^{-1}hk$ with $1\neq h_1\in H$;
    \item[(2)] the length $\ell(u)$ of $u$ is $\geq4$; 
    \item[(3)] replacing the normal forms of $x$ and $u$ in $y=[h,k]ux$ gives the normal form of $y$;
    \item[(4)] in the first half (from the left) of the normal form of $y$, there is exactly one of the subwords $[h,k]$ and $x$.
\end{itemize}

We claim that $C_G(y)=\langle y\rangle$, or, in other words, that every element of $C_G(y)$ is a power of~$y$. Let $c$ be any element of $G$ (which we consider written in its normal form) such that $cy=yc$. We use induction on the length~$\ell(c)$ of~$c$. If either~$\ell(yc)$ or $\ell(cy)$ is strictly smaller than $\ell(c)$, then we may apply induction. Assume that this is not the case. Now, if $c$ does not completely cancel out when multiplied on the right of~$y$, then this does not happen also when $c$ is multiplied on the left of~$y$, so the first element of $c$ must be a non-trivial element of $H$, while the latter is a non-trivial element of $K$. But since~$c^{-1}$ centralizes $y$ too, the same argument yields that~$c^{-1}$ cancels out when multiplied on the left and on the right of~$y$. So, we may assume without loss of generality that~$c$ has the latter property. Therefore, $y=h^{-1}\ldots c^{-1}=c^{-1}\ldots k$, and consequently~\hbox{$c=k^{-1}\ldots h$.} Since $\ell(cy)\geq\ell(c)$, so $\ell(c)\leq\ell(y)/2$ and hence the normal form of~$y$ is given either by $c^{-1}c^{-1}$ or by~$c^{-1}z_3c^{-1}$ for a non-trivial word $z_3$.

If $\ell(c^{-1})\geq4$, then $c^{-1}$ starts with $[h,k]$ and ends with $x$ (although there may be some overlapping). Thus, if $\ell(c^{-1})>6$, then there is no overlapping, so $c^{-1}=[h,k]\ldots x$ and hence property (4) gives $c^{-1}=y$, and we are done. If $\ell(c^{-1})=6$, then $c^{-1}=[h,k]hk=h^{-1}khkhk$, $k=k^{-1}$, and $h_1=h$, contradicting property (4). Finally, if $4\leq \ell(c^{-1})< 6$, then $\ell(c^{-1})=4$, so~\hbox{$c^{-1}=[h,k]$} and hence $cy=yc$ yields $z_3=[h,k]\ldots$ contradicting (4) because of (2).

It is therefore possible to assume $\ell(c^{-1})<4$, so  $h=h^{-1}$, $k=k^{-1}$ and~$c=kh$. Since $c^{-1}$ centralizes $y$, conditions (1)--(3) yield $y=hkhkhk\ldots k$, contradicting~(4). The claim is proved.

\medskip

Now, let $U=[H,K]$, and note that $U$ is a free group on the generators $[h,k]$, \hbox{$h\in H\setminus\{1\}$} and $k\in K\setminus\{1\}$ (see, for example, \cite{acourse}, Exercise 7 of Section 6.2). We divide the proof into cases. In each of them, we are going to find infinitely many~\hbox{$G$-in}\-va\-riant subgroups $\{Y_i\}_{i\geq0}$ of $U$ whose intersection is trivial, and elements~\hbox{$y_i\in Y_i$} satisfying the above four properties (1)--(4). This shows that $$P(G)\leq \bigcap_{i\geq0} C_G(y_i)^G\leq\bigcap_{i\geq0}Y_i=\{1\}$$ and completes the proof. The most natural candidates for our subgroups $Y_i$ are terms of the lower central series of $U$, and terms of the form $[U,G,\ldots, G]$.

\medskip





Suppose first that $H$ has non-trivial elements $h$ and $h_1$ such that $h_1\not\in\{h,h^{-1}\}$, and that~$K$ has non-trivial elements $k$ and $k_1$ such that $k_1\not\in\{k,k^{-1}\}$. Define $$v_0=[k_1,h_1],\, v_1=\big[v_0^{-1},[k_1,h]\big],\, v_{i+1}=\big[v_i^{-1},[k_1,h_1]\big]$$ for every $i\geq1$. It is easy to see that $v_i=k_1^{-1}h_1^{-1}k_1h_1h^{-1}k_1^{-1}\ldots k_1^{-1}h^{-1}k_1h$ for every $i\geq1$. Moreover, the only elements of $K$ appearing in the normal form of~$v_i$ are $k_1^{\pm1}$. Now, put~\hbox{$u=[h,k]$,} choose $i\geq1$ and define $$y_i=[v_iu^{-1},v_i^{-1}u]=uv_i^{-1}u^{-1}v_iv_iu^{-1}v_i^{-1}u.$$ Then $y_i$ satisfies (1)--(4), and $y_i\in \gamma_{i+1}([H,K])$, so we are done.

\smallskip

This shows that we may assume without loss of generality that one between $H$ and $K$, say $K$, has order at most $3$.

\smallskip





Suppose now that we can choose non-trivial elements $h$ and $h_1$ of $H$ in such a way that $h_1\neq 1,h,h^{-1},h^2$. Let $k$ be a non-trivial element of $K$. Define $$v_0=[h_1,k],\; v_{i+1}=[v_i^{-1},[h,k]]$$ for every $i\geq0$. Then for every $i\geq1$ we can write $v_i=h_1^{-1}k^{-1}\ldots k(h_1 h^{-1})k^{-1}hk$, and note that the normal form of $v_i$ does not contain expressions of the form $[h,k]^{\pm1}$. Put~\hbox{$u=[h,k]$,} and choose $i\geq1$. Then $y_i=[v_i^{-1}u^{-1},v_iu]=uv_iu^{-1}v_i^{-1}v_i^{-1}u^{-1}v_iu$ satisfies properties (1)--(4), and belongs to $\gamma_{i+1}([H,K])$. Again, we are done. 

\smallskip

This shows that $H$ is cyclic of order either $3$ or $4$.

\smallskip






Now, assume that $H$ has order $4$ and that $|K|=2$. It follows from \cite{Grue},~Co\-rol\-la\-ry at~p.44, that $G$ is residually nilpotent, so the intersection of the terms $[[H,K],G,\ldots,G]$ is $\{1\}$. Let $h$ be the element of order $2$ of $H$, $h_1$ a generator of $H$, and $k$ be a non-trivial element of $K$. Define $v_0=[h_1,k]$, $v_1=[[k,h_1],k]$, $v_2=[h_1,[[k,h_1],k]]$. In general, if $v_i$ is defined for $i\geq1$, then we put $$v_{i+1}=\begin{cases}
[h_1,v_i] & \textnormal{if $i$ is odd}\\
[v_i^{-1},k] & \textnormal{if $i$ is even}
\end{cases}$$ It is easy to see that $v_i=h_1^{-1}k^{-1}\ldots h_1k$ for every $i\geq0$. Put $u=[h,k]$. Then $y_i=[v_i^{-1}u^{-1},v_iu]$, which we can write as $uwu$, satisfies (1)--(4) because there is only one occurrence of $h$ in the first half of $w$, namely in $v_iu^{-1}=\ldots k(h_1h)kh$. Clearly, $$y_i\in [[H,K],\underbrace{G,\ldots,G}_{\textnormal{$i$ times}}]$$ and so also this case is done.

\smallskip





Suppose now that both $H$ and $K$ have order $3$. Again, $G$ is residually nilpotent. Let $h,k$ be non-trivial elements of $H$ and $K$, respectively. Define $v_0=[kh^{-1}k,h^{-1}]$ and note that~$v_0$ does not contain any subwords of type $[h,k]^{\pm}$. Note also that $$v_0=k^{-1}hk^{-1}h\cdot kh^{-1}kh^{-1}\in [H,K].$$  Put $v_1=[v_0^{-1},k^{-1}]=k^{-1}hk^{-1}h\ldots k^{-1}h^{-1}kh^{-1}$. Assume we have defined by induction a word $$v_i=k^{-1}hk^{-1}h\ldots k^{-1}h^{-1}kh^{-1}\in \gamma_i([H,K],G)=[[H,K],\underbrace{G,\ldots,G}_{\textnormal{$i$ times}}]$$ that does not contain any subwords of type $[h,k]$ and $[k,h]$. Put $v_{i+1}=[v_i^{-1},k^{-1}]$. Then one easily sees that $$v_{i+1}=k^{-1}hk^{-1}h\ldots k^{-1}h^{-1}kh^{-1}\in\gamma_{i+1}([H,K],G)=[[H,K],\underbrace{G,\ldots,G}_{\textnormal{$i+1$ times}}]$$ does not contain any subwords of type $[h,k]^{\pm}$. Now, put $u=[h,k]$ and, for each \hbox{$i\geq0$,} consider the word $y_i=[v_iu^{-1},v_i^{-1}u]$ (which belongs to $\gamma_i([H,K],G)$). Observe that there is only one occurrence of $u$ in the first half (from the left) of the normal form of $y_i$. Thus, again properties (1)--(4) are satisfied and we are done.

\smallskip







Assume finally that $H\simeq\mathbb Z_4$ and $|K|=3$. Let $h$ be an element of order $4$ of $H$, and let $k$ be a non-trivial element of $K$. Define the words $v=[k,h^2]$ and $w=[k^{-1},h^2]$. Put $$v_1=[v,w]=h^2k^{-1}h^2k\ldots kh^2k^{-1}h^2$$ and $$v_2=[v_1^{-1},v]=h^2k^{-1}h^2k\ldots k^{-1}h^2kh^2.$$
Suppose that by recursion we have defined words $$v_{2i+1}=h^2k^{-1}h^2k\ldots kh^2k^{-1}h^2$$ and $$v_{2i+2}=h^2k^{-1}h^2k\ldots k^{-1}h^2kh^2$$ for $i\geq0$. Then we put $$v_{2i+3}=[v_{2i+2}^{-1},w]=h^2k^{-1}h^2k\ldots kh^2k^{-1}h^2$$ and $$v_{2i+4}=[v_{2i+3}^{-1},v]=h^2k^{-1}h^2k\ldots k^{-1}h^2kh^2.$$ Now, it is easy to see that for any odd $i$, the only occurrence of the subword $u=[h,k]$ in the word $y_i=[v_iu^{-1},v_i^{-1}u]\in\gamma_i([H,K])$ is at the beginning, while the only occurrence of the subword $hk^{-1}hk$ is at the end. Thus, again the conditions are satisfied, and we are done.
\end{proof}

\begin{corollary}\label{trivialembedding}
Every group can be embedded in a group with trivial pseudocentre.
\end{corollary}



\medskip

In order to prove that the pseudocentre of the modular group $\operatorname{PSL}(2,\mathbb Z)$ is trivial, we need to recall some facts and notation about the so-called {\it Pell's equations}, that is, any Diophantine equation of the form $x^{2}-dy^{2}=1$, where $d$ is a given positive non-square integer. The {\it fundamental solution} of these equations is a pair of positive integers $(x_1,y_1)$ with $x_1^2-dy_1^2=1$ and the value $u=x_1+y_1\sqrt{d}$ is smallest possible. Let $$u^n=x_n+y_n\sqrt d,\quad n=0,1,2,\ldots$$ It is well-known that $(\pm x_n,\pm y_n)$, $n=0,1,2,\ldots$, is a complete set of (integers) solutions to~Pell's equation $x^{2}-dy^{2}=1$.

\medskip

For our purposes, we will also need a matrix form for the solutions. First note that if $(x,y)$ is a solution with $x,y>0$, then there is $n>0$ with $(x,y)=(x_n,y_n)$ and $$
\begin{array}{c}
x_n+y_n\sqrt d=(x_1 + y_1\sqrt{d})^n = (x_1 + y_1\sqrt{d})^{n-1}(x_1 + y_1\sqrt{d})=\\[0.2cm]
(x_{n-1}+y_{n-1}\sqrt d)(x_1+y_1\sqrt d)= x_{1}x_{n-1} + y_{n-1}y_1d + (x_1y_{n-1} + y_1x_{n-1})\sqrt{d},
\end{array}
$$ which in turns gives $$\begin{pmatrix}
    x \\
    y
\end{pmatrix}= \begin{pmatrix}
    x_1 & dy_1 \\
    y_1 & x_1
\end{pmatrix}^{n} \begin{pmatrix}
    1 \\
    0
\end{pmatrix}$$ Assume $y<0<x$, and note the following easy fact.

\begin{itemize}
    \item[($\star$)] {\it If $(x,y)$ is a solution with $y<0<x$, then there exists a positive integer $n$ such that $x + y\sqrt{d} = (x_1 + y_1\sqrt{d})^{-n}$.}
\end{itemize}
\begin{proof}
Write $y = -z$, where $z$ is a positive integer. Then $x + z\sqrt{d} = (x_1 + y_1\sqrt{d})^n$ for some positive integer $n$. Then $(x_1 + y_1\sqrt{d})^{-n} = x-z\sqrt{d}=x+y\sqrt{d}$. 
\end{proof}

\medskip

Now, ($\star$) implies the existence of $n>0$ with $x+y\sqrt d=(x_1+y_1\sqrt d)^{-n}=(x_1-y_1\sqrt d)^n$. As above, this yields 
\begin{equation}
\begin{pmatrix}
    x \\
    y
\end{pmatrix}= \begin{pmatrix}
    x_1 & -dy_1 \\
    -y_1 & x_1
\end{pmatrix}^{n} \begin{pmatrix}
    1 \\
    0
\end{pmatrix}= \begin{pmatrix}
    x_1 & dy_1 \\
    y_1 & x_1
\end{pmatrix}^{-n} \begin{pmatrix}
    1 \\
    0
\end{pmatrix}
\end{equation} Thus, every solution $(x,y)$ with positive $x$ can be obtained from the matrix equation $$\begin{pmatrix}
    x \\
    y
\end{pmatrix}= \begin{pmatrix}
    x_1 & dy_1 \\
    y_1 & x_1
\end{pmatrix}^{n} \begin{pmatrix}
    1 \\
    0
\end{pmatrix},$$ where $n\in\mathbb Z$. Consequently, all the solutions $(x,y)$ with $x<0$ are given by $$\begin{pmatrix}
    x \\
    y
\end{pmatrix}= -\begin{pmatrix}
    x_1 & dy_1 \\
    y_1 & x_1
\end{pmatrix}^{n} \begin{pmatrix}
    1 \\
    0
\end{pmatrix},$$ where $n\in\mathbb Z$.

Recall also that if $n\geq1$ is an integer, then there is a homomorphism $$\pi_n:\operatorname{SL}_2(\mathbb Z)\rightarrow\operatorname{SL}_2(\mathbb Z_{n})$$ induced by the reduction modulo $n$. The {\it principal congruence subgroup of level $n$} of $\operatorname{SL}_2(\mathbb Z)$ is the kernel of $\pi_n$, and is usually denoted by $\Gamma(n)$. Clearly, $\Gamma(n)$ is the set of all matrices $x\in\operatorname{SL}(2,\mathbb Z)$ such that $x-1\equiv_n0$.

\begin{theorem}\label{pellequation}
The pseudecentre of $\operatorname{SL}(2,\mathbb{Z})$ coincides with its centre, while the psudocentre of $\operatorname{PSL}(2,\mathbb{Z})$ is trivial.
\end{theorem}
\begin{proof}
Let $n\ge 2$ be an integer, set $G=\operatorname{SL}(2,\mathbb Z)$, and define the aperiodic element $$A_n=A = \begin{pmatrix}
        1 & 2n \\
        n & 2n^2 +1 
    \end{pmatrix}$$ of $G$.
    Let $C_n=C= \begin{pmatrix}
        x_1 & x_2 \\
        x_3 & x_4 
    \end{pmatrix}\in C_G(A)$. Then $$\begin{cases}
         x_2 - 2x_3 = 0 \\
         x_1 + nx_2 - x_4 = 0
    \end{cases},$$ so $$
    C=\left\{
    \begin{pmatrix}
        a & 2b \\
        b & a + 2nb
    \end{pmatrix}\,:\, a^2 - 2b^2 + 2n ab - 1 = 0,\,\, a,b \in \mathbb{Z}
    \right\}.
    $$ Putting $x = a + nb$ and $y = b$, the equation $a^2 - 2b^2 + 2n ab - 1\!=\! 0$ becomes \hbox{$x^2 - (n^2 + 2)y^2\!=\! 1$.} 
    In order to solve this Pell's equation we use the so-called continuum fraction method to find the fundamental solution $(x_1,y_1)$. Thus, we write
    $$ \sqrt{n^2 + 2} = \lfloor\sqrt{n^2 + 2}\rfloor + (\sqrt{n^2 + 2} - \lfloor\sqrt{n^2 + 2}\rfloor) = n + (\sqrt{n^2 + 2} -n) = $$
    $$n + \frac{1}{\frac{\sqrt{n^2 + 2} +n}{2}} = n + \frac{1}{n + \frac{1}{\sqrt{n^2 + 2} + n}}= n + \frac{1}{n + \frac{1}{ 2n + \frac{2}{\sqrt{n^2 + 2} + n}}}.$$ Then $\sqrt{n^2+2}\sim \frac{n^2+1}{n}$, and so the fundamental solution is $x_1=n^2+1$ and $y_1=n$. Therefore, all the solutions of the equation $x^2 - (n^2 + 2)y^2 = 1$ are given by $$ \pm\begin{pmatrix}
        n^2 +1 & (n^2 +2)n \\
        n & n^2 +1
    \end{pmatrix}^m \cdot \begin{pmatrix}
        1 \\
        0
    \end{pmatrix},$$ where $m\in\mathbb Z$. Set $$D=\begin{pmatrix}
        n^2 +1 & (n^2 +2)n \\
        n & n^2 +1
    \end{pmatrix}\quad\textnormal{ and }\quad L = \begin{pmatrix}
    1 & n \\
    0 & 1
\end{pmatrix}.$$ Note that the equation $a^2 - 2b^2 + 2nab = 1$ can be written as \begin{equation}\label{PellSL}
    \begin{pmatrix}
    a & b
\end{pmatrix} \cdot L^t\begin{pmatrix} 1 & 0 \\
0 & -(n^2 +2)
    
\end{pmatrix}L \cdot \begin{pmatrix}
    a \\
    b
\end{pmatrix} = 1,\end{equation} where $L^t$ is the transpose of $L$. 

\medskip

We claim that all the solutions of $a^2 - 2b^2 + 2nab = 1$ are given by $$\pm L^{-1} D^m L\cdot\binom{1}{0}=\pm\big(L^{-1}DL\big)^m\cdot\binom{1}{0}=\pm\begin{pmatrix}
1 & 2n\\
n & 2n^2+1
\end{pmatrix}^m\cdot\binom{1}{0}=\pm A^m\cdot\binom{1}{0},$$ where~\hbox{$m\in\mathbb Z$.} Let $K=(k_1,k_2)$ be a solution of the equation $a^2 - 2b^2 + 2nab = 1$. Then $(LK^t)^t$ is a solution of~\hbox{$x^2 - (n^2 +2)y^2 = 1$} by (\ref{PellSL}). If the first component of $(LK^t)^t$ is positive, then there is an integer $\ell$ such that $$LK^t = D^\ell\cdot\begin{pmatrix}
    1 \\ 
    0
\end{pmatrix}.$$ Then $$K^t=L^{-1}\cdot D^\ell\cdot \begin{pmatrix}
    1 \\ 
    0
\end{pmatrix}=L^{-1}\cdot D^\ell\cdot L\cdot\begin{pmatrix}
    1 \\ 
    0
\end{pmatrix}$$ Similarly, we deal with the case in which the first component of $(LK^t)^t$ is negative. The claim is proved.


Set $$B = \begin{pmatrix}
    0 & 2 \\
    1 & 2n
\end{pmatrix}.$$ Then every element of $C$ can be written as $aI_2+bB$ for some integers $a,b$ satisfying the equation $a^2 - 2b^2 + 2nab = 1$. In particular, $A=I_2+nB$ and $A^{-1}=(2n^2+1)I_2-nB$. Put $$\begin{pmatrix}
    a_m \\
    b_m
\end{pmatrix} = A^m\cdot \begin{pmatrix}
    1 \\
    0
\end{pmatrix}$$ for every integer $m$. Clearly, $a_0=1$, $b_0=0$, $a_1=1$, $b_1=n$, $a_{-1}=2n^2+1$ and $b_{-1}=-n$.

We claim that $A^m=a_mI_2+b_mB$ for every $m\in\mathbb Z$. First note that $B^2=2I_2+2nB$. Thus, if~\hbox{$m\geq0$,} then $$
\begin{array}{c}
A^{m+1} = A^m\cdot A = (a_mI_2 + b_mB)(a_1I_2 + b_1B)\\[0.2cm]
= (a_ma_1 + 2b_mb_1)I_2 + (a_1b_m+a_mb_1 + 2nb_1b_m)B = a_{m+1}I_2+b_{m+1}B.
\end{array}
$$ The case $m<0$ can be dealt with in a similar way.  The claim is proved. 

The previous claim shows that all the elements of $C$ can be written as $\pm A^m$ for $m\in\mathbb Z$, so~\hbox{$C=Z(G)\times A$.}

\medskip

Now, $$\bigcap_{n\in\mathbb N}\Gamma(n)=\{1\},$$ so  $$\bigcap_{n\in\mathbb N}\big(Z(G)\Gamma(n)\big)=Z(G)$$ and consequently $$P(G)\leq \bigcap_{n\in\mathbb N}C_n=\bigcap_{n\in\mathbb N} \big(Z(G)\langle A_n\rangle^G\big)\leq\bigcap_{n\in\mathbb N} \big(Z(G)\Gamma(n)\big)=Z(G).$$ Therefore $P(G)=Z(G)$.

Finally, the centralizer $D_n/Z(G)$ of $A_nZ(G)/Z(G)$ in $G/Z(G)$ is infinite cyclic (see \cite{acourse}, p.172), so $D_n$ is abelian and is contained in $C_n=C_G(A_n)$. It follows that \hbox{$P(G/Z(G))=\{1\}$.}
\end{proof}

\medskip

In connection with Theorem \ref{pellequation}, we are now going to prove that also the pseudocentre of~$\operatorname{GL}(2,\mathbb Z)$ is trivial. In order to do this, we need some insight on the {\it negative} Pell's equation $x^2-y^2\sqrt d=-1$, where $d$ is a positive non-perfect square. First, note that this equation is solvable only for certain values of $d$, so we now assume it is solvable and we explain how a general solution is constructed (see \cite{Titu}, Theorem 3.4.1). Let $(u_0,v_0)$ be the minimal positive solution of the equation $x^2-y^2\sqrt d=-1$, and let $(x_n,y_n)$ be the positive solutions of \hbox{$x^2-y^2\sqrt d=1$} obtained from a fundamental solution $(x_1,y_1)$ as explained just before proving Theorem \ref{pellequation}. Then the general solution of $x^2-y^2\sqrt d=1$ is given by $(\pm u_n,\pm v_n)$, where $$u_n=u_0x_n+dv_0y_n\quad\textnormal{ and }\quad v_n=v_0x_n+u_0y_n.$$ The connection between $(u_0,v_0)$ and $(x_1,y_1)$ is explained by Remark 1 at p.137 of \cite{Titu} shows that $$x_1 + y_1\sqrt{d} = \big(u_0 - v_0\sqrt{d}\big)^2.$$


\medskip

\begin{theorem}\label{gl2z}
The pseudocentre of $\operatorname{GL}(2,\mathbb Z)$ coincides with the centre, while the pseudocentre of $\operatorname{PGL}(2,\mathbb Z)$ is trivial.
\end{theorem}
\begin{proof}
Set $G=\operatorname{GL}(2,\mathbb Z)$. Let $p$ be an odd prime number, and set $$A = \begin{pmatrix}
        1 & 2p \\
        p & 2p^2 +1
    \end{pmatrix}.$$
Similarly to what has been done in the proof of Theorem \ref{pellequation}, the centralizer $C$ of $A$ in $G$ turns out to be $$C = \left\{\begin{pmatrix}
        a & 2b \\
        b & a + 2pb 
    \end{pmatrix}: a^2 - 2b^2 + 2pab = \pm1, a,b \in \mathbb{Z}\right\}$$
    Put $x = a + pb$ and $y = b$. Then $a^2 - 2b^2 + 2pab = \pm 1$ becomes \hbox{$x^2 - (p^2 +2)y^2 = \pm 1$.}
    We claim that the equation $x^2 - (p^2 + 2)y^2 = -1$ has no integer solutions and from this it follows that $C_G(A) = C_{\operatorname{SL}(2,\mathbb{Z})}(A)$. 

    Assume by contradiction that $x^2 - (p^2 +2)y^2 = -1$ has an integer solution, let $(u_0,v_0)$ be the minimal positive solution, and let $(p^2+1,p)$ be a fundamental solution to Pell's equation \hbox{$x^2-(p^2+2)y^2=1$} (see again the proof of Theorem \ref{pellequation}). The above remark shows that $$(p^2+1) + p\sqrt{p^2+2} = \big(u_0 - v_0\sqrt{p^2+2}\big)^2.$$ Therefore, $u_0^2 + (p^2 +2)v_0^2 = p^2 +1$ and $p = 2u_0v_0$, which is a contradiction because $p$ is odd. The claim is proved.

    Now, note that $\Gamma(p)$ is a normal subgroup of $G$, so the same argument we employed in the final part of the proof Theorem \ref{pellequation} yields that $P(G)=Z(G)$.

\medskip

Finally, we show that the pseudocentre of $G/Z(G)$ is trivial. In order to do this, we note that the centralizer $D/Z(G)$ of $AZ(G)/Z(G)$ in $G/Z(G)$ is abelian and it is an extension of index at most $2$ of $AZ(G)/Z(G)$. If $D/Z(G)$ is infinite cyclic, then $D$ is abelian and so $D\leq C$. If $D/Z(G)$ is not infinite cyclic, then there is a subgroup $X$ of $G$ such that $D=AXZ(G)$ and $XZ(G)/Z(G)$ has order $2$. In particular, $XZ(G)$ has order at most $4$, and there are two possibilities for $D$: either $X$ has order $2$, so $D=X\ltimes \big(A\times Z(G)\big)$, or $X$ has order $4$, so $D=X\ltimes A$ --- in both cases, $A^X=-A$. In both cases, we obtain a contradiction as follows. Suppose first that $X^2=I_2$ and write $$X = \begin{pmatrix}
        a & b \\
        c & d
    \end{pmatrix}.$$ First, assume $X^2=I_2$, so in particular $$\begin{cases}
        a^2 + bc = 1 \\
        b(a+d) = 0 \\
        c(a+d) = 0 \\
        d^2 + cb = 1
    \end{cases}$$ If $\operatorname{det}(X)=1$, then $ad -1 = bc$ and the system becomes $$\begin{cases}
        a(a+d) = 2 \\
        b(a+d) = 0 \\
        c(a+d) = 0 \\
        d(a+d) = 2
    \end{cases}.$$ 
    Therefore $a + d \neq 0$ and so $b = c = 0$ and $a+ d = \pm 1,\pm 2$. Since $ad = 1$, it follows that $a+d = \pm 2$ and hence $a=d=\pm1$, contradicting $XA = -AX$. 

    If $\operatorname{det}(X)=-1$, so $bc = ad +1$ and the system becomes $$\begin{cases}
    a(a+d) = 0 \\
    b(a+d)= 0 \\
    c(a+d) = 0 \\
    d(a +d ) = 0
    \end{cases}
    $$ This implies $a +d = 0$ and so $$X = \begin{pmatrix}
        a & b \\
        c & -a
    \end{pmatrix}.$$ Now, equating the $(2,1)$-entries of the equation $XA = -AX$ gives $c -ap = -pa- (2p^2 +1)c$, which yields $c = 0$. On the other hand, equating the $(1,2)$-entries of the same equation gives $2pa + b(2p^2 +1) = -b + 2pa$, and so $b = 0$. Therefore $$X = \begin{pmatrix}
        a & 0 \\
        0 & -a
    \end{pmatrix},$$ so $a^2 = 1$ and hence $a=\pm1$, a contradiction.

\medskip

Now, assume $X^2 = -I_2$, so in particular $$\begin{cases}
        a^2 + bc = -1 \\
        b(a+d) = 0 \\
        c(a+d) = 0 \\
        d^2 +bc = -1
    \end{cases}.$$ If $\operatorname{det}(X)=1$, then $bc = ad -1$, and the system becomes $$\begin{cases}
        a(a+d) = 0 \\
        b(a+d) = 0 \\
        c(a+d) = 0 \\
        d(a+d) = 0
    \end{cases},$$ while if $\operatorname{det}(X)=-1$, it becomes $$\begin{cases}
        a(a+d) = -2 \\
        b(a+d) = 0\\
        c(a+d) = 0 \\
        d(a+d) = -2
    \end{cases}.$$ In both cases, we obtain a contradiction as above.

Thus, $D\leq C$ and we may prove that the pseudocentre of $G/Z(G)$ is trivial as we did in the final part of Theorem \ref{pellequation}.
\end{proof}

\medskip

Now, we are going to give another proof of Theorem \ref{pellequation} that is interesting in its own respect and uses the~Che\-byshev polynomials of the first and second kind. Recall that the~{\it Chebyshev polynomials of the first kind} are recursively defined as follows: $T_0(x)=1$, $T_1(x)=x$ and $T_{n+1}(x)=2xT_n(x)-T_{n-1}(x)$. It can be shown that $$T_n(x)=\frac{1}{2}\left(\left(x-\sqrt{x^2-1}\right)^n+\left(x+\sqrt{x^2-1}\right)^n\right)$$ for every $x\in\mathbb R$ and every $n\in\mathbb N$. On the other hand, {\it Chebyshev polynomials of the second kind} are recursively defined as follows: $U_i(x)=0$ for all $i<0$, $U_0(x)=1$, $U_1(x)=2x$ and $U_{n+1}(x)=2xU_n(x)-U_{n-1}(x)$. The following two auxiliary results are interesting in their own and show some intriguing link between Chebyshev polynomials and powers of matrices over arbitrary commutative rings.

\begin{lemma}\label{cheb} Let $R$ be a commutative ring with identity and let $g$ be an element of $\operatorname{SL}(2,R)$. Then $\operatorname{trace}(g^n)=2\cdot T_n\big(\operatorname{trace}(g)/2\big)$ for any non-negative integer $n$.
\end{lemma}
\begin{proof}
Set $t=\operatorname{trace}(g)$. By the Cayley--Hamilton theorem, $g$ is a solution of its own characteristic polynomial, so $g^2-tg+I_2=0$. Then $g^n=tg^{n-1}-g^{n-2}$ for every integer~\hbox{$n\geq2$.} By induction, $$
\begin{array}{c}
\operatorname{trace}(g^n)=\operatorname{trace}(tg^{n-1}-g^{n-2})=t\cdot\operatorname{trace}(g^{n-1})-\operatorname{trace}(g^{n-2})\\[0.2cm]
=2\cdot 2\cdot t/2\cdot T_{n-1}(t/2)-2\cdot T_{n-2}(t/2)
=2\cdot T_n(t/2)
\end{array}
$$ and the lemma is proved.
\end{proof}

\begin{lemma}\label{cheb2} Let $R$ be a commutative ring with identity and let $g$ be an element of $\operatorname{SL}(2,R)$. Then $g^n=U_{n-1}\big(\operatorname{trace}(g)/2\big)g-U_{n-2}\big(\operatorname{trace}(g)/2\big)I_2$ for any positive integer $n$.
\end{lemma}
\begin{proof}
Set $t=\operatorname{trace}(g)$. By the Cayley--Hamilton theorem, $g$ is a solution of its own characteristic polynomial, so $g^2-tg+I_2=0$. From this it follows that the statement is true for~\hbox{$n\leq 2$.} Assume $n>2$ and the usual induction hypothesis. Then
$$\begin{array}{c}
g^{n+1}=U_{n-1}(t/2)g^2-U_{n-2}(t/2)g
=U_{n-1}(t/2)(tg-I_2)-U_{n-2}(t/2)g\\[0.2cm]
=tU_{n-1}(t/2)g-U_{n-1}(t/2)I_2-U_{n-2}(t/2)g\\[0.2cm]
=\big(tU_{n-1}(t/2)-U_{n-2}(t/2)\big)g-U_{n-1}(t/2)I_2\\[0.2cm]
=U_{n}(t/2)g-U_{n-1}(t/2)I_2.
\end{array}$$
The lemma is proved.
\end{proof}

\begin{theorem}\label{pellequationcheb}
The pseudecentre of $\operatorname{SL}(2,\mathbb{Z})$ coincides with its centre and the psudocentre of $\operatorname{PSL}(2,\mathbb{Z})$ is trivial.
\end{theorem}
\begin{proof}
Let $p$ be an odd prime, and let $$g=\Big(\begin{array}{cc}
   1  &  p\\
    p &  p^2+1
\end{array}\Big)\in G=\operatorname{SL}(2,\mathbb Z).$$ We claim that $C_G(g)=\langle g\rangle Z(G)$. Assume there is an element
$$h=\Big(\begin{array}{cc}
   a  &  b\\
    c &  d
\end{array}\Big)
$$
of $G$ with $h^m=\pm g$ for some positive integer $m>1$ --- clearly, $m$ can be chosen to be a prime.  If $m=2$, then Lemma \ref{cheb} yields $\operatorname{trace}(h)^2-2=\pm(p^2+2)$, which is impossible. Thus, we may assume $m$ is odd and consequently that $h^m=g$.

Set $t=\operatorname{trace}(h)$, and note that $t>2$ because $\operatorname{trace}(g)>2$. Put $r=U_{m-1}\big(t/2\big)$ and \hbox{$s=U_{m-2}\big(t/2\big)$}. By Lemma \ref{cheb2}, we have that $g=h^m=rh-sI$, so $p=rb=rc$ and hence $b=c=1$ and $r=p$. However, we have $1=ra-s$ and $p^2+1=rd-s$, from which $p^2=r(d-a)$. Thus, $d-a=1$ because $\operatorname{det}(h)=ad-1=1$, and hence we have the contradiction $p^2=p$.

This contradiction shows that $C_G(g)=\langle g\rangle Z(G)$, and then the proof proceeds as in the final part of Theorem~\ref{pellequation}.
\end{proof}

\medskip

It now possible to state and prove the main result of this section.

\begin{theorem}\label{theoremafinalefreeproduct}
Let $G=H\ast K$ be the free product of the non-trivial groups $H$ and $K$ with $|H|>2$. Then $P(G)=\{1\}$.
\end{theorem}
\begin{proof}
This follows at once from Theorems \ref{freeproducttrivial}, \ref{pellequation} and/or \ref{pellequationcheb}.
\end{proof}

\begin{corollary}
Free groups have trivial pseudocentre.
\end{corollary}

\begin{corollary}
Let $N$ be a group, and let $H$ and $K$ be groups containing an isomorphic copy of $N$ as a proper normal subgroup. Let $G=H\ast_N K$ be the free product of $H$ and $K$ with amalgamated subgroup~$N$. If $|H/N|>2$, then $P(G)\leq N$.
\end{corollary}
\begin{proof}
Clearly, $G/N \simeq (H/N)\ast (K/N)$ is the free product of $H/N$ and $K/N$. Thus, by~The\-o\-rem~\ref{theoremafinalefreeproduct}, $P(G/N) = \{1\}$ and so $P(G) \le N$.
\end{proof}

\medskip

The previous corollary is the only real result we have managed to obtain for the pseudocentre of an amalgamated free product, and in the final part of this section, we explain why for such a construction the pseudocentre can be very hard to compute, even in some circumstances in which the group is known to behave well.

First, observe that in general a ‘‘non-trivial’’ amalgamated free product of two groups may be a simple group (see for example \cite{Burg} and \cite{Camm}). In particular,~The\-o\-rem~5.5 of \cite{Burg} proves the existence of a simple group which is the free product of finitely generated free groups with an amalgamated subgroup of finite index --- by The\-o\-rem~1.3 of~\cite{Wise}, the amalgamated subgroup cannot be malnormal (see also~\cite{DO}, Corollary 2.2 and Theorem 3.7, in which it is shown that there may even exist uncountably many normal subgroups). Recall that a subgroup $M$ of a group~$G$ is {\it malnormal} whenever $M\cap M^x=\{1\}$ for every~\hbox{$x\in G\setminus M$.} Clearly, the trivial subgroup and the whole group are malnormal subgroups, and in fact they are the only normal subgroups of $G$ which are malnormal at the same time. When $G$ is locally finite, a proper non-trivial malnormal subgroup $H$ is also called a {\it Frobenius complement}. In this case, the set $N=G\setminus \bigcup_{x\in G}\big(M^x\setminus\{1\}\big)$ is a normal subgroup of~$G$, called the {\it Frobenius kernel}, and $G=M\ltimes N$ (see \cite{Kegel}, The\-o\-rem~1.J.2) --- note that this is not true for arbitrary infinite groups (see the example on page~51 of~\cite{Kegel}). Groups with a malnormal subgroup that ‘‘splits’’ in the previous sense have been studied in \cite{dela}. It has been proved for example that if a semidirect product $G=M\ltimes N$ satisfies $C_G(x) \le N$ for all $x\in N$, then $M$ is malnormal in~$G$, while the converse does not hold.

Frobenius complements (and more generally ‘‘splitting’’ malnormal subgroups) can be very tempting subgroups to amalgamate (see for example \cite{karrass}). In fact, if $G_1=M\ltimes N_1$ and $G_2=M\ltimes N_2$ are Frobenius groups, then we easily see that $G_1\ast_MG_2\simeq M\ltimes N$, where $N=N_1\ast N_2$. Thus, one could guess that it would be possible to reproduce the proof Theorem \ref{freeproducttrivial} considering words inside $N$. However, this is not possible, as shown by the following example.

\begin{example}
There exists a group $G$ satisfying the following properties:
\begin{itemize}
    \item $G=M \ltimes N$, where $M$ is malnormal in~$G$ and $N$ is a certain  free product of two groups.
    \item There is a non-trivial element $z\in N$ such that $C_G(z)\not\leq N$.
\end{itemize}

\noindent Moreover, $N$ can be chosen to be infinitely generated  with periodic elements.
\end{example}
\begin{proof}
Let $G_i = M \ltimes N_i$, $i = 1,2$, be arbitrary Frobenius groups such that $M$ has an element $g$ of prime order $p\geq5$, and set $G = G_1\ast_M G_2$. Then $G=M\ltimes N$, where $N=N_1\ast N_2$. 

We claim that $M$ is a malnormal subgroup of $G$.  Let $y \in N\setminus M$, and suppose there is $1 \neq u\in M^y \cap M$. Then there is $x\in M$ such that $x^y=u$ and consequent\-ly~$[x,y]=x^{-1}u\in M\cap N=\{1\}$. Thus, $y^x=y$, which is only possible if $x=1$ (by the normal form of the amalgamated free product). This contradiction shows that~$M$ is malnormal in $G$.

Write $p=2m-1$ for some positive integer $m$, and choose non-trivial elements $h\in N_1$ and $k\in N_2$. Set $x =[h,k]$, $w=kh$\,\,\footnote{Here, $w$ can be replaced by any word (possibly the empty word) starting with an element of $N_2$ and ending with any element of $N_1$.} and $$
\begin{array}{c}
r_1 = \big(x^{2g^2}w^{-g^2} \ldots x^{2g^m}w^{-g^m}\big)\cdot \big(x^{2g^{-(m-2)}}w^{-g^{m-2}} \ldots x^{2g^{-1}}w^{-g^{-1}}\big)\cdot w^{g^{-1}}.
\end{array}
$$ Set $p_1\!=\! w^{-g}r_1w^{-g^{-1}}$, $u\!=\!w^{-1}x^{2g}p_1$, and $z\!=\!xux$. We claim that \hbox{$t = g\cdot (x^{-g} w x^{-1})$} centralizes $z$, that is $$z^g \cdot x^{-g} \cdot w \cdot x^{-1} = x^{-g} \cdot w \cdot x^{-1}\cdot z.$$

Since $o(g)=p=2m-1$, so $$
\begin{array}{c}
r_1^g = \big(x^{2g^3}w^{-g^3} \ldots x^{2g^{m+1}}w^{-g^{m+1}}\big)\cdot \big(x^{2g^{-(m-2)+1}}w^{-g^{(m-2)+1}} \ldots x^{2}w^{-1}\big)\cdot w\\[0.2cm]
=\big(x^{2g^3}w^{-g^3} \ldots x^{2g^{m}}w^{-g^{m}} x^{2g^{-(m-2)}}w^{-g^{-(m-2)}}\big)\cdot \big(x^{2g^{-(m-3)}}w^{-g^{(m-3)}} \!\!\!\ldots x^{2}w^{-1}\big)\cdot w\\[0.2cm]
=\big(x^{2g^3}w^{-g^3}\ldots x^{2g^{m}}w^{-g^{m}}\big)\cdot \big(x^{2g^{-(m-2)}}w^{-g^{-(m-2)}} x^{2g^{-(m-3)}}w^{-g^{(m-3)}} \!\!\!\ldots x^{2}w^{-1}\big)\cdot w\\[0.2cm]
\end{array}
$$ and hence $r_1w^{-g^{-1}} x^2=x^{2g^2} w^{-g^2}r_1^g$. Then $$
\begin{array}{c}
u^g = w^{-g}x^{2g^2}p_1^g = w^{-g}\cdot\big(x^{2g^2}w^{-g^2}r_1^g\big)\cdot w^{-1}\\[0.2cm]
=w^{-g}\cdot\big(r_1w^{-g^{-1}} x^2\big)\cdot w^{-1}=p_1x^2w^{-1}.
\end{array}
$$ and consequently, $$
\begin{array}{c}
z^g \cdot (x^{-g} \cdot w \cdot x^{-1}) = (x^g u^g x^g) \cdot (x^{-g} \cdot w \cdot x^{-1})\\[0.2cm]
= x^g\cdot u^g\cdot wx^{-1} =x^g\cdot p_1x^2w^{-1}\cdot wx^{-1}
=x^g p_1 x.
\end{array}
$$ On the other hand, $$
\begin{array}{c}
(x^{-g} w x^{-1}) \cdot z = (x^{-g} w x^{-1}) \cdot (x u x)\\[0.2cm]
= x^{-g} w\cdot u\cdot x =x^{-g} w\cdot (w^{-1} x^{2g} p_1)\cdot x = x^g p_1 x.
\end{array}
$$ The claim is proved.

In order to prove the moreover part of the statement, we observe that $G_1$ and~$G_2$ can be chosen to be infinite and locally finite (see \cite{Kegel}, Theorem 1.J.3).    
\end{proof}

\medskip

The question concerning free amalgamating products with respect to malnormal subgroups is completely solved in Section \ref{appendice} by Anthony Genevois.


\section{Free soluble groups}

The aim of this very short section is to exhibit our ignorance --- probably widely shared --- by proving that every non-abelian free soluble group has a trivial pseudocentre. This result was given without a proof in the paper of Wiegold \cite{Weig} together with other very easy statements. As the reader can see, our proof is not very difficult but it is still not easy at all.

\begin{theorem}\label{wiegoldcrepa}
Every non-abelian free soluble group $G$ has a trivial pseudocentre.
\end{theorem}
\begin{proof}
Since every non-abelian free soluble group is residually of rank $2$ (see for example \cite{hanna}, Problem 9, or \cite{guptalevin}), we may assume $G$ is freely generated by $a$ and $b$. Moreover, by \cite{hanna}, Corollary 26.32, $G$ is residually a finite $2$-group.

First, suppose $G'$ is abelian. Let $N$ be a normal subgroup of $G$ such that $G/N$ is a finite $2$-group of order $2^n$ (which, if needed, can be assumed large enough), and set $g=a^{2^n}[a,b^{2^n}]$. Clearly, $g$ is a non-trivial element of $G$ that is contained in $N$. If $x$ is any element of $G$ with $x^p=g$ for some prime $p$, then looking at the images in the abelianization shows that $p=2$ and that $x=a^{2^{n-1}}u$ for some~\hbox{$u\in G'$.} Now, consider the quotient of $G$ isomorphic to the standard wreath product $\langle b\rangle\wr\langle a\rangle\simeq\mathbb Z_{2^{n+1}}\wr \mathbb Z_{2^{n+1}}$. It is easy to see that the image of $g$ in such a quotient does not have a square root of type $a^{2^{n-1}}u$. This contradiction shows that there is no element $x\in G$ with $x^p=g$, where $p$ is a prime. By Theorem 1 of~\cite{malcev} (see also \cite{auslander}), we have that $C_G(g)=\langle g\rangle$, so $C_G(g)^G\leq N$. The arbitrariness of $N$ proves that~\hbox{$P(G)=\{1\}$.}


Now, suppose $G$ is not metabelian. Note that $G'$ is freely generated by the commutators $[a^{k_1},b^{k_2}]$, where $k_1,k_2\in\mathbb Z$ and $k_1,k_2\neq0$ (see for example \cite{putman}). Let~$N$ be a normal subgroup of $G$ such that $G/N$ is a finite $2$-group of order $2^n$, and set $g=[a,b^{2^n}]$. Then $g\in N\setminus\{1\}$ and we easily see that $C_G(g)=\langle g\rangle$ (see again \cite{malcev}, Theorem 1, or \cite{auslander}). Thus, $C_G(g)^G\leq N$ and again $P(G)=\{1\}$.
\end{proof}

\section{Burnside groups}

Let $m$ be a cardinal number, and let $n$ be a positive integer. If $F_m$ is the free group of rank $m$, then the {\it free $m$-generator Burnside group} $B(m,n)$ of exponent~$n$ is defined to be the quotient group of $F_m$ by the subgroup $F_m^n$ of $F_m$ generated by all $n$th powers of elements of $F_m$. It is well-known that if $n$ is a large enough ($>10^{78}$) odd positive integer, then $B(m,n)$ is infinite (see for example \cite{oldshanski}). The aim of this short section is to show that similarly to free groups and free products, also many free Burnside groups have a trivial pseudocentre.

\begin{lemma}\label{Burn}
Let $m\geq2$ be a cardinal number, $n$ an odd positive integer with $n > 10^{78}$, and put $G=B(m,n)$. Then $G'$ contains elements of order $n$.
\end{lemma}
\begin{proof}
Assume that $G$ is freely generated by $\mathcal A$ in the variety of groups of exponent $n$. Let $a_1,a_2$ be distinct elements of $\mathcal A$, let $B=B(2,n)$ be the Burnside group freely generated by~\hbox{$\{g_1,g_2\}$.} The universal property of $G$ yields that the assignations $$a_1\mapsto g_1,\;\, a_2\mapsto g_2,\;\, a\mapsto 1\quad\forall a\in\mathcal A\setminus\{a_1,a_2\}$$ can be uniquely extended to the surjective homomorphism $\overline{\varphi}: G \to B$. Set \hbox{$\varphi=\overline\varphi|_{\langle a_1,a_2\rangle}$.}  Similarly, the assignations $g_i\mapsto a_i$, $i=1,2$, can be uniquely extended to the surjective homomorphism $\psi: B \to \langle a_1,a_2\rangle$. 

By universal property of $B$, it follows that $\psi \circ \varphi = \operatorname{Id}_B$, so $\psi$ is injective. Therefore, $\psi: B \to \langle a_1,a_2 \rangle$ is an isomorphism. 

Since $B$ is not soluble (otherwise it would be finite), its derived subgroup contains a copy of $G$ (see the main theorem of \cite{Iva}) and so contains elements of order~$n$. It follows that $\langle a_1,a_2\rangle'\leq G'$ contains elements of order~$n$ as well.
\end{proof}

\begin{lemma}\label{arbitrarycardinal}
Let $m\geq2$ be a cardinal number, $n$ an odd positive integer with \hbox{$n > 10^{78}$,} and $G = B(m, n)$. Then, every non-cyclic subgroup $H$ of $G$ contains a subgroup $K\simeq B(m,n)$ such that if $N \trianglelefteq K$, then $N^G \cap K = N$.
\end{lemma}
\begin{proof}
Let $B \simeq {B}(2,n)$. By Theorem 39.1 \cite{oldshanski}, it follows that $G$ can be embedded in $B$. Now, if $H$ is a non-cyclic subgroup of $G$, then there exists a subgroup $K \simeq B(\omega,n)$ of $P$ such that whenever $N \trianglelefteq K$, then $N^B \cap K = N$ (see the main theorem of \cite{Iva}). Since $N^G \le N^B$, so $N^G \cap K = N$ and the result follows.
\end{proof}

\begin{theorem}\label{burnsidegroups}
Let $m\geq2$ be a cardinal number, $n$ an odd positive integer with \hbox{$n > 10^{78}$,} and put $G =B(m,n)$. Then $P(G)=\{1\}$.
\end{theorem}
\begin{proof}
Let $P$ be the pseudocentre of $G$. Assume that $P\neq \{1\}$. If $P$ is cyclic, then $C_G(P)$ is a finite-index cyclic subgroup of $G$ (see page 26 of \cite{Ata}) and so $G$ is finite, a contradiction. Thus, $P$ is non-cyclic, and it follows from~Lem\-ma~\ref{arbitrarycardinal} that~$P$ contains a subgroup $H\simeq B(\omega,n)$, properly, satisfying the following condition:
\begin{itemize}
    \item if $K \trianglelefteq H$, then $K^G \cap H = K$.
\end{itemize}
By Lemma \ref{Burn}, $H'$ contains an element $a$ of order $n$.  Then $C\!=\!C_G(a)$ is cyclic, and so~\hbox{$C_G(a)\!=\!\langle a \rangle$} because $\operatorname{exp}(G)=n$. Set $K = C^H $, so in particular~\hbox{$K \leq H'<H$.} Since $H\leq P\leq C^G=K^G$, so \hbox{$H\cap K^G=H$.} On the other hand, the above property of $H$ yields that $H\cap K^G=K$. Thus,~\hbox{$H=K$}
which is impossible. This contradiction shows that \hbox{$P=\{1\}$.}
\end{proof}


\section{Appendix: acylindrically hyperbolic groups}\label{appendice}

In this appendix, we use some tools coming from geometric group theory in order to show that the pseudocenter of any member of a broad family of groups, namely acylindrically hyperbolic groups, is always finite. Recall from \cite{MR3430352} that a group $G$ is \emph{acylindrically hyperbolic} if it admits an action on a Gromov-hyperbolic space $X$ that is non-elementary and \emph{acylindrical}, i.e.\ for every $d \geq 0$, there exist $L,N \geq 0$ such that
$$\forall x,y \in X, \ d(x,y) \geq L \Rightarrow \# \{ g \in G \mid d(x,gx),d(y,gy) \leq d \}.$$
Roughly speaking, the action of $G$ on $X$ is allowed to have infinite (quasi-)stabilisers of points, but the (quasi-)stabiliser of a pair of points that are sufficiently far away apart must be (uniformly) finite. Examples of acylindrically hyperbolic groups include hyperbolic groups (e.g.\ non-abelian free groups, fundamental groups of surfaces of genus $\geq 2$, various small cancellation groups, various random groups), relatively hyperbolic groups (e.g.\ free products that are not virtually cyclic, fundamental groups of finite-volume hyperbolic manifolds), mapping class groups of non-sporadic surfaces, outer automorphism groups of free groups of finite rank $\geq 2$, many $3$-manifold groups, groups of deficiency $\geq 2$, Cremona groups, most graph braid groups. We refer the reader to \cite{MR3430352, MR3966794, MR4057355} for more information. 

\medskip 
The main result of the appendix is the following statement.

\begin{theorem}\label{thm:AcylHyp}
Let $G$ be an acylindrically hyperbolic group. The pseudocentre of $G$ is contained in the finite radical $E(G)$. 
\end{theorem}

Here, the finite radical of $G$ refers to the unique maximal normal finite subgroup $E(G)$ of $G$ (see \cite[Theorem~2.24]{DGO}). As a consequence of Theorem~\ref{thm:AcylHyp}, pseudocentres of torsion-free acylindrically hyperbolic groups are always trivial. In full generality, the structure of the pseudocentre of $G$ will depend on the pseudocentre of $E(G)$ but also on how $G$ acts on $E(G)$ by conjugation. 

\medskip

In order to prove Theorem~\ref{thm:AcylHyp}, the following preliminary observation will be needed. 

\begin{lemma}\label{lem:FreeAvoid}
Let $F$ be a free group of finite rank $\geq 2$. For every finite subset $B \subset F \backslash \{1\}$, there exists a primitive element $h \in F$ such that $\langle\langle h \rangle\rangle \cap B = \emptyset$. 
\end{lemma}

\noindent
In our lemma, and in the rest of the appendix, we denote by $\langle\langle \cdot \rangle \rangle_G$ the normal closure in a group $G$, or just $\langle\langle \cdot \rangle\rangle$ if the group under consideration is clear. 

\medskip \noindent
It is worth mentioning that, because the centraliser of a non-trivial primitive element $h$ coincides with $\langle h \rangle$, Lemma~\ref{lem:FreeAvoid} provides an alternative proof of the triviality of pseudocentres of non-abelian free groups. Our proof is based on classical small cancellation theory. We refer the reader to \cite[Chapter~V]{MR1812024} for more information on the subject.

\medskip

\noindent{\sc Proof of Lemma~\ref{lem:FreeAvoid}} ---
Fix a basis $\{a_1, \ldots, a_r\}$ of $F$ and an integer $n \geq 10$ greater than the length of every element in $B$ (taken with respect to our basis). Consider the primitive element
$$h:= a_1a_2 a_1^2 a_2^2 \cdots a_1^n a_2^n.$$
Because we took $n \geq 10$, the presentation
$$\langle a_1, \ldots, a_r \mid a_1a_2 a_1^2 a_2^2 \cdots a_1^n a_2^n = 1 \rangle$$
of $F/ \langle\langle h \rangle\rangle$ satisfies the small cancellation condition C'(1/6). Indeed, the longest piece we can find is $a_1^{n-1}a_2^{n-1}$, whose length $2(n-1)$ is smaller than the sixth of the total length $n(n+1)$ of our relation as soon as $n \geq 10$. As a consequence, if $g \in F$ is a non-trivial element that belongs to $\langle\langle h \rangle\rangle$, then it follows from Dehn's algorithm that a reduced word representing $g$ must contain a subword of our relation of length $>n(n+1)/2 \geq n$. Thus, $g$ must have length $>n$, and consequently it cannot belong to $B$. We conclude that $\langle\langle h \rangle\rangle \cap B=\emptyset$, as desired. \hfill\qedbox

\medskip

\noindent{\sc Proof of Theorem~\ref{thm:AcylHyp} --- }Since $G \twoheadrightarrow G/E(G)$ sends the pseudocentre of $G$ inside the pseudocentre of $G/E(G)$, and since $G/E(G)$ is an acylindrically hyperbolic group with trivial finite radical \cite{Correction}, it suffices to show that every acylindrically hyperbolic group with trivial finite radical has a trivial pseudocentre in order to deduce Theorem~\ref{thm:AcylHyp}. So, from now on, we assume that $G$ has trivial finite radical. 

\medskip \noindent
Given a non-trivial element $g \in G$, our goal is to find an element $h \in G$ such that $g \notin \langle\langle C_G(h) \rangle\rangle$. Our arguments are based on the notion of \emph{hyperbolically embedded subgroups}, as introduced in \cite{DGO}. Given a group $G$, a subgroup $H$, and a generating set $X \subset G$, we write $H \hookrightarrow_h (G,X)$ to mean that $H$ is hyperbolically embedded in $G$ with respect to $X$. 

\medskip \noindent
We know from \cite[Theorem~6.14]{DGO} that there exist a generating set $X \subset G$ and a free subgroup $F \leq G$ of rank two such that $F \hookrightarrow_h (G,X)$. According to \cite[Corollary~4.27]{DGO}, we also have $F \hookrightarrow_h (G, X \cup \{g^{\pm 1}\})$. As a consequence of \cite[Theorem~2.27]{DGO}, there exists a finite subset $\mathcal{F} \subset F\backslash \{1\}$ such that, for every $N \lhd F$ avoiding $\mathcal{F}$, every element of $\langle\langle N \rangle\rangle_G$ either has a conjugate in $N$ or has unbounded orbits in $\mathrm{Cayl}(G, X \cup \{g^{\pm 1}\} \cup F)$. It follows from Lemma~\ref{lem:FreeAvoid} that there exists a primitive element $h \in F$ such that $\langle\langle h \rangle\rangle_F$ avoids $\mathcal{F}^+$, where $\mathcal{F}^+:= \mathcal{F}$ if no conjugate of $g$ belongs to $F$ and otherwise $\mathcal{F}^+:= \mathcal{F} \cup \{aga^{-1}\}$ for an arbitrary conjugate $aga^{-1}$ of $g$ belonging to $F$ which we fix once for all. Since $g$ has bounded orbits in $\mathrm{Cayl}(G, X \cup \{g^{\pm 1}\} \cup F)$, we know that, if $g \in \langle\langle h \rangle\rangle_G$, then $g$ must have a conjugate in $\langle\langle h \rangle\rangle_F$. We claim that this is not possible, which will prove that $g \notin \langle\langle h \rangle\rangle_G$.

\medskip \noindent
If $g$ has a conjugate in $\langle\langle h \rangle\rangle_F \leq F$, say $bgb^{-1}$, then
$$aga^{-1} = (ab^{-1}) bgb^{-1} (ab^{-1})^{-1} \in F \cap  (ab^{-1}) F (ab^{-1})^{-1}.$$
But $F$ is almost malnormal according to \cite[Proposition~2.10]{DGO}, so this implies that $ab^{-1} \in F$. Since $aga^{-1}$, which belongs to $\mathcal{F}^+$, does not belong to $\langle\langle h \rangle\rangle_F$ by construction, it follows that $bgb^{-1}$ cannot belong to $\langle\langle h \rangle\rangle_F$ either, a contradiction.

\medskip \noindent
So far, we have proved that there exists a primitive element $h \in F$ such that $g \notin \langle\langle h \rangle\rangle_G$. In order to conclude our proof, it suffices to verify that $C_G(h)= \langle h \rangle$. 

\medskip \noindent
Notice that, for every $k \in C_G(h)$, 
$$\langle h \rangle \leq kFk^{-1} \cap F,$$
which implies, again because $F$ is almost malnormal, that $k \in F$. We deduce from this observation that $C_G(h)=C_F(h)$. But we chose $h$ primitive in $F$, hence $C_F(h)= \langle h \rangle$. We conclude that $C_G(h)= \langle h \rangle$, as desired. \hfill\qedbox

\medskip

As concrete applications of Theorem~\ref{thm:AcylHyp}, let us mention a couple of examples for which we can describe the pseudocentres. For our first application, given a group $G$ and a subgroup $H \leq G$, we say that $H$ is \emph{slightly malnormal} if there exists some $g \in G$ such that $H \cap gHg^{-1}= \{1\}$. 

\begin{corollary}
Let $A,B$ two groups and $C \leq A,B$ a common proper subgroup. If $C$ is slightly malnormal in $G:=A \ast_C B$, then the pseudocentre of $G$ is trivial; unless $C=\{1\}$ and $A=B= \mathbb{Z}_2$, in which case the pseudocentre is an index-$4$ infinite cyclic subgroup.
\end{corollary}

\begin{proof}
We know from \cite[Corollary~2.2]{DO} that $G$ is either virtually cyclic or acylindrically hyperbolic. In the latter case, it follows from Theorem~\ref{thm:AcylHyp} that the pseudocentre of $G$ is contained in its finite radical $E(G)$, so it suffices to verify that $G$ has trivial finite radical. According to Claim~\ref{claim:AmalgamFiniteNormal} below, $E(G) \leq C$. But we know by assumption that $C \cap gCg^{-1}= \{1\}$ for some $g \in G$, hence
$$E(G) = gE(G)g^{-1} \leq C \cap gCg^{-1} = \{1\}.$$
We conclude, as desired, that the pseudocentre of $G$ is trivial.

\begin{claim}\label{claim:AmalgamFiniteNormal}
Every finite normal subgroup in a non-trivial amalgamated product $P \ast_R Q$ is contained in $R$.
\end{claim}

\noindent
Let $N$ be a finite normal subgroup in $P \ast_R Q$. According to \cite[Corollary p.\ 36]{MR1954121}, $N$ must be contained in a conjugate of $P$ or $Q$, say $P$. Since $N$ is normal, actually $N$ must be contained in $P$ itself. Given some $q \in Q \backslash R$, $N=qNq^{-1}$ must be contained in $P \cap qPq^{-1}$. But $P \cap qPq^{-1} \leq R$. This follows from the normal form in amalgamated product, or from Bass-Serre theory since $(P,Q,qP)$ defines a geodesic in the Bass-Serre tree and that then an isometry fixing the vertices $P$ and $qP$ must stabilise $Q$ as well. This concludes the proof of Claim~\ref{claim:AmalgamFiniteNormal}. 

\medskip \noindent
It remains to consider the case where $G$ is virtually cyclic. According to Claim~\ref{claim:AmalgamVirtCyclic} below, this happens precisely when $A,B,C$ are all finite and $C$ has index $2$ in both $A$ and $B$. But then, $C$ is normal in both $A$ and $B$, which implies that $C$ is normal in $G$. Since $C$ is slightly malnormal, it follows that $C=\{1\}$. In other words, $G$ coincides with the infinite dihedral group $\operatorname{Dih}(\infty) \simeq \mathbb{Z}_2 \ast \mathbb{Z}_2$. In this case, the pseudocentre can be computed by a direct computation.

\begin{fact}\label{fact:Dinfty}
The pseudocentre of $\operatorname{Dih}(\infty)$ is an infinite cyclic subgroup of index $4$. 
\end{fact}

\noindent
If we think of $\operatorname{Dih}(\infty)$ as presented by $\langle a,b \mid a^2=b^2=1 \rangle$, then the centraliser of $a$ is $\langle a \rangle$, whose normal closure is $\{ (ab)^na(ab)^{-n}, (ab)^{2n} \mid n \in \mathbb{Z} \}$. Similarly, the normal closure of the centraliser of $b$ is $\{ (ab)^nb(ab)^{-n},(ba)^{2n} \mid n \in \mathbb{Z}\}$. Finally, for every $k \in \mathbb{Z} \backslash \{0\}$, the centraliser of $(ab)^k$ is the normal subgroup $\langle ab \rangle$. The intersection between all these normal closures, i.e.\ the pseudocentre of $\operatorname{Dih}(\infty)$, is $\langle (ab)^2 \rangle$, which is an infinite cyclic subgroup of index $4$. This concludes the proof of our corollary.

\begin{claim}\label{claim:AmalgamVirtCyclic}
A non-trivial amalgamated product $P \ast_R Q$ is virtually cyclic if and only if $P,Q,R$ are all finite and $R$ has index $2$ in both $P$ and $Q$.
\end{claim}

\noindent
If $P,Q,R$ are all finite and $R$ has index $2$ in both $P$ and $Q$, then $R$ is a normal subgroup in $P \ast_R Q$ and $(P \ast_R Q)/ R$ is isomorphic to the free product $(P/R) \ast (Q/R) \simeq \mathbb{Z}_2 \ast \mathbb{Z}_2$. Consequently, $P \ast_R Q$ is virtually cyclic. Conversely, assume that $P \ast_R Q$ is virtually cyclic. Then, it has to stabilise a bi-infinite line in its Bass-Serre tree (see \cite{MR1954121,MR564422} for more information on Bass-Serre theory), which amounts to saying that the tree itself must be a bi-infinite line. It follows that $R$ has index $2$ in both $P$ and $Q$. As a consequence, $R$ is normal in $P \ast_R Q$. Since $P \ast_R Q$ surjects onto the free product $(P/R) \ast (Q/R) \simeq \mathbb{Z}_2 \ast \mathbb{Z}_2$, which is virtually $\mathbb{Z}$, the fact that $P \ast_R Q$ is virtually cyclic imposes that the kernel of this projection, namely $R$, must be finite. Since $R$ has finite index in both $P$ and $Q$, this implies that all $P,Q,R$ are finite. 
\end{proof}

\noindent
For our second application, recall that, given a graph $\Gamma$ and a collection of groups $\mathcal{G}= \{ G_u \mid u \in V(\Gamma)\}$ indexed by the vertex-set $V(\Gamma)$ of $\Gamma$, the \emph{graph product} $\Gamma \mathcal{G}$ is defined by the relative presentation
$$\langle G_u \ (u \in V(\Gamma)) \mid [G_u,G_v]=1 \ (u,v \in V(\Gamma) \text{ adjacent}) \rangle$$
where $[H,K]=1$ means that $[h,k]=1$ for all $h \in H$ and $k \in K$. Usually, one says that graph products interpolate between free products (when $\Gamma$ has no edge, i.e.\ ``nothing commute'') and direct sums (when $\Gamma$ is complete, i.e.\ ``everything commute''). Graph products of infinite cyclic groups (resp.\ cyclic groups of order two) coincide with right-angled Artin groups (resp.\ right-angled Coxeter groups). 

\begin{corollary}\label{cor:GraphProd}
Let $\Gamma$ be a finite graph and $\mathcal{G}$ a collection of non-trivial groups indexed by $V(\Gamma)$. The pseudocentre of $\Gamma \mathcal{G}$ is isomorphic to
$$\mathbb{Z}^p \oplus \bigoplus\limits_{u \text{ central}} P(G_u)$$
where $p \geq 0$ denotes the number of $\ast$-factors of $\Gamma$ are pairs of non-adjacent vertices both labelled by $\mathbb{Z}_2$. 
\end{corollary}

\noindent
Recall that a graph $\Xi$ decomposes as a \emph{join} $\Xi_1 \ast \Xi_2$ if $\Xi_1, \Xi_2 \leq \Xi$ are two non-empty subgraphs such that every vertex of $\Xi_1$ is adjacent to every vertex of $\Xi_2$. A \emph{central vertex} is a vertex that is adjacent to all the other vertices. A graph with a central vertex is automatically a join.

\noindent{\sc Proof of Corollary~\ref{cor:GraphProd} --- }We start by proving the following particular case:

\begin{claim}\label{claim:IrreducibleCase}
If $\Gamma$ is not a join and contains at least two vertices, then the pseudocentre of $\Gamma \mathcal{G}$ is trivial; unless $\Gamma$ is a pair of non-adjacent vertices both labelled by $\mathbb{Z}_2$, in which case its pseudocentre is infinite cyclic. 
\end{claim}

\noindent
We know from \cite[Corollary~2.13]{DO} that $\Gamma \mathcal{G}$ is either virtually cyclic or acylindrically hyperbolic. According to Claim~\ref{claim:GPvirtCyclic} below, $\Gamma \mathcal{G}$ is virtually cyclic if and only if $\Gamma \mathcal{G} \simeq \operatorname{Dih}(\infty)$, which has infinite cyclic pseudocentre according to Fact~\ref{fact:Dinfty}. So we can assume that $\Gamma \mathcal{G}$ is acylindrically hyperbolic, in which case we know from Theorem~\ref{thm:AcylHyp} that its pseudocentre is contained in its finite radical $E(\Gamma \mathcal{G})$. Notice that, for every vertex $u \in V(\Gamma)$, $\Gamma \mathcal{G}$ splits as 
$$\langle \mathrm{star}(u) \rangle \ast_{\langle \mathrm{link}(u) \rangle} \langle \Gamma \backslash \{u\} \rangle.$$
Here, given a subgraph $\Xi \leq \Gamma$, we set $\langle \Xi \rangle := \langle G_v \mid v \in V(\Xi) \rangle$. Recall that the \emph{link} (resp.\ \emph{star}) of a vertex $v$ refers to the subgraph induced by all the neighbours of $v$ (resp.\ by $v$ and all its neighbours). It follows from Claim~\ref{claim:AmalgamVirtCyclic} that $E(\Gamma \mathcal{G})$ is contained in $\langle \mathrm{link}(u) \rangle$. Hence
$$E(\Gamma \mathcal{G}) \leq \bigcap\limits_{u \in V(\Gamma)} \langle \mathrm{link}(u) \rangle = \left\langle \bigcap\limits_{u \in V(\Gamma)} \mathrm{link}(u) \right\rangle = \langle \emptyset \rangle = \{1\},$$
where the first equality is justified by \cite[Lemma~3.3]{MR3365774}. This concludes the proof of Claim~\ref{claim:IrreducibleCase} modulo:

\begin{claim}\label{claim:GPvirtCyclic}
If $\Gamma$ is not a join and contains at least two vertices, then $\Gamma \mathcal{G}$ is virtually cyclic if and only if it is isomorphic to the infinite dihedral group $\operatorname{Dih}(\infty)$. 
\end{claim}

\noindent
Because $\Gamma$ contains at least two vertices and is not complete, it must contain at least two vertices that are not adjacent, say $u,v \in V(\Gamma)$. If $\Gamma \mathcal{G}$ is virtually cyclic, it follows that $G_u \simeq G_v \simeq \mathbb{Z}_2$ since $\Gamma \mathcal{G}$ contains the free product $\langle G_u,G_v \rangle \simeq G_u \ast G_v$. If $\Gamma$ does not contain a vertex distinct from both $u$ and $v$, we are done. Otherwise, let $w \in V(\Gamma)$ be a third vertex. If $w$ is not adjacent to $v$, then $\Gamma \mathcal{G}$ contains the free product $\langle G_u,G_v,G_w \rangle \simeq \langle G_u , G_w \rangle \ast G_v$, and a fortiori a non-abelian free subgroup, which is impossible. So $w$ must be adjacent to $v$. Similarly, we know that $w$ must be adjacent to $u$ as well. In other words, all the vertices of $\Gamma$ distinct from $u$ and $v$ must be adjacent to both $u$ and $v$. But $\Gamma$ is not a join, so necessarily there cannot be vertices distinct from $u$ and $v$. This concludes the proof of Claim~\ref{claim:GPvirtCyclic}. 

\medskip \noindent
Now, let us consider the general case, without any assumption on $\Gamma$. Decompose $\Gamma$ as a join $\Gamma_0 \ast \Gamma_1 \ast \cdots \ast \Gamma_n$ where $\Gamma_0$ is the subgraph of $\Gamma$ spanned by all the central vertices and where $\Gamma_1, \ldots, \Gamma_n$ are not joins and each contains at least two vertices. Clearly, 
$$\Gamma \mathcal{G}= \left( \bigoplus\limits_{u \text{ central}} G_u \right) \oplus \langle \Gamma_1 \rangle \oplus \cdots \oplus \langle \Gamma_n \rangle,$$
where we note $\langle \Gamma_i \rangle:= \langle G_u \mid u \in V(\Gamma_i) \rangle$ for every $1 \leq i \leq n$.  Each $\langle \Gamma_i \rangle$ is naturally a graph product (see \cite[Proposition~3.1]{Green}), and Claim~\ref{claim:IrreducibleCase} implies that its pseudocentre is infinite cyclic when our graph product is infinite dihedral and trivial otherwise. Hence
$$P(\Gamma \mathcal{G})= \bigoplus\limits_{u \text{ central}} P(G_u) \oplus P(\langle \Gamma_1 \rangle) \oplus \cdots \oplus P(\langle \Gamma_n \rangle) \simeq \mathbb{Z}^p \oplus \bigoplus\limits_{u \text{ central}} P(G_u)$$
where $p$ denotes the number of indices $i$ for which $\Gamma_i$ is a pair of non-adjacent vertices both labelled by $\mathbb{Z}_2$. The first equality above is justified by Lemma~2.6.\hfill\qedbox

\chapter{Other groups of matrices}\label{sectmatrice}

In Section \ref{secmclain} we have dealt with the pseudocentre of the group of unitriangular matrices over a field, while in Section \ref{freeproductsect} we have studied the pseudocentre of invertible groups of matrices over the integers. In this section we  study other relevant groups of matrices over a field. First, note that $\operatorname{GL}(2,2)\simeq\operatorname{Sym}(3)$, $\operatorname{PGL}(2,3)\simeq\operatorname{Sym}(4)$ and $\operatorname{PSL}(2,3)\simeq\operatorname{Alt}(4)$, so the following result takes care of all ‘‘small’’ cases.

\begin{theorem}\label{GL23}
Let $G = \operatorname{GL}(2,3)$ and $S=\operatorname{SL}(2,3)$. Then $P(G)=S$ and $P(S)=O_2(S)$.
\end{theorem}
\begin{proof}
If $g  = \begin{pmatrix}
        1 & 0 \\
        1 & 1
    \end{pmatrix}$, then $$C_G(g) =\left\{\begin{pmatrix}
        a & 0 \\
        b & a
    \end{pmatrix} : a \in \mathbb{F}_3^\times,\, b \in \mathbb{F}_3\right\}\simeq\mathbb Z_6.$$ Then $C_G(g)^G=S$ and so $P(G)\leq S$. If $x$ is a non-central element of $G$ such that $\langle x\rangle$ is not subnormal in $G$, then $\langle x\rangle^G$ contains $S$, so also $S\leq C_G(x)^G$. Suppose that $x$ is a non-central element of $G$ such that $\langle x\rangle$ is subnormal in $G$. Then $x\in O_2(G)\simeq Q_8$ can be one of the following matrices $$\pm\begin{pmatrix}
        0 & 1 \\
        -1 & 0
    \end{pmatrix},\quad \pm\begin{pmatrix}
        -1 & -1 \\
        -1 & 1
    \end{pmatrix}\quad\textnormal{ and }\quad \pm\begin{pmatrix}
        -1 & 1 \\
         1 & 1
    \end{pmatrix}$$ But these matrices are centralized respectively by the matrices $$\begin{pmatrix}
        1 & 1 \\
        -1 & 1
    \end{pmatrix},\quad\begin{pmatrix}
        0 & 1 \\
        1 & 1
    \end{pmatrix}\quad\textnormal{ and }\quad \begin{pmatrix}
        0 & 1 \\
        1 & -1
    \end{pmatrix},$$ which are not contained in $O_2(G)$. Consequently, $C_G(x)^G\geq S$ and so $P(G)=S$.

Now, $S/Z(S)\simeq\operatorname{Alt}(4)$ and $P(S/Z(S))=O_2(G)/Z(S)$, so $P(S)\leq O_2(G)$. On the other hand, it is easy to see that $P(S)\geq O_2(G)$ and so $P(S)=O_2(G)$. 
\end{proof}

\medskip

The following all-in-one result deals with all the other cases.

\begin{theorem}\label{GL}
Let $\mathbb{F}$ be a field and $n$ a positive integer with either $|\mathbb{F}| > 3$ and $n = 2$, or $n >2$. Let~\hbox{$G=\operatorname{GL}(n,\mathbb F)$} and $S=\operatorname{SL}(n,\mathbb F)$.
    \begin{itemize}
        \item[\textnormal{(1)}] $P(G)=Z(G)\cdot S$, that is, the set of all matrices whose determinant is an $n$-th power in~$\mathbb F$. 
        \item[\textnormal{(2)}] $P(S) = S$.
        \item[\textnormal{(3)}] $P(S/Z(S))=S/Z(S)$.
        \item[\textnormal{(4)}] $P(G/Z(G))=SZ(G)/Z(G)$.
        \item[\textnormal{(5)}] $P(P(G))=P(G)$.
    \end{itemize}
\end{theorem}
\begin{proof}
(1)\quad Let $g\in G\setminus Z(G)$ and set $C = C_G(g)$. Then~\hbox{$[g,S] \neq \{1\}$} by~The\-o\-rem~3.2.5 of~\cite{acourse}. If $C^G$ contains a non-central element of $S$, then $C^G \ge S$ (see~The\-o\-rem~3.2.8 in~\cite{acourse}). Assume this is not the case, so $$[g,S]\leq C^G\cap S\leq Z(S)\leq Z(G)$$ and $g\not\in S$. Now, let $L$ be a subgroup of $S$ such that $L\cap Z(S)=\{1\}$, so in particular,~\hbox{$C_L(g)=\{1\}$.} Since $[g,S]\leq Z(G)$, so the map $$\varphi: x\in L \mapsto [g,x]\in [g,L]\leq Z(G)\simeq \mathbb F^\times$$ is an injective homomorphism. Therefore $L$ is isomorphic to a subgroup of $\mathbb F^\times$. 

Now, the consideration of $L$ as the group of all unitriangular matrices $\operatorname{Tr}_1(n,\mathbb F)$ immediately yields $n=2$ and so $|\mathbb{F}|>3$. If $\mathbb F^\times$ has an element $b$ whose order is either infinite or a power of~$2$, then $K\ltimes\operatorname{Tr}_1(2,\mathbb F)$, where $$K=
\left\langle\begin{pmatrix}
b & 0\\
0 &b^{-1}
\end{pmatrix}\right\rangle,
$$ has trivial intersection with $Z(S)$ and is non-abelian. Therefore $\mathbb F^\times$ is a $2$-group. Consequently, $\operatorname{char}(F)\neq2$ and hence $\operatorname{Tr}_1(2,\mathbb F)$ is a $2'$-group. But since this group trivially intersects $Z(S)$, again we obtain a contradiction.

Therefore $S\leq C_G(g)^G$ for every $g\in G$, and hence $S\cdot Z(G)\leq P(G)$. Let $$a=\begin{pmatrix}
      1 & 0 & \cdots & 0 \\
      1 & 1 & \cdots & 0\\ 
     \vdots & \vdots & \ddots & \vdots \\
     1 & 1 & \ldots & 1
    \end{pmatrix}.$$ We claim that $C_G(a)^G = Z(G) \cdot S$. Let  $$x = \begin{pmatrix}
        x_{1,1} & \cdots & x_{1,n} \\
        \vdots & \ddots & \vdots \\
        x_{n,1} & \cdots & x_{n,n}
    \end{pmatrix}$$ be such that $xa=ax$. Clearly, $$ax = \begin{pmatrix}
        x_{1,1} & x_{1,2} & \ldots & x_{1,n} \\
        \sum_{i=1}^2 x_{i,1} & \sum_{i=1}^2 x_{i,2} & \cdots & \sum_{i=1}^2x_{i,n} \\
        \vdots & \cdots & \cdots & \vdots \\
        \sum_{i=1}^n x_{i,1} & \sum_{i=1}^nx_{i,2} & \cdots & \sum_{i=1}^nx_{i,n}
    \end{pmatrix}$$ and $$xa= \begin{pmatrix}
        \sum_{i=1}^nx_{1,i} & \sum_{i=2}^nx_{1,i} & \cdots & x_{1,n} \\
        \sum_{i=1}^nx_{2,i} & \sum_{i=2}^nx_{2,i} & \cdots & x_{2,n} \\
        \vdots & \cdots & \cdots & \vdots \\
        \sum_{i=1}^n x_{n,i} & \sum_{i=2}^n x_{n,i} & \cdots & x_{n,n}
    \end{pmatrix}$$
    Therefore $$\sum_{j = 1}^i x_{j,k} - \sum_{j = k}^n x_{i,j} = 0$$ for $i=1,\ldots,n$ and $k= 1,\ldots,n$. If $i=k=1$, then $x_{1,1}=\sum_{j=1}^nx_{1,j}$ and so~$\sum_{j=2}^nx_{1,j}=0$. Consequently, $x_{1,2}=0$ by choosing $i=1$ and $k=2$. Continuing in this way, we obtain \hbox{$x_{1,2}=\ldots=x_{1,n}=0$.} In particular, the second row of the matrix $ax$ is $(x_{1,1}+x_{2,1}\;\, x_{2,2}\;\,\ldots\;\, x_{2,n})$. Arguing as for the first row, we obtain that $x_{2,3}=\ldots=x_{2,n}=0$. Now, the second row of the matrix $xa$ becomes $(x_{2,1}+x_{2,2}\;\, x_{2,2}\;\,\ldots\;\,0)$, so $x_{1,1}=x_{2,2}$. Repeating the first half of the procedure on the third row gives $x_{3,4}=\ldots=x_{3,n}=0$, so the equality in the third row is $$
    \begin{array}{c}
    (x_{1,1}+x_{2,1}+x_{3,1}\;\;x_{2,2}+x_{3,2}\;\; x_{3,3}\;\; 0\;\;\ldots\;\; 0)\qquad\qquad\qquad\qquad\\[0.2cm]
    \qquad\qquad\qquad=(x_{3,1}+x_{3,2}+x_{3,3}\;\; x_{3,2}+x_{3,3}\;\; x_{3,3}\;\;0\;\;\ldots\;\;0),
    \end{array}$$ and hence $x_{3,3}=x_{2,2}=x_{1,1}$ and $x_{2,1}=x_{3,2}$. Now, an easy induction argument gives $$x=\begin{pmatrix}
    x_{1,1} & 0 &0& \cdots &  0 &0 \\
    x_{2,1} & x_{1,1} & 0 & \cdots &  0 &0 \\
    x_{3,1} & x_{2,1} & x_{1,1} & \cdots  &0 &
    0 \\
    \vdots & \ddots &  \ddots & \ddots & \ddots & \vdots \\
    x_{n-1,1} & \ddots & \ddots & x_{2,1} &x_{1,1} &  0\\
    x_{n,1} &x_{n-1,1} & \cdots & x_{3,1} & x_{2,1} & x_{1,1}
\end{pmatrix}$$



Thus, the $C_G(a)=Z(G)\cdot K$, where $K$ is the subgroup of all  lower unitriangular matrices whose elements of each $i$-lower diagonal coincide for $i=1,\ldots,n$. Since~$K$ is a non-central subgroup of $S$, $C_G(a)^G =Z(G)\cdot S$ (see \cite{acourse}, Theorem 3.2.8). Therefore $P(G)=S\cdot Z(G)$. 

\medskip

\noindent(2)\quad If $g\notin Z(S)$, then $C_S(g)\not\leq Z(S)$ and so $C_S(g)^S = S$ by Theorem 3.2.8 of~\cite{acourse}. 

\medskip

\noindent(3)\quad This follows at once from  Theorem 3.2.8 of \cite{acourse}.

\medskip

\noindent(4)\quad By Lemma \ref{homomorphicimages} and (1), we have  $P(G/Z(G))\geq SZ(G)/Z(G)$. Again, set $$a=\begin{pmatrix}
      1 & 0 & \cdots & 0 \\
      1 & 1 & \cdots & 0\\ 
     \vdots & \vdots & \ddots & \vdots \\
     1 & 1 & \ldots & 1
    \end{pmatrix}$$ and let $x=(x_{i,j})\in G$ be such that $ax=d\cdot (xa)$ for some $d\in\mathbb F^\times$. Then $$\sum_{j = 1}^i x_{j,k} - d\cdot \sum_{j = k}^n x_{i,j} = 0$$ for $i=1,\ldots,n$ and $k= 1,\ldots,n$. If $i=1$ and $d\neq1$, then the above equations become $$x_{1,k} - d\cdot\sum_{j = k}^nx_{1,j} = 0$$ for $k = 1,\ldots,n$. This gives $x_{1,n}=x_{1,n-1}=\ldots=x_{1,1}=0$ and so a contradiction. Thus, $d=1$ and consequently $$C_G(aZ(G))=C_G(a)/Z(G)\leq S/Z(G).$$ Therefore $P(G/Z(G))=SZ(G)/Z(G)$.

\medskip

\noindent(5)\quad It follows from (1), (2) and Theorem \ref{Bern1} that $$P(P(G))=P\big(S\cdot Z(G)\big)=P(S)\cdot Z(G)=S\cdot Z(G)=P(G).$$ The statement is proved. 
\end{proof}

\medskip

The group of triangular matrices seems to be more difficult to deal with. Here we only describe the pseudocentre of the triangular matrices of degree $2$.

\begin{theorem}
Let $G =\operatorname{Tr}(2,\mathbb{F})$ be the group of all triangular matrices over the field $\mathbb F$. Then $P(G) = Z(G) \times \operatorname{Tr}_1(2,\mathbb F)$.
\end{theorem}
\begin{proof}
It is well-known that $U=\operatorname{Tr}_1(2,\mathbb F)$ is minimal normal in $G$. Thus, $P(G)\geq Z(G)\times U$. On the other hand, the centralizer of any non-trivial element of $U$ is contained in $Z(G)\times U$ and so $P(G)=Z(G)\times U$.
\end{proof}

\medskip

In the final part of this section, we deal with the pseudocentre of the affine general and special linear groups.

\begin{theorem}\label{affinegeneral}
Let $n$ be a positive integer and $\mathbb{F}$ a field such that either $n = 2$ and $|\mathbb{F}| > 3$, or $n \ge 3$. Let $V$ be the vector space of dimension $n$ over $\mathbb F$.\linebreak Set~$S=\operatorname{SL}(n,\mathbb F)$, $G = \operatorname{GL}(n,\mathbb{F})$ and $H = G \ltimes V$, where $G$ acts in the natural way on~$V$. Then $P(H)=P(S\ltimes V)=S\ltimes V$.
\end{theorem}
\begin{proof}
Let $R=S\ltimes V$. By using the upper and lower unitriangular matrices, we easily see that~$V$ is a minimal normal subgroup of $R$. Thus, $V\leq P(R)$ by~The\-o\-rem~\ref{normsempl}. Let~$v=(v_i)_{i=1,\ldots,n}\in V\setminus\{1\}$ and let $s\in Z(S)$. By suitably conjugating $sv$ by elements of $S$, we may assume that $v_1=\ldots=v_{n-1}=0$ and $v_n\neq0$. But then $sv$ is clearly centralized by an element of~$S\setminus Z(S)$. Consequently, $S\leq C_R(sv)^R$ (see \cite{acourse},~The\-o\-rem~3.2.8) and $C_R(sv)^R=R$. 

Now, let $g\in R\setminus\{1\}$, and write $g=sv$, where $s\in S\setminus Z(S)$ and $v\in V$. Since $$v\in V\leq P(R)\leq C_R(g)^R,$$ so $s\in C_R(g)^R$. Then~\hbox{$S\leq C_R(g)^R$} by Theorem 3.2.8 of \cite{acourse}, so again $R=C_R(g)^R$. Therefore $P(R)=R$.

Finally, we compute the pseudocentre of $H$. As before, $V$ is minimal normal in $H$, so $V\leq P(H)$. Let $b = (1,\ldots,1)\in V$ and $$a=\begin{pmatrix}
      1 & 0 & \cdots & 0 \\
      1 & 1 & \cdots & 0\\ 
     \vdots & \vdots & \ddots & \vdots \\
     1 & 1 & \ldots & 1
    \end{pmatrix}.$$ If $cd\in C_H(ab)$, where $c\in G$ and $d\in V$, then $c\in C_G(a)$ and $b^cd=d^ab$. Now,~$C_G(a)$ has been computed in the proof of Theorem \ref{GL}, and it is the set of all matrices whose elements of each lower diagonal coincide. But then the first equation deduced from~\hbox{$b^cd=d^ab$} yields that the main diagonal of $c$ consists of $1$s, so $c\in S$. Consequently, $P(H)\leq SV$. Since the same argument as in the second paragraph shows that $P(H)\geq S$, so $P(H)=S\ltimes V$.
\end{proof}

\begin{theorem}
\begin{itemize}
    \item[(1)] Let $G = \operatorname{GL}(2,2) \ltimes V$, where $V$ is the vector space of dimension $2$ over the field of order $2$. Then $P(G) = \operatorname{Alt}(3) \ltimes V$.
    \item[(2)] Let $G = \operatorname{SL}(2,3) \ltimes V$, where $V$ is the vector space of dimension $2$ over the field of order $3$. Then $P(G) = O_2(\operatorname{SL}(2,3)) \ltimes V$; 
    \item[(3)] Let $G = \operatorname{GL}(2,3) \ltimes V$, where $V$ is the vector space of dimension $2$ over the field of order $3$. Then $P(G) = \operatorname{SL}(2,3) \ltimes V$.
\end{itemize}
\end{theorem}
\begin{proof}
(1)\quad We just need to note that $G\simeq\operatorname{Sym}(4)$, so the result follows from Theorem \ref{thsymalg}.

    \medskip
(2)\quad Now, $V$ is a minimal normal subgroup of $G_2$, so is contained in $P(G)$ by~The\-o\-rem~\ref{normsempl}. Let $g = xb$, where $x\in S=\operatorname{SL}(2,3)$ and $b\in V$. Since $V\leq P(G)$, so $x\in C_G(g)^G$. Now, recall that $S/Z(S)\simeq\operatorname{Alt}(4)$, and set $Q = O_2(S)$. Then we easily see that if~\hbox{$x\not\in Z(S)$,} then $Q\leq \langle x\rangle^G\leq C_G(g)^G$. On the other hand, if~\hbox{$x\in Z(S)$,} then we can argue as in the first paragraph of Theorem \ref{affinegeneral} and have that~$x$ is centralized by an element of $S\setminus Z(S)$, which again means that~\hbox{$Q\leq C_G(g)^G$.} Therefore $QV\leq P(G)$. Finally, a combination of Lemma \ref{homomorphicimages} and Theorem \ref{GL23} yields that $P(G)\leq QV$, thus completing the proof of this case.

\medskip

\noindent(3)\quad Again, $V$ is a minimal normal subgroup of $G$, while by Theorem \ref{GL23}, $P(G) \le S \ltimes V$, where $S=\operatorname{SL}(2,3)$. Let $g = xb$, where $x\in H=\operatorname{GL}(2,3)$ and \hbox{$b=(b_1,b_2)\in V$.} Since~\hbox{$V\leq P(G)$,} so $x\in C_G(g)^G$. If $x\in H\setminus O_2(H)$, then $S\leq \langle x\rangle^G\leq C_G(g)^G$. Assume $$x\in O_2(H)=\left\{\pm\begin{pmatrix}
        0 & 1 \\
        -1 & 0
    \end{pmatrix},\pm\begin{pmatrix}
        -1 & -1 \\
        -1 & 1
    \end{pmatrix},\pm\begin{pmatrix}
        -1 & 1 \\
         1 & 1
    \end{pmatrix}\right\}.$$ We do a case-by-case analysis to show that $C_G(g)$ contains an element of the form~$yc$, where $y\not\in O_2(H)$ and $c\in V$. If  $x = \begin{pmatrix}
    0 & 1 \\
    -1 & 0
\end{pmatrix}$, then we put $y = \begin{pmatrix}
    1 & 1 \\
    -1 & 1
\end{pmatrix}$ and $c = (b_2 - b_1, -b_1 - b_2)$. If $x = \begin{pmatrix}
    0 & -1 \\
    1 & 0
\end{pmatrix}$, then we put $y = \begin{pmatrix}
    1 & 1 \\
    -1 & 1
\end{pmatrix}$ and $c = (b_1+b_2,b_2 -b_1)$. If $x = \begin{pmatrix}
    -1 & -1 \\
    -1 & 1
\end{pmatrix}$, then $y = \begin{pmatrix}
    0 & 1 \\
    1 & 1
\end{pmatrix}$ and $c = (-b_1,-b_2)$. If~$x = \begin{pmatrix}
    1 & 1 \\
    1 & -1
\end{pmatrix}$, then $y = \begin{pmatrix}
    0 & 1 \\
    1 & 1
\end{pmatrix}$ and $c = (-b_1 -b_2,b_2 -b_1)$. If $x = \begin{pmatrix}
    -1 & 1 \\
    1 & 1
\end{pmatrix}$, then we put $y=\begin{pmatrix}
    0 & 1 \\
    1 & -1
\end{pmatrix}$ and $c =(b_2-b_1, b_1+b_2)$. If $x = \begin{pmatrix}
    1 & -1 \\
    -1 & -1
\end{pmatrix}$, then we put $y = \begin{pmatrix}
    0 & -1 \\
    -1 & 1
\end{pmatrix}$ and $c = (b_2,b_1-b_2)$. Now, since $y\in C_G(g)^G$, so $S\leq C_G(g)^G$, which proves that $P(G)=S\ltimes V$.
\end{proof}

\chapter{Trees and weakly branch groups}\label{treesect}

We start by recalling some basic definitions and results about rooted trees; this is mostly taken from \cite{Grig2} and \cite{Grig1}. Let $T$ be a {\it rooted} tree, that is, a tree with a distinguished {\it root} vertex $r$. The set of vertices of the tree is $V(T)$. The {\it norm} of a vertex $u$ is the number of edges in a geodesic path from $r$ to $u$. Clearly, $\left|u\right| = 0$ if and only if $u=r$. The cardinality of the set of vertices that are both adjacent to~$u$ and of norm $\left|u\right|+1$ is the \textit{degree} of $u$. The tree is \textit{spherically homogeneous} if the vertices with the same norm have the same degree.
The {\it$n$-th level}~$V_n$ of $T$ is the set of vertices of $T$ with norm $n$. A spherically homogeneous tree is determined by its {\it branch index}, that is, a sequence $\overline{m} = \{m_n\}_{n=1}^N$, where $N$ can be finite or infinite, and~$m_n$ is the degree of a vertex at level $n$. If~\hbox{$m_n = d$} for every $n$, then $T$ is said to be the {\it rooted regular tree of degree~$d$}; if $d = 2$, then $T$ is called the \textit{binary tree}. Note that in the following we will usually assume $m_n \ge 2$ for each~\hbox{$n \in \{1,\ldots,N\}$.}

A bijective map $\varphi: V(T) \to V(T)$ preserving edge incidence and the root vertex is said to be an \textit{automorphism} of $T$. The set of all automorphisms of $T$ is denoted by $\operatorname{Aut}(T)$ and it is a group with the composition.
If $G \le\operatorname{Aut}(T)$ and $u \in V(T)$, then $\operatorname{Stab}_G(u) = \{g \in G : g(u) = u\}$ is  \textit{stabilizer} of $u$ in $G$. If~\hbox{$n\in\{1,\ldots,N\}$,} then the \textit{$n$-level stabilizer of $G$} is $$\operatorname{Stab}_G(n) = \bigcap_{\left| u\right| = n}\operatorname{Stab}_G(u).$$
Let $T[n]$ be the subtree of $T$ consisting of the vertices of norm at most $n$. Clearly, if $\varphi\in\operatorname{Aut}(T)$, then $\varphi|_{T[n]}\in\operatorname{Aut}(T[n])$ and  $\operatorname{Aut}(T)/\operatorname{Stab}_{\operatorname{Aut}(T)}(n) \simeq \operatorname{Aut}(T[n])$; actually, it turns out that $\operatorname{Aut}(T)= K\ltimes \operatorname{Stab}_{\operatorname{Aut}(T)}(n)$ for some obvious subgroup $K\simeq\operatorname{Aut}(T[n])$. If $\varphi\in G\leq\operatorname{Aut}(T)$, $u_1,\ldots,u_{m_n}$ are the vertices of $T$ at level $n$, and $\varphi_i = \varphi|_{T_{u_i}}$, where $T_{u_i}$ is the subtree of $T$ with root $u_i$, then the map $$\tau_n: \varphi\in\operatorname{Stab}_G(n)\mapsto (\varphi_1,\ldots,\varphi_{m_n})\in\operatorname{Aut}(T_{u_1}) \times \ldots \times \operatorname{Aut}(T_{u_{m_n}}),$$ is a monomorphism; if $G = \operatorname{Aut}(T)$, then $\tau_n$ is an isomorphism; it is clear that if $T$ is a spherically homogeneous, then the $T_{u_i}$'s are isomorphic as trees, and so we can simply write~$T_n$. Recall also that $\{\operatorname{Aut}(T[n]),\, \theta_{m,n}\}_{m\geq n\in N}$ is an inverse system with $$\theta_{m,n}:\operatorname{Aut}(T[m]) \to \operatorname{Aut}(T[n])$$ defined by restriction of the action of~$\operatorname{Aut}(T[m])$ to~$\operatorname{Aut}(T[n])$, and that the group $\operatorname{Aut}(T)$ is the inverse limit of this inverse system. It is also well-known that if $T$ is a finite spherically homogeneous rooted tree with end level $n$, then $$\operatorname{Aut}(T[n]) \simeq \big(\ldots\big(\big(\operatorname{Sym}(m_n)\,\operatorname{\it wr}\, \operatorname{Sym}(m_{n-1})\big)\,\operatorname{\it wr}\, \ldots \big)\,\operatorname{\it wr}\, \operatorname{Sym}(m_1).$$ Thus, if $T$ is a infinite, then $$\operatorname{Aut}(T) \simeq \varprojlim \big(\ldots\big(\big(\operatorname{Sym}(m_n )\,\operatorname{\it wr}\, \operatorname{Sym}(m_{n-1})\big)\,\operatorname{\it wr}\, \ldots \big)\,\operatorname{\it wr}\, \operatorname{Sym}(m_1).$$

\medskip





\smallskip

Now, we are in a position to introduce the main definition of this section, and to prove our main result. Let $T$ be a rooted tree, and $G$ a subgroup of $\operatorname{Aut}(T)$. If $v \in V(T)$ then the \textit{rigid vertex stabilizer}~$\operatorname{rst}_G(v)$ of $v$ is the subgroup of all elements $g$ of $G$ stabilizing every vertex \hbox{$u \in V(T)\setminus V(T_v)$;} note that $\operatorname{rst}_G(v)\leq\operatorname{Stab}(|v|)$. The group~$G$ is said to be \textit{weakly branch} if it is {\it level-transitive} (i.e. it acts transitively on every level) and if $\operatorname{rst}_G(v) \neq \{1\}$ for every $v \in V(T)$. Clearly, if $G$ is weakly branch, then $T$ must be spherically homogeneous.

\begin{lemma}\label{CentralizeAut(T)}
    Let $K$ be a transitive permutation group on a finite non-empty set~$X$, $H$ a non-trivial group, and~\hbox{$G = H\,\operatorname{\it wr}\, K$}. If $b=\prod_{x \in X}b_x$ is an element of the base group $B=\operatorname{Dr}_{x\in X}H_x$ of $G$, where $b_x\in H_x\simeq H$ for $x\in X$, and $b_x$ is not conjugate to $b_y$ in $G$ for $x\neq y$, then $C_G(b) = \{ \prod_{x\in X}a_x : a_x \in C_{H_x}(b_x)\}$.
\end{lemma}
\begin{proof}
Let $g \in C_G(b)$ and write $g = ka$, where $k \in K$ and $a=\prod_{x \in X}a_x$ with $a_x\in H_x$. Now, $$b = b^g = \prod_{x \in X}b_{k(x)}^{a_{k(x)}}$$ implies $k=1$ and $a_x\in C_G(b_x)$.
\end{proof}

\begin{theorem}\label{thpetto}
Let $G$ be a weakly branch group on a rooted tree $T$. Then $P(G) = \{1\}$.\footnote{Our thanks go to Jan Moritz Petschick for his insightful discussions on the topic, which have significantly enhanced the exposition of this chapter.}
\end{theorem}
\begin{proof}
Let $(v_i)_{i \in \mathbb{N}}$ be a sequence of vertices of $T$ such that $v_{i+1}$ is a {\it children} of~$v_i$ (that is,~\hbox{$(v_{i+1},v_i)$} is an edge, and $v_{i+1}$ belongs to the $(|v_{i}|+1)$-th level). Put $u_0=v_0$, and suppose $u_i$ is defined as a certain $v_{\ell_i}$. Since $\operatorname{rst}_G(u_i)\neq\{1\}$ by definition, it contains a non-trivial element~$g$ which must act non-trivially on some level of $T$ which is strictly larger that $|u_i|$. Thus, there is~\hbox{$\ell_{i+1}\in\mathbb N$} such that $\ell_i<\ell_{i+1}$ and $v_{\ell_{i+1}}^g\neq v_{\ell_{i+1}}$. Put $u_{i+1}=v_{\ell_{i+1}}$. This procedure defines a subsequence $(u_i)_{i\in\mathbb N}$ of $(v_i)_{i \in \mathbb{N}}$ such that $\operatorname{rst}_G(u_{i})   \not\leq\operatorname{Stab}(|u_{i+1}|)$ for every $i\in\mathbb N$.

Let $i,n \in \mathbb{N}$ such that $|u_i|=n$. Say that the $n$-th level has cardinality $m$, and write $V_n=\{u_i^0,u_i^1,\ldots,u_i^{m-1}\}$. Let $x_j \in \operatorname{rst}_G(u_{i+j}) \setminus \operatorname{Stab}(|u_{i+j+1}|)$ for $j\in\{0, 1, \dots, m-1\}$. These elements are not conjugate in $\operatorname{Aut}(T)$ by construction. For each $j \in\{0,1, \ldots,m-1\}$, let $g_j$ be an element of $G$ mapping~$u_i$ to $u_i^j$. Set $z=x_0^{g_0}\ldots x_{m-1}^{g_{m-1}}\in\operatorname{Stab}(n)$. By Lemma \ref{CentralizeAut(T)}, the centralizer of $z$ in~$\operatorname{Aut}(T)$ is contained in $\operatorname{Stab}(n)$, so $C_G(z)\leq\operatorname{Stab}(n)$, and a fortiori $C_G(z)^G\leq\operatorname{Stab}(n)$.  Thus, $P(G)\leq\bigcap_{n\in\mathbb N}\operatorname{Stab}(n)=\{1\}$.
\end{proof}

\medskip\medskip

Finally, we briefly review which groups are covered by~The\-o\-rem~\ref{thpetto}: 

\begin{itemize}

    \item $\operatorname{Aut}(T)$, where $T$ is spherically homogeneous. Since $\operatorname{Aut}(T)$ is a profinite group, this shows that the non-triviality of the pseudocentre of finite groups cannot be extended to profinite groups.

    \item The branch groups as defined in \cite{PET1} and \cite{PET2}.

    \item The Gri\-gor\-chuk $p$-group \cite{Grig2},\cite{Grig1} and, for an odd prime $p$, the~\hbox{Gupta--Sid}\-ki $p$-group as defined in \cite{Gup}.

    \item The generalization of Grigorchuk $2$-groups given by \v Suni\'c in \cite{Sun} (see~De\-finition 2, Lemma 1 and Lemma 6).

    \item The Basilica groups, the generalized Brunner--Sidki--Vieira groups (see~\cite{Ele} and \cite{Far}), and the very general Basilica group provided in \cite{PET3}.

    \item Some of the Multi Extended Gupta--Sidki groups~$G_E$ (see \cite{Klo}, and Lemmas 3.3 and 3.4 of \cite{Anita}). 

    \item Groups acting on non-regular trees (see for example \cite{PET4}, \cite{PET5}, and \cite{PET6}).

    \item The {\it $p$-adic tree automorphism group} $\mathcal{A}_{p}$, which is the subgroup of the tree automorphism group where every permutation induced on a vertex belongs to~$\langle (1,\ldots,p)\rangle$.
\end{itemize}

\begin{remark}
{\rm The Grigorchuk groups are just-infinite, so in particular the normal closures of the centralizers have finite index in the group. Nevertheless, its pseudocentre is trivial.}
\end{remark}

\chapter{Thompson-like groups}

The Thompson groups were first introduced by Richard J. Thompson in an unpublished manuscript from the 1960s (see \cite{mckenzie}). Thompson’s original aim was to apply what is now called the group $V$ (in fact a larger monoid of which $V$ is the group of units) to questions in algebraic logic. Although this application never materialized, he soon realised that $V$ is an infinite, finitely presented, and simple group --- an unexpected combination of properties that made it a striking discovery (see also \cite{thompsonthing}). The group $F$, introduced in the same context, later became central in the study of the {\it von Neumann conjecture} (a group is non-amenable if and only if it contains a countable free subgroup): in fact, $F$ is a torsion-free, finitely presented simple group containing no free subgroups of rank greater than $1$. Since then, the Thompson group $F$ has been rediscovered in many different contexts and has attracted a lot of attention for its remarkable properties. The wider group-theoretical context in which $F$ is studied is that of $\operatorname{PL_o}(\mathbb{R})$ (or~$\operatorname{PLF}(\mathbb{R})$), namely the group of all orientation-preserving piecewise linear homeomorphisms of the real line with only finitely many points of non-differentiability. In this chapter, we study the pseudocentre of a family of relevant subgroups of $\operatorname{PL_o}(\mathbb{R})$, including $F$ itself.

Following \cite{Bier} (see also \cite{BrinBrin} and \cite{Brin1}), let $I$ be an interval of the real line, $P$ a subgroup of the multiplicative group of the positive real numbers, and $A$ a $\mathbb{Z}[P]$-submodule of $(\mathbb{R},+)$. For any piecewise linear function $f$, we define the \textit{slopes} of $f$ as the set of all slopes of its linear components, the \textit{breaks} as the set of all points where $f$ is non-differentiable, and the \textit{support} of $f$ as the set $\operatorname{supp}(f)=\{x\in I\mid f(x)\neq x\}$. Then we define~$G(I;A,P)$ to be the subgroup of $\operatorname{PL_o}(\mathbb{R})$ with $\operatorname{supp}\!f$ in $I$, slopes in~$P$ and breaks in $A$, with the further condition that, if $I=\mathbb R$, then every function of $G(I;A,P)$ has to map $A$ onto itself.

\smallskip

We begin recalling some properties of the Thompson group $F$, which is defined here as $G(I;A,P)$, where $I=[0,1]$, $A=\mathbb Z[1/2]$ and $P=\{2^n\mid n\in \mathbb{Z}\}$.  It is known that $F$ is generated by the maps
$$
A(x) =  \begin{cases}
    \frac{x}{2},& 0 \le x \le \frac{1}{2} \\ 
    x - \frac{1}{4}, & \frac{1}{2} \le x \le \frac{3}{4} \\
    2x - 1 & \frac{3}{4} \le x \le 1
\end{cases} \qquad\textnormal{and}\quad
\quad B(x) = \begin{cases}
    x, & 0 \le x \le \frac{1}{2}\\
    \frac{x}{2} + \frac{1}{4}, & \frac{1}{2} \le x\leq \frac{3}{4}\\
    x - \frac{1}{8}, & \frac{3}{4} \le x \le \frac{7}{8} \\
    2x - 1 & \frac{7}{8} \le x \le 1 
\end{cases}$$ (see \cite{Cann1}, Corollary 2.6). Moreover, $F'$ is a simple group (see \cite{Cann1}, Theorem 4.5),  $F/F'$ is a free abelian group of rank $2$ (see \cite{Cann1}, Theorem 4.1), and every proper quotient of $F$ is abelian (see \cite{Cann1},~The\-o\-rem~4.3). Let $f\in F$ and let $D$ be a closed interval of~$[0,1]$. Then $D$ is said to be a \textit{bump domain} for $f$ if the only fixed points of $f$ in $D$ are the first and the end points of $D$. A sequence of intervals $\mathcal{D}$ of $[0,1]$ is a \textit{bump chain} for $f$ if every interval $D \in \mathcal{D}$ is a bump domain for $f$. Of course, there exists a {\it partition} of $[0,1]$ (that is, a set of subsets of $[0,1]$ with pairwise disjoint interiors whose union is $[0,1]$) consisting of a bump chain of $f$ and maximal connected open sets of fixed points of $f$. Moreover, $f$ is said to be a \textit{one bump function} if~$[0,1]$ is a bump domain for~$f$. 

\medskip
\begin{lemma}\label{divisibleroot}
Let $f$ and $g$ be one bump functions of $F$, and let $n$ be a positive integer. 
\begin{itemize}
    \item The centralizer $C_F(f)$ of $f$ in $F$ is infinite cyclic.
    \item If $f$ and $g$ are conjugate, then their initial and final slopes coincide.
    \item If $g^n = f$, then $n$ divides the base $2$ logarithm of the first and final slope of $f$.
\end{itemize}
\end{lemma}
\begin{proof}
The first point follows from \cite{Gill}, Theorem 7.2, while the second is a trivial computation. Suppose now that $g^n=f$.  There exists a positive rational number $\epsilon$ such that $f$ and $g$ are both piecewise on $[0,\epsilon]$; let $a_1$ and $a_2$ be the slopes of $f$ and $g$ in $[0,\epsilon]$, respectively. Then $a_2^n=a_1$. Since~$a_1$ and $a_2$ are powers of $2$, so~$n$ divides the base $2$ logarithm of $a_1$. A similar argument works for the final slope of $f$.
\end{proof}

\begin{lemma}\label{magicfunctions}
For all positive integers $m,n$, there exists one bump functions $f_{m,n}$ and $g_{m,n}$ such that:
\begin{itemize}
    \item the initial slope of $f_{m,n}$ is $2^m$, and the final slope of $f_{m,n}$ is $2^{-n}$;
    \item the initial slope of $g_{m,n}$ is $2^{-m}$, and the final slope of $g_{m,n}$ is $2^n$.
\end{itemize}
\end{lemma}
\begin{proof}
For each $f\in F$, we need the following two auxiliary functions: $$\theta_f:\, x\in \left[0,\frac{1}{2}\right]\mapsto \frac{1}{4}f(2x)\in \left[0,\frac{1}{4}\right]$$  and 
$$\rho_f:\, x\in \left[\frac{3}{4},1\right]\mapsto \frac{1}{2}\big(f(4x-3)\big) + \frac{1}{2}\in \left[\frac{1}{2},1\right]$$ Note that if $f(x)<x$ for every $x\in (0,1)$, then $\rho_f(y)<y$ for every $y\in [3/4,1)$: in fact, $\rho_f(y)=1/2f(4y-3)+1/2<1/2(4y-3)+1/2=2y-1< y$. Similarly, if $f(x)<x$ for~\hbox{$x\in(0,1)$,} then $\theta_f(y)<y$ for $y\in(0,1/2]$. Moreover, $\theta_f(0)=0$, $\theta_f(1/2)=1/4=A(1/2)$, $\rho_f(1)=1$, and $\rho_f(3/4)=1/2=A(3/4)$.

Since the function $A(x)$ is a one bump function, we can put $g_{1,1}(x)=A(x)$. Now, suppose $g_{m,1}$ and $g_{1,n}$ have been defined and let $$g_{1,n+1}(x)=\begin{cases}
    A(x), & 0 \le x \le \frac{3}{4} \\
    \rho_{g_{1,n}}(x), & \frac{3}{4} \le x \le 1
\end{cases}\quad\textnormal{and}\quad g_{m+1,1}(x)=\begin{cases}
    \theta_{g_{_{m,1}}}(x), & 0 \le x \le \frac{1}{2} \\
    A(x), & \frac{1}{2} \le x \le 1
\end{cases}$$ Since $A(x)<x$ for $x\in(0,1)$, so $g_{1,n+1}$ and $g_{m+1,1}$ are one bump functions whose initial and final slopes are as in the statement. Finally, we define $$
g_{m,n}(x)=\begin{cases}
    g_{m,1}(x), & 0\leq x\leq \frac{1}{2}\\
    A(x), & \frac{1}{2}\le x \le \frac{3}{4} \\
    g_{1,n}(x), & \frac{3}{4} \le x \le 1
\end{cases}
$$ and we see that its final slope is $2^n$, while its initial slope is $2^{-m}$. In order to complete the proof, it is enough to define $f_{m,n}=g_{m,n}^{-1}$.
\end{proof}

\begin{theorem}\label{thompson}
The pseudocentre of the Thompson group $F$ is the commutator subgroup $F'$.
\end{theorem}
\begin{proof}
Let $ f \in F$. By Theorem 4.3 in \cite{Cann1}, $C_F(f)^F \ge F'$. Therefore, \hbox{$P(F) \ge F'$.} Let $p$ and $q$ be distinct prime numbers. By Lemma \ref{magicfunctions}, there exists $f_{p,q}\in F$ such that the initial and final slopes are $2^p$ and $2^{-q}$, respectively. By~Lem\-ma~\ref{divisibleroot}, $C_F(f_{p,q})$ is infinite cyclic, so we may find a generator $g$ of $C_F(f_{p,q})$ such that $g^n=f$ for a positive integer $n$. However, the same lemma yields that $n$ divides both $p$ and $q$, so $n=\pm1$, and $C_F(f_{p,q})=\langle f_{p,q}\rangle$.

Now, the map $\varphi: F \to \mathbb{Z} \times \mathbb{Z}$ assigning to each $f \in F$ the ordered pair $\big(\log_2(m_f),\log_2(n_f)\big)$, where $m_f$ and $n_f$ are respectively the initial and final slopes of $f$, is an epimorphism with kernel $F'$ (see the proof of Theorem 4.1 of \cite{Cann1}). It follows that $$P(F)\leq C_F(f_{2,3})^F\cap C_F(f_{5,7})^F\leq \langle f_{2,3}\rangle F'\cap\langle f_{5,7}\rangle F'=F',$$ and hence $P(F)=F'$, as claimed.
\end{proof}

\medskip

Now that we are done with the Thompson group, we are going to deal with a family of larger subgroups of $\operatorname{PL_o}(\mathbb{R})$.

In \cite{Brin}, the structure of the centralizers of the elements of $G(S;A,P)$ is provided. These results will be used many times in the next part of this section, so, for the reader's convenience, we recall some basic concepts and notation that will be useful in the following.  Let $g \in G(S;A,P)$ and $x\in S$. The \textit{indicator function} $i_g$ is defined as follows: $$i_{g} : s\in S \mapsto i_{g}(s) = \begin{cases}
    1 & g(s) > s\\
    0 & g(s) = s \\
    -1 & g(s) < s
\end{cases}$$ The \textit{orbit} of $g$ containing $x$ is \[\operatorname{orb}(g,x) = \{ g^{n}(x)\}_{n \in \mathbb{Z}},\] while the \textit{orbital} of $g$ containing $x$ is \[\operatorname{orbl}(g,x) = \underset{n \in \mathbb{Z}}{\bigcup}[g^{n}(x),g^{n+i_{g}(x)}(x)].\] Also, recall that $B(S;A,P)$ is the subgroup of $G(S;A,P)$ consisting of all functions that are the identity at the beginning and at the end (see page iv of \cite{Bier}). 

\begin{lemma}\label{B}
Let $I$ be an infinite interval of the real line, $P$ a non-trivial subgroup of $(\mathbb{R}^+,\cdot)$ containing a rational number different from $1$, and $A$ a $\mathbb{Z}[P]$-module of $(\mathbb{R},+)$. If $A$ is divisible, then $B=B(I;A,P)$ is a non-abelian simple group.
\end{lemma}
\begin{proof}
Since the action of $P$ on $A$ is fixed-point-free, the $0$-th homology group $H_0(P,A)$ is trivial, so we may apply Theorem A of \cite{Rob76} to obtain that \hbox{$H_1(P,A)=\{0\}$.} Moreover, since~$P$ contains a rational number different from $1$ and $A$ is divisible, the submodule $IP\cdot A$ of $A$ defined at page iii of \cite{Bier} equals $A$. Then, Theorem 2.1 of \cite{Bier} yields that~$B$ is perfect and hence~$B$ is a non-abelian simple group by Corollary C10.3 of \cite{Bier} and because~$B$ is clearly not trivial.~
\end{proof}

\begin{theorem}\label{PLF0l}
Let $F$ be a subfield of $\mathbb{R}$, and set $A=(F,+)$ and $P=(F^+,\cdot)$. Let $a$ and $b$ be real numbers contained in $A$. Then the pseudocentre of $G=G((a,b);A,P)$ is its commutator subgroup.
\end{theorem}
\begin{proof}
Clearly, we can transform $G$ via the affine mapping $x\in\mathbb{R}\mapsto x-a\in\mathbb{R}$ and suppose that $G=G((0,l);A,P)$ with $l=b-a$. Let $\varphi: G \to P \times P$ be the epimorphism mapping $g \in G$ into the ordered pair $(a,b)$, where $a$ is the initial slope of $g$ and $b$ the final one in~$(0,l)$, and notice that $\operatorname{Ker}(\varphi)$ is exactly the group $B\big((0,l);A,P\big)$ defined at page~iv of \cite{Bier}, call it~$B$. Then $B$ is a non-abelian simple group by Lemma \ref{B}, and so $B=G'$, which is hence contained in $P(G)$ by~The\-o\-rem~\ref{normsempl}.

Choose $s\in (0,l)\cap A$ and $1<a\in A$, and define
\begin{align*}
\gamma=\gamma(s,a)(x) =  \begin{cases}
    \frac{ax}{1+l^{-1}s(a-1)},& 0 < x \le s \\[0.3cm] 
    \frac{x+s(a-1)}{1+l^{-1}s(a-1)}, & s \le x < l \\
\end{cases} &
\end{align*}

Now, if we write $$b=\frac{1}{1+l^{-1}s(a-1)}\quad\textnormal{ and }\quad c=\frac{s(a-1)}{1+l^{-1}s(a-1)},$$ we have that, for any $x\in(s,l)$, $$\gamma^n(x)=b^nx+c\sum_{i=0}^{n-1}b^i$$ so the convergence formula for geometric series yields that $\lim_{n\to +\infty}\gamma^{n}(x)=l$. On the other hand, for any $x$ in $(0,s)\subseteq\big(0,(ab)s\big)$, we have that $\gamma^{-n}(x)=d^nx$, where $d=\frac{1+l^{-1}s(a-1)}{a}=(ab)^{-1}$, so $\lim_{n\to +\infty}\gamma^{-n}(x)=0$ because $a>1>d$. This yields that $\operatorname{orbl}(\gamma,s) = (0,l)$. Now, let $G_0$ be the subgroup of all the affine functions of~$G$. Since $f(0)=0$ and $f(l)=l$ for any element $f$ of $G$, $G_0$ is the trivial subgroup. From this and Lemma 3.17 (b) of \cite{Brin}, it follows that the identity is the only so-called $(0,0)$-conjugator of any $\gamma(s,a)$. Hence we may use Proposition 3.24 of \cite{Brin} to show that $\langle\gamma(s,a)\rangle=C_G(\gamma(s,a))$ for any choice of $s\in(0,l)$ and $a>1$.

Take now into account $\gamma_1=\gamma(\frac{l}{2},2)$ and $\gamma_2=\gamma(\frac{l}{3},2)$ and let $n$ be a non-negative number. Now, $\varphi(\gamma_1)=\big(\frac{4}{3},\frac{2}{3}\big)$ and $\varphi(\gamma_2)=\big(\frac{3}{2},\frac{3}{4}\big)$. Thus, $$\langle \varphi(\gamma_1 G')\rangle\cap \langle \varphi(\gamma_2G')\rangle=\{(0,0)\}$$ and hence $\langle \gamma_1G'\rangle\cap\langle \gamma_2G'\rangle=G'$, which means that \[P(G) \le \langle \gamma_1 \rangle^{G}\cap \langle \gamma_2 \rangle^{G} \le G'.\]
The statement is proved.
\end{proof}

\medskip


\begin{theorem}\label{PLFR+}
Let $F$ be a subfield of $\mathbb{R}$, and set $A=(F,+)$ and $P=(F^+,\cdot)$. Let $a$ be an element of $A$ and let $I$ be either $(-\infty,a)$ or $(a,+\infty)$. Then the pseudocentre of $G=G(I;A,P)$ is its commutator subgroup.
\end{theorem}
\begin{proof}
If $I=(-\infty,a)$, we can transform $G$ via the affine mapping $$x\in\mathbb{R}\mapsto -x+a\in\mathbb{R}$$ and suppose that $I=\mathbb{R}^+$. Let $\varphi: G \to P \times P$ be the epimorphism mapping $g \in G$ to the ordered pair $(a,b)$, where $a$ and $b$ are respectively the initial and final slopes of~$g$, and set~\hbox{$K=\operatorname{Ker}(\varphi)$.} Since any two linear functions with slope~$1$ commute, we also have that $K'$ is contained in $B=B(\mathbb R^+;A,P)$, which is a non-abelian simple group by Lemma \ref{B}. Then $B=K'=G''$ and this is contained in $P(G)$ by~The\-o\-rem~\ref{normsempl}.

We claim that $\langle g\rangle^{G} K' = K$. To this aim, let $\vartheta : K \to A$ be the epimorphism mapping every element of~$K$ to the intercept of its last linear piece, and let $g$ and $g_1$ be elements of~$K\setminus K'$ and~$K$ with final intercepts $k$ and $k_1$, respectively. Without loss of generality, we may assume that $k$ and $k_1$ are both positive. Note also that, using induction and the fact that each function in $G$ has to pass through $0$, it is easy to show that the intercept of any linear piece of a function in $G$ is still an element of $A$. Then, if we choose $h\in G$ with final linear piece given by $k^{-1}k_1x$, then $g^h (x) = x + k_1$ for any large enough $x$. This means that $\vartheta(g^h)=\vartheta(g_1)$, so  $$g^hg_1^{-1}\in\operatorname{Ker}(\vartheta)=K'.$$ The arbitrariness of $g_1$ shows that~\hbox{$\langle g\rangle^{G} K' = K$.} Therefore, $K/K'$ is a minimal normal subgroup of $G/K'$, so in particular $G'= K$ and $\operatorname{Ker}(\vartheta) = G''$.

Let $g\in G$. If $g$ is a non-trivial element of $G''$ which is the identity function from a point~$\overline{x}$ on, then $g$ is centralized by any element of $G'\setminus G''$ which is the identity before $\overline{x}$ and this shows that $G'\leq C_G(g)^G$, since $G''$ is simple and $G'/G''$ is a minimal normal subgroup of $G/G''$. If $g\in G'\setminus G''$, then there is $h\in G'$ such that $[g,h]$ is a non-trivial element of $G''$, and hence $G'\leq\langle g\rangle^G\leq C_G(g)^G$. Finally, let $g$ be an element of $G\setminus G'$. Then the final linear piece of $g$ can be written as  $ax+b$, where $a,b\in A$ with $a\neq1$. Choose~$h$ in $G$ with final linear piece $x+c$ with $c\in A\setminus\{0\}$. Then the final linear piece of $[g,h]$ is $x+c(1-a)$,
so $[g,h]\in (G'\setminus G'')\cap C_G(g)^G$, and, as before, also in this case $G'$ is contained in~$C_G(g)^G$. This shows that $G'\leq P(G)$.

It only remains to show that $P(G)\leq G'$. Consider now the following family $\mathcal{F}$ of functions \begin{align}
\alpha(a,b,s)(x) =  \begin{cases}
    ax & 0 \le x \le s \\ 
    bx+s(a-b) & s \le x \\
\end{cases} &
\end{align}    
where $a,b$ and $s$ are elements of $A$ such that $s\geq0$ and $a > b>1$. Fix $a,b$ and $s$ and let $\alpha=\alpha(a,b,s)$. Then, since every element of $G$ is order-preserving, we have that $\lim_{n\to +\infty}\alpha^n(s)=+\infty$ and that $\lim_{n\to +\infty}\alpha^{-n}(s)=0$, so $\operatorname{orbl}(\alpha,s) = \mathbb{R^{+}}$. Now, let~$G_0$ be the subgroup of all the affine functions of $G$. Since $f(0)=0$ for any element $f$ of $G$, $G_0$ is the subgroup of $G$ consisting of all functions which, if not trivial, have in $0$ the only point of non-differentiability. In particular, the identity is the only element of $G_0$ centralizing any element of $\mathcal{F}$. Moreover, from~Lem\-ma~3.17~(b) of~\cite{Brin}, it follows that any $(0,0)$-conjugator of any function of $\mathcal{F}$ is an element of $G_0$. Take 
$\alpha\in\mathcal{F}$. By Proposition~3.24 of~\cite{Brin} we have that $C_G(\alpha)\subseteq\langle\alpha\rangle \{f\}G_0=\langle\alpha\rangle G_0$, where $f$ is a~\hbox{$(0,0)$-con}\-ju\-ga\-tor of $G$. By the Dedekind identity, this means that $C_G(\alpha)=\langle\alpha\rangle$. Therefore, 
 \[P(G)\le \underset{\alpha\in\mathcal{F}}{\bigcap}\langle\alpha\rangle^G  \le G',\] and hence $P(G)=G'$.
\end{proof}

%


\begin{theorem}\label{PLFR}
Let $F$ be a subfield of $\mathbb{R}$, and set $A=(F,+)$ and $P=(F^+,\cdot)$. Then the pseudocentre of $G=G(\mathbb{R};A,P)$ is equal to its commutator subgroup.
\end{theorem}
\begin{proof}
The commutator subgroup of $G=G(\mathbb R;A,P)$ consists of the elements of $G$ which have initial and final slopes equal to $1$, while the second commutator is the set of all elements of $G$ which are the identity in a neighborhood of $\pm\infty$ and is a non-abelian simple group by Lemma \ref{B}. In particular, $G''\leq P(G)$.

We claim that $P(G)\leq G'$. Consider the following element of $G$:\begin{align}
        f(x) = \begin{cases}
          x-4  & x\le -1 \\
         3x-2   & -1 \le x \le 1 \\
         5x-4 & 1 \le x \le 2 \\
         x+4 & 2 \le x
        \end{cases}
    \end{align}
and its inverse \begin{align}
        f^{-1}(x) = \begin{cases}
          x+4  & x\le -5 \\
         x/3+2/3   & -5 \le x \le 1 \\
         x/5+4/5 & 1 \le x \le 6 \\
         x-4 & 6 \le x
        \end{cases}
    \end{align}

It is easy to see that $J_1=\operatorname{orbl}(f,2)=(1,+\infty)$, $J_2=\operatorname{orbl}(f,0)=(-\infty,1)$ and $f(1)=1$. Now, define \[
        f_{J_1}(x) = \begin{cases}
          f(x)  & x\in J_1 \\
         x   & x\not\in J_1 \\
        \end{cases}\quad\textnormal{and}\quad 
        f_{J_2}(x) = \begin{cases}
          f(x)  & x\in J_2 \\
         x   & x\not\in J_2 \\
        \end{cases}
\] and let $g_1$ and $g_2$ be the restrictions of $f_{J_1}$ and $f_{J_2}$ to $J_1$ and $J_2$, respectively. It follows from Theorem 5.5 of \cite{Brin} that $C_G(f)\simeq C_{\operatorname{PLF}(J_1)}(g_1)\times C_{\operatorname{PLF}(J_2)}(g_2)$. We need to show that $C_{\operatorname{PLF}(J_i)}(g_i)=\langle g_i\rangle$ for~\hbox{$i=1,2$.} We prove it for $i=1$. To this aim, let $G_0$ be the subgroup of all elements of $\operatorname{PLF}(J_1)$ with no break points. Then~$G_0$ is the set of all functions of $J_1$ of the form $ax+(1-a)$. By Theorem 3.24 of \cite{Brin}, $C_{\operatorname{PLF}(J_1)}(g_1)\leq \langle g_1\rangle G_0$, so $$C_{\operatorname{PLF}(J_1)}(g_1)\cap \langle g_1\rangle G_0=\langle g_1\rangle\big(C_{\operatorname{PLF}(J_1)}(g_1)\cap G_0\big)=\langle g_1\rangle$$ because $C_{\operatorname{PLF}(J_1)}(g_1)\cap G_0=\{1\}$ (one just needs to check the commutativity relation in the interval $x\geq2$). Similarly, we deal with the case $i=2$. Thus, $C_G(f)=\langle f_{J_1}\rangle\times\langle f_{J_2}\rangle$. Since $G'$ consists of all elements of $G$ which have initial and final slopes equal to $1$, it follows that $C_G(f)\leq G'$, so also $P(G)\leq G'$, and the claim is proved.

Recall that the intercepts of the linear pieces of any element of $G$ are always contained in $A$, since, by definition, every element of $G(\mathbb{R};A,P)$ has to map $A$ onto itself. We claim that $\langle g\rangle^{G} G'' = G'$. To this aim, observe first that $G'$ is generated by the functions with positive initial and final intercepts. Now, let $\vartheta : G' \to A\times A$ be the epimorphism mapping every element $g$ of~$G'$ to the pair of the  intercepts of the initial and final linear pieces of $g$, and let $g$ and $g_1$ be elements of~$G'\setminus G''$ and $G'$ such that $\vartheta(g)=(u,v)$ and $\vartheta(g_1)=(u_1,v_1)$, where $u,v,u_1$ and $v_1$ are positive. If we choose $h\in G$ with initial linear piece $u^{-1}u_1x$ and final linear piece $v^{-1}v_1x$, then $g^h (x) = x + v_1$ for any large enough $x$, and $g^h(x)=x+u_1$ for any small enough~$x$. This means that $\vartheta(g^h)=\vartheta(g_1)$, so  $g^hg_1^{-1}\in\operatorname{Ker}(\vartheta)=G''$. The arbitrariness of $g_1$ shows that~\hbox{$\langle g\rangle^{G} G'' = G'$.} Therefore, $G'/G''$ is a minimal normal subgroup of~$G/G''$.

Let $g\in G$. If $g$ is a non-trivial element of $G''$ which is the identity function from a point~$\overline{x}$ on, then $g$ is centralized by any element of $G'\setminus G''$ which is the identity before $\overline{x}$ and this shows that $G'\leq C_G(g)^G$, since $G''$ is simple and $G'/G''$ is a minimal normal subgroup of $G/G''$. If $g\in G'\setminus G''$, then there is $h\in G'$ such that $[g,h]$ is a non-trivial element of $G''$, and hence $G'\leq\langle g\rangle^G\leq C_G(g)^G$. Let $g$ be an element of $G\setminus G'$. Then either its initial or final slopes are not both $1$. Without loss of generality, we assume the final slope is not $1$, so the final linear piece of $g$ can be written as  $ax+b$, where $a,b\in P$ with $a\neq1$. Choose~$h$ in $G$ with final linear piece $x+c$ with $c\in A\setminus\{0\}$. Then the final linear piece of $[g,h]$ is $x+c(1-a)$,
so $[g,h]\in (G'\setminus G'')\cap C_G(g)^G$, and, as before, also in this case $G'$ is contained in~$C_G(g)^G$. This shows that $G'\leq P(G)$, and consequently that $G'=P(G)$.
\end{proof}

\begin{theorem}\label{PLF0l/2}
Let $F$ be a subfield of $\mathbb{R}$, and set $A=(F,+)$ and $P=(F^+,\cdot)$. Let $a<b$ be real numbers with~$a$ or $b$ not contained in $A$. Then the pseudocentre of $G=G((a,b);A,P)$ is its commutator subgroup.
\end{theorem}
\begin{proof}
If both $a$ and $b$ are not contained in $A$, then $G$ equals $B((a,b);A,P)$, which is a non-abelian simple group, so we may assume that precisely one between~$a$ and $b$ is not contained in $A$. If $b$ is not contained in $A$, then we can transform $G$ via the affine mapping $x\in\mathbb{R}\mapsto -x+a+b\in\mathbb{R}$ and suppose that $b$ is in $A$ and $a$ is not, which means that we may assume $G=\operatorname{Ker}(\varphi)$, where $$\varphi: G((0,l),A,P) \to P$$ is the homomorphism mapping $g \in G((0,l),A,P)$ to the initial slope of $g$, and $l\in A$ (see Proposition E16.6 in \cite{Bier}).

The proof is similar to that of Theorem \ref{PLF0l}, so we are going to give only a sketch. As in the previous case, $G'$ is contained in $P(G)$ and for every $s\in(0,l)\cap A$ and $1<r\in A$ we can define $\gamma(s,r)$ as the function which is the identity before $l/2$, and the suitably shrunk and shifted version of $\gamma(s,r)$ defined in Theorem \ref{PLF0l} from $l/2$ on. Fix $s$ and $r$ and call $f$ the first linear piece of $\gamma(s,r)$ and $g$ the remaining part. Moreover, let $H$ and $K$ be the subgroups of $G$ consisting of all the functions which are the identity after and before $l/2$, respectively. Now, Theorem 5.5 of \cite{Brin} gives that $C_G(\gamma(s,r))=C_H(f)\times C_K(g)=H\times C_K(g)$. Clearly, $C_H(f)=H\leq G'$, while just as in Theorem \ref{PLF0l} it is not difficult to see that $\operatorname{orbl}(g,\frac{s+l}{2})=(l/2,l)$. Then, choosing carefully $s$ and $a$, we may find $\gamma_1$ and $\gamma_2$ such that $\langle\gamma_1 \rangle^{G}\cap \langle \gamma_2 \rangle^{G} \le G'$ and hence $P(G) \le G'$.
\end{proof}

\begin{theorem}\label{PLFR+/2}
Let $F$ be a subfield of $\mathbb{R}$, and set $A=(F,+)$ and $P=(F^+,\cdot)$. Let $a$ be a real number not contained in $A$ and let $I$ be either $(-\infty,a)$ or $(a,+\infty)$. Then the pseudocentre of $G=G((a,+\infty)),A,P)$ is its commutator subgroup.
\end{theorem}
\begin{proof}
If $I=(-\infty,a)$, we can again transform $G$ via the affine mapping $x\in\mathbb{R}\mapsto -x+a\in\mathbb{R}$ and assume $G=\operatorname{Ker}(\varphi)$, where $\varphi: G((0,+\infty),A,P) \to P$ is the homomorphism mapping $g \in G((0,l),A,P)$ to the initial slope of $g$ (see~Pro\-po\-sition~E16.6 in \cite{Bier}). Just as in Theorem \ref{PLFR+}, we have that~$G''$ is contained in $P(G)$, that~$G'/G''$ is a minimal normal subgroup of~$G/G''$ and that $G'\leq P(G)$, because one can still choose functions final slopes and intercepts in the same way. 

Let now $a,b$ and $s$ be elements of $A$ such that $s\geq1$, $a>b>1$ and define the set $\mathcal F$ of all functions $\alpha(a,b,s)(x)$ which are the identity before $1$, and the suitably shrunk and shifted version of $\alpha(a,b,s)$ defined in Theorem \ref{PLFR+} from $1$ on. Fix~$a,b$ and $s$ and call $f$ the first linear piece of~$\alpha(a,b,s)$ and $g$ the remaining part. Moreover, let $H$ and $K$ be the subgroups of $G$ consisting of all the functions which are the identity after and before $1$, respectively. Now,~The\-o\-rem~5.5 of~\cite{Brin} gives that $$C_G(\alpha(a,b,s))=C_H(f)\times C_K(g)=H\times C_K(g).$$ Clearly, $C_H(f)=H\leq G'$, while just as in Theorem \ref{PLFR+} it is not difficult to see that $\operatorname{orbl}(g,s)=(1,+\infty)$ and that, for instance, the group of all the affine functions of~$G$ is trivial, since they need to start being the identity. Then, we again have \[P(G)\le \underset{\alpha\in\mathcal{F}}{\bigcap}\langle\alpha\rangle^G  \le G',\] and hence $P(G)=G'$.
\end{proof}

\medskip

Let $\{S_i\}_{i\in I}$ be a family of disjoint intervals of $\mathbb R$ and set $S=\bigcup_{i\in I} S_i$. Let~$P$ be a subgroup of the multiplicative group of the positive real numbers, and let~$A$ be a $\mathbb{Z}[P]$-submodule of $(\mathbb{R},+)$. In order to study $G(S;A,P)$, one can clearly assume that $\operatorname{cl}(S_i)\cap\operatorname{cl}(S_j)=\emptyset$ for any $i\neq j\in I$. Then, we have that $$G(S;A,P)=\operatorname{Dr}_{i\in I} G(S_i^\circ;A,P).$$ Combining this observation and the last five theorems, we obtain the main theorem of this chapter.

\begin{theorem}\label{finaltheoremthompson}
Let $\{S_i\}_{i\in I}$ be a family of intervals of $\mathbb R$ with $\operatorname{cl}(S_i)\cap\operatorname{cl}(S_j)=\emptyset$ for any $i\neq j\in I$, $F$ a subfield of $\mathbb{R}$, and set $A=(F,+)$, $P=(F^+,\cdot)$ and $S=\bigcup_{i\in I} S_i$. Then the pseudocentre of~$G(S;A,P)$ is its commutator subgroup.
\end{theorem}

\part{Further topics}

\chapter{Groups whose proper subgroups are pseudocentral}\label{propersubgroupsarepseudocentral}

Since the class of pseudocentral groups is not closed with respect to forming subgroups, an interest may arise in understanding what a group whose all (proper) subgroups are pseudocentral would look like. This is a very natural question for properties that are not inherited by subgroups, as for example happens for the class of groups in which normality is a transitive relation (see \cite{Ro64} and \cite{Ro69}).  As one may expect, it turns out that these kinds of groups mostly resemble minimal non-abelian groups. In fact, since every  locally graded minimal non-abelian group is metabelian, then every locally graded minimal non-abelian group is a minimal non-pseudocentral group. The following result shows that the converse holds in the universe of periodic groups.

\begin{theorem}\label{minimalnonabelian}
Let $G$ be a locally graded periodic group whose proper subgroups are pseudocentral. Then $G$ is metabelian. In particular, $G$ is either abelian or minimal non-abelian.
\end{theorem}
\begin{proof}
Suppose first that $G$ is finite. We use induction on the order of $G$ to prove that~$G$ is metabelian. If $H$ is a proper subgroup of $G$, then $H$ is metabelian because all proper subgroups of $H$ are pseudocentral. But then $H$ is abelian by Corollary \ref{pseudodercentr}. Then every proper subgroup of~$G$ is abelian, and hence~$G$ is either abelian or minimal non-abelian. In any case,~$G$ is metabelian. 

Assume that $G$ is arbitrary periodic. Since the class of metabelian groups is local, we only need to show that $G$ is locally finite. Let $H$ be an infinite finitely generated  subgroup of $G$. Then~$H/H''$ is finite, so there is a proper $H$-invariant subgroup $K$ of finite index of $H''$. However, the finite case yields that~$H/K$ is metabelian and hence that $H''\leq K$, a contradiction.~\end{proof}

\begin{quest}
Can we get rid of the periodicity assumption in the statement of \textnormal{Theorem \ref{minimalnonabelian}}?
\end{quest}

We are not able to answer this question even if the whole group is pseudocentral. In fact, as shown by the following consequences of Theorem \ref{minimalnonabelian}, answering the previous question could be a difficult task.

\begin{corollary}
Let $G$ be a locally graded group with an ascending series whose factors are locally \textnormal(soluble-by-periodic\textnormal). If all proper subgroups of $G$ are pseudocentral, then $G$ is metabelian. In particular,~$G$ is either abelian or minimal non-abelian.
\end{corollary}
\begin{proof}
By the main result of \cite{locallygraded}, we may assume that $G$ is locally (soluble-by-periodic). As in the proof of~The\-o\-rem~\ref{minimalnonabelian}, we only need to show that $G$ is soluble. If $G$ is finitely generated, then $G$ is soluble-by-periodic, and hence soluble by Theorem \ref{minimalnonabelian} (and again \cite{locallygraded}). If $G$ is not finitely generated, then every finitely generated subgroup of $G$ is soluble pseudocentral, and so even abelian by Co\-rol\-la\-ry~\ref{pseudodercentr}. Consequently, $G$ is abelian.
\end{proof}

\begin{corollary}
Let $G$ be a locally graded group whose commutator subgroup is finitely generated. If all proper subgroups of $G$ are pseudocentral, then $G$ is metabelian. In particular,~$G$ is either abelian or minimal non-abelian.
\end{corollary}
\begin{proof}
Suppose $G$ is not abelian. Since $G'$ is finitely generated and $G$ is locally graded, so there is a proper normal subgroup $N$ of $G'$ such that $G'/N$ is finite. By Theorem \ref{minimalnonabelian}, $G'/N$ is soluble, and hence $G''<G'$. Since all proper subgroups of $G/G''$ are pseudocentral, they are abelian, and hence $G/G''$ is minimal non-abelian. Thus, $G/G''$ is finite, and $G''$ is finitley generated.  Since $G'$ is a proper pseudocentral subgroup of $G$, so $G''=G'''$ by Corollary \ref{pseudodercentr}, and hence $G''=\{1\}$ by Theorem \ref{minimalnonabelian}.
\end{proof}

\medskip

 The consideration of a {\it Tarski group} (that is, an infinite group whose proper non-trivial subgroups are cyclic of prime order) shows that there exist periodic (but not locally graded) minimal non-abelian groups that are not minimal non-pseudocentral. Thus, the following question is a natural one.

\begin{quest}
Is there a minimal non-pseudocentral group that is not minimal non-abelian?
\end{quest}

\chapter{Pseudo-integrability}\label{secpseudoint}

Every time a new characteristic subgroup of a group is introduced, the question arises as to which groups can or cannot appear as this subgroup. A simple example in this context is the centre: every abelian group is its own centre, so every abelian group can be the centre of a group. The situation dramatically changes when we consider the derived subgroup: which groups can be derived subgroups of some groups? This problem first appeared in \cite{neumann11}, where it was also noted that not all groups can be derived subgroups, and has been investigated (mostly for finite groups) in~\cite{casolo1}, \cite{casolo2}, and \cite{guralnick}. It has been shown for example that every finite abelian group is {\it integrable} (that is, can be a derived subgroup of a group), and that every integrable group has a finite {\it integral} (that is, a group realizing the integrability).

The goal of this chapter is to provide a brief insight into the analogous problem concerning the pseudocentre of a group, which seems to be even more difficult to tackle due to the further complication that the pseudocentre contains the centre. Thus, we give the following definition. A group $G$ is said to be \textit{pseudo-integrable} if there exists a group $H$ such that $P(H)\simeq G$. Clearly, every pseudocentral group is pseudo-integrable, so in particular every abelian group is pseudo-integrable. As the reader will see, it appears from our results that if a group is pseudo-integrable, then it is also integrable, although the converse does not hold (see Example \ref{exnonfunziona}), and that most results for pseudo-integrability are analogous to the results for integrability. This justifies a bit the name pseudo-integrability, and naturally leads to the following question.

\begin{quest}
Is there any relation between integrable and pseudo-integrable groups?
\end{quest}

It has been proved in \cite{casolo1}, Theorem 5.1, that a complete group is integrable if and only if it is perfect. Recall that a group $G$ is {\it complete} if $\operatorname{Aut}(G)=\operatorname{Inn}(G)$ and $Z(G)=\{1\}$. Our first result gives a general criterion for a complete group to be pseudo-integrable. This criterion will be employed in Example \ref{exnonfunziona} to show that there exist integrable groups that are not pseudo-integrable.

\begin{theorem}\label{thpseudointegra}
    Let $G$ be a complete group. Then $G$ is pseudocentral if and only if $G$ is pseudo-integrable.
\end{theorem}
\begin{proof}
Clearly, if $G$ is pseudocentral, then it is also pseudo-integrable. Assume $P(H)=G$ for some group $H$. By Theorem 13.5.8 of \cite{acourse}, \hbox{$H = L \times P(H)$} for some subgroup $L$ of $H$. Then Lemma $\ref{directproducts}$ yields that $P(H)$ is the direct product of the pseudocentres of the direct factors, so $P(L)=\{1\}$ and $P(P(H))=P(H)$, as required.
\end{proof}

\begin{corollary}
$\operatorname{Sym}(n)$ is not pseudo-integrable for every finite $3\leq n\neq 6$.
\end{corollary}
\begin{proof}
This follows at once from Theorems \ref{thsymalg} and \ref{thpseudointegra}. 
\end{proof}

\medskip

We briefly deal with the case of the symmetric group of degree $6$ as follows.

\begin{theorem}\label{thauths6}
Let $G=\operatorname{Aut}\big(\operatorname{Sym}(6)\big)$. Then $P(G)=\operatorname{Alt}(6)$.
\end{theorem}
\begin{proof}
Let $x$ be the outer automorphism of $\operatorname{Sym}(6)$ defined by \[
\begin{array}{c}
(12)\mapsto (15)(23)(46),\qquad
(13)\mapsto (14)(26)(35),\qquad
(14)\mapsto (13)(24)(56),\\[0.2cm]
(15)\mapsto(12)(36)(45),\qquad
(16)\mapsto(16)(25)(34).
\end{array}
\] Then $G=\langle x\rangle\ltimes\operatorname{Sym}(6)$ and $o(x)=2$. Now, $\operatorname{Alt}(6)$ is simple, so The\-o\-rem~\ref{normsempl} yields that $P(G)\geq\operatorname{Alt}(6)$. Moreover, it is easy to see that $(12345)^x=(12345)$, so~$C_G((12345))^G\leq\langle x\rangle\operatorname{Alt}(6)$. On the other hand, $C_G((123))^G\leq \operatorname{Sym}(6)$. Consequently $P(G)=\operatorname{Alt}(6)$.
\end{proof}

\begin{corollary}
$\operatorname{Sym}(6)$ is not pseudo-integrable.
\end{corollary}
\begin{proof}
Suppose by way of contradiction that there is a group $G$ such that $P=P(G)\simeq\operatorname{Sym}(6)$. Since $Z(P)=\{1\}$, so $Z(G)=\{1\}$. Moreover, $G/C_G(P)$ is either isomorphic to $\operatorname{Sym}(6)$ or to $\operatorname{Aut}(\operatorname{Sym}(6))$. In both cases, we have a contradiction  by Theorems \ref{thsymalg} and \ref{thauths6}.
\end{proof}

\medskip

Further applications of Theorem \ref{thpseudointegra} can be obtained by combining our next lemma and Theorem \ref{affinegeneral}.

\begin{lemma}\label{complete}
   Let $2\leq n\neq3$ be an integer and set $G =\operatorname{GL}(n,2) \ltimes V$, where~$V$ is the vector space of dimension $n$ over the field of order $2$. Then $\operatorname{Aut}(G) \simeq G$. 
\end{lemma}
\begin{proof}
The automorphism group of $\operatorname{GL}(n,2)$ is defined by the short exact sequence $$1 \to \operatorname{CAut}(\operatorname{GL}(n,2)) \to \operatorname{Aut}(\operatorname{GL}(n,2)) \to \operatorname{Aut}(\operatorname{PGL}(n,2)) \to 1 $$ where $\operatorname{CAut}(\operatorname{GL}(n,2))$ is the central automorphism group. Now, $$\operatorname{Aut}(\operatorname{PGL}(n,2)) = \mathbb Z_2 \ltimes\operatorname{PGL}(n,2),$$ where the cyclic group of order $2$ is the transpose-inverse map. However, in our case, $$\operatorname{PGL(n,2)} = \operatorname{GL(n,2)}$$  and the transpose-inverse map is an inner automorphism. Also, in our case, $$\operatorname{CAut}(\operatorname{GL}(n,2))=\{1\},$$ so $\operatorname{GL}(n,2)\simeq\operatorname{Aut}(\operatorname{GL}(n,2))$ is complete.


Since $V$ is characteristic in $G$, so every automorphism $\alpha$ of $G$ induces an automorphism $\alpha^\varphi$ of $G/V$. Then the map $$\varphi:\operatorname{Aut}(G)\rightarrow\operatorname{Aut}(G/V)\simeq \operatorname{GL}(n,2)$$ is a well-defined homomorphism. Set $K = \operatorname{Ker}(\varphi)$, and note that the set $L$ of all inner automorphisms determined by elements of $V$ is a subset of $K$. In order to prove that $\operatorname{Aut}(G)\simeq G$, it remains to show that $K/L$ is trivial. Thus, we need to consider the automorphisms of $G$ acting trivially on $V$ and $G/V$ modulo the inner automorphisms determined by $V$. These automorphisms correspond to the first cohomology group $\operatorname{H^1(\operatorname{GL(n,2)},V)}$. By Table 1 of \cite{Bell}, we have that this cohomology group is trivial, so that $\operatorname{Aut}(G) = G$.
\end{proof}

\begin{example}
The natural semidirect product $\operatorname{SL}(5,2)\ltimes \mathbb Z_2^5$ is complete and pseudocentral, so is pseudo-integrable.
\end{example}
\begin{proof}
This follows from Lemma \ref{complete} and Theorem \ref{affinegeneral}.
\end{proof}

\medskip

\begin{example}\label{exnonfunziona}
There exists an integrable group that is not pseudo-integrable.
\end{example}
\begin{proof}
Let $H=\langle x,y \rangle$ be the subgroup of $\operatorname{SL}(10,2)$ generated by the matrices $$x= \begin{pmatrix}
      0 & 0 & 1 & 1 & 0 & 0 & 1 & 0 & 0 & 1\\
      0 & 0 & 0 & 0 & 1 & 0 & 1 & 0 & 1 & 0\\
      1 & 1 & 0 & 1 & 1 & 0 & 0 & 0 & 1 & 0\\
      0 & 1 & 0 & 1 & 1 & 0 & 0 & 1 & 0 & 1\\
      0 & 1 & 0 & 0 & 0 & 0 & 1 & 1 & 1 & 0\\
      0 & 1 & 1 & 1 & 1 & 1 & 1 & 0 & 1 & 0\\
      0 & 0 & 0 & 0 & 0 & 0 & 1 & 0 & 0 & 0\\
      0 & 0 & 0 & 0 & 0 & 0 & 0 & 1 & 0 & 0\\
      0 & 0 & 0 & 0 & 0 & 0 & 1 & 0 & 1 & 0\\
      0 & 0 & 0 & 0 & 0 & 0 & 0 & 1 & 0 & 1
\end{pmatrix}$$ and $$y=\begin{pmatrix}
      1 & 0 & 0 & 0 & 1 & 1 & 1 & 1 & 0 & 1\\
      0 & 0 & 1 & 1 & 0 & 0 & 0 & 0 & 1 & 1\\
      0 & 1 & 1 & 1 & 0 & 0 & 0 & 1 & 1 & 1\\
      0 & 0 & 0 & 1 & 0 & 0 & 0 & 0 & 1 & 0\\
      0 & 1 & 0 & 1 & 1 & 1 & 1 & 1 & 0 & 1\\
      0 & 0 & 1 & 1 & 1 & 0 & 1 & 1 & 0 & 1\\
      0 & 0 & 0 & 0 & 0 & 1 & 1 & 0 & 1 & 0\\
      0 & 0 & 0 & 0 & 0 & 1 & 0 & 0 & 0 & 0\\
      0 & 0 & 0 & 0 & 0 & 1 & 1 & 1 & 1 & 0\\
      0 & 0 & 0 & 0 & 0 & 0 & 0 & 0 & 0 & 1 
\end{pmatrix}.$$ Consider the natural action of $H$ on the vector space $V$ of dimension $10$ over the field of order $2$, and put $G=H\ltimes V$. It is possible to check on GAP that $G$ is perfect and complete, so is integrable by \cite{casolo1}, Theorem 5.1. Similarly, we can check that $G$ is not pseudocentral (see also Remark \ref{cuborubik}), so is not pseudo-integrable by~The\-o\-rem~\ref{thpseudointegra}. Note that this was computed using the ‘‘Perfect Groups’’ library and that the Id of the group is (344064,33).
\end{proof}


The following easy observation shows that no group with a non-central locally cyclic characteristic subgroup can be pseudo-integrable, so in particular, no group with a non-central locally cyclic derived subgroup can be pseudo-integrable.

\begin{lemma}\label{cycliccentral}
Let $G$ be a group and $H$ a normal subgroup such that $G/C_G(H)$ is Dedekind. Then $H\cap P(G)\leq Z(P(G))$. 
\end{lemma}
\begin{proof}
Let $h\in H\cap P(G)$. Then $C_G(H)\leq C_G(h)=C_G(h)^G$. Consequently,~\hbox{$P(G)\leq C_G(h)$} and we are done.
\end{proof}

\medskip

Thus, for example the following groups are not pseudo-integrable:
\begin{itemize}
    \item Every finite non-abelian dihedral group.
    \item The infinite dihedral group. 
    \item The locally dihedral $2$-group.
    \item The generalized quaternion groups $Q_{2^n}$ for $n\geq4$.
    \item The quasi-dihedral groups $QD_{2^n}$ for $n\geq4$.
    \item Every non-abelian group of order $pq$, where $p$ and $q$ are primes.
\end{itemize}

\begin{quest}
Let $G$ be a pseudo-integrable finite group. Is it possible to find a finite group $H$ such that $G\simeq P(H)$?
\end{quest}

\bigskip\bigskip

\begin{flushleft}
\rule{8cm}{0.4pt}\\
\end{flushleft}

{
\sloppy
\noindent
Bernardo Giuseppe Di Siena

\noindent
Dipartimento di Matematica e Fisica

\noindent
Università degli Studi della Campania  ``Luigi Vanvitelli''

\noindent
viale Lincoln 5, Caserta (Italy)

\noindent
e-mail: bernardogiuseppe.disiena@unicampania.it
}

\bigskip
\bigskip

{
\sloppy
\noindent
Mattia Brescia, Ernesto Ingrosso, Marco Trombetti

\noindent 
Dipartimento di Matematica e Applicazioni ``Renato Caccioppoli''

\noindent
Università di Napoli Federico II

\noindent
Complesso Universitario Monte S. Angelo

\noindent
Via Cintia, Napoli (Italy)

\noindent
e-mail: mattia.brescia@unina.it; ernesto.ingrosso2@unina.it; marco.trombetti@unina.it 

}

\end{document}